\documentclass[leqno,11pt]{amsart}

\usepackage[dvipsnames]{xcolor}

\usepackage{multirow,booktabs}

\usepackage{amsmath, amssymb, amsfonts, latexsym, mdwlist, amsthm}
\usepackage{subfig}
\usepackage{graphicx}
\usepackage{tikz-cd}
\usepackage{wrapfig}

\usepackage{comment}
\usepackage{mathtools}
\usepackage{tikz}
\usetikzlibrary{calc,trees,positioning,arrows,chains,shapes.geometric,%
	decorations.pathreplacing,decorations.pathmorphing,shapes,%
	matrix,shapes.symbols}

\tikzset{
	>=stealth',
	punktchain/.style={
		rectangle,
		rounded corners,
		draw=black, thick,
		minimum height=3em,
		text centered,
		on chain},
	line/.style={draw, thick, <-},
	element/.style={
		tape,
		top color=white,
		bottom color=blue!50!black!60!,
		minimum width=8em,
		draw=blue!40!black!90, very thick,
		text width=10em,
		minimum height=3.5em,
		text centered,
		on chain},
	every join/.style={->, thick,shorten >=1pt},
	decoration={brace},
	tuborg/.style={decorate},
	tubnode/.style={midway, right=2pt},
}

\usepackage[all]{xy} 

\usetikzlibrary{patterns}

\usepackage{enumerate}

\usepackage{amsmath, amssymb, amsfonts, latexsym, mdwlist, amsthm}
\usepackage{subfig}
\usepackage{graphicx}
\usepackage{wrapfig}

\usepackage{mathtools}
\usepackage{bm}

\definecolor{lightred}{HTML}{ff4d4d}
\definecolor{lightblue}{HTML}{1F88CD}
\definecolor{lightgrey}{HTML}{727272}
\definecolor{lightblue2}{HTML}{009EC1}
\definecolor{mypink}{HTML}{FD00B0}

\usepackage[bookmarks, colorlinks, breaklinks, pdftitle={},
pdfauthor={Augustinas Jacovskis, Zhiyu Liu and Shizhuo Zhang}]{hyperref}
\hypersetup{linkcolor=OliveGreen,citecolor=lightred,filecolor=black,urlcolor=blue}

\usepackage{tikz}
\usetikzlibrary{calc,trees,positioning,arrows,chains,shapes.geometric,%
    decorations.pathreplacing,decorations.pathmorphing,shapes,%
    matrix,shapes.symbols}

\tikzset{
>=stealth',
  punktchain/.style={
    rectangle,
    rounded corners,
    draw=black, thick,
    minimum height=3em,
    text centered,
    on chain},
  line/.style={draw, thick, <-},
  element/.style={
    tape,
    top color=white,
    bottom color=blue!50!black!60!,
    minimum width=8em,
    draw=blue!40!black!90, very thick,
    text width=10em,
    minimum height=3.5em,
    text centered,
    on chain},
  every join/.style={->, thick,shorten >=1pt},
  decoration={brace},
  tuborg/.style={decorate},
  tubnode/.style={midway, right=2pt},
}

\usepackage{paralist}
\setdefaultenum{(a)}{(i)}{}{}
\usepackage[shortlabels]{enumitem} 

\newcommand{\bv}{\mathbf{v}}
\newcommand{\bw}{\mathbf{w}}

\makeatletter
\newtheorem*{rep@theorem}{\rep@title}
\newcommand{\newreptheorem}[2]{%
\newenvironment{rep#1}[1]{%
 \def\rep@title{#2 \ref{##1}}%
 \begin{rep@theorem}}%
 {\end{rep@theorem}}}
\makeatother

\newtheorem{Thm}{Theorem}[section]
\newreptheorem{Thm}{Theorem}
\newtheorem{Prop}[Thm]{Proposition}

\newtheorem{Lem}[Thm]{Lemma}

\newtheorem{Cor}[Thm]{Corollary}
\newreptheorem{Cor}{Corollary}

\newreptheorem{Con}{Conjecture}

\newtheorem{thm-int}{Theorem}

\theoremstyle{definition}
\newtheorem{Def-s}[Thm]{Definition}
\newtheorem{Def}[Thm]{Definition}
\newtheorem{Rem}[Thm]{Remark}

\newtheorem{Ex}[Thm]{Example}

\newcommand{\ignore}[1]{}

\usepackage[numbers]{natbib}
\newcommand{\ra}{\rightarrow}
\newcommand{\xra}{\xrightarrow}

\newcommand{\sst}{\subset}

\DeclareMathOperator{\bL}{\bm{\mathrm{L}}}
\DeclareMathOperator{\bR}{\bm{\mathrm{R}}}

\newcommand{\ZZ}{\mathbb{Z}}

\newcommand{\CC}{\mathbb{C}}
\newcommand{\PP}{\mathbb{P}}

\newcommand{\Db}{\mathrm{D^b}}

\newcommand{\ch}{\mathrm{ch}}

\newcommand{\pr}{\mathrm{pr}}

\renewcommand{\Re}{\operatorname{Re}}
\renewcommand{\Im}{\operatorname{Im}}

\DeclareMathOperator{\Aut}{Aut}

\DeclareMathOperator{\identity}{id}
\DeclareMathOperator{\im}{im}

\DeclareMathOperator{\rk}{rk}
\DeclareMathOperator{\cok}{cok}

\DeclareMathOperator{\Coh}{\mathsf{Coh}}

\DeclareMathOperator{\Ext}{Ext}

\DeclareMathOperator{\Hom}{Hom}
\DeclareMathOperator{\RHom}{RHom}
\DeclareMathOperator{\ext}{ext}

\DeclareMathOperator{\Pic}{Pic}

\DeclareMathOperator{\cone}{cone}

\DeclareMathOperator{\Gr}{Gr}

\newcommand{\cC}{\mathcal{C}}
\newcommand{\cA}{\mathcal{A}}
\newcommand{\cE}{\mathcal{E}}
\newcommand{\cU}{\mathcal{U}}
\newcommand{\cH}{\mathcal{H}}
\newcommand{\cB}{\mathcal{B}}
\newcommand{\cK}{\mathcal{K}}
\newcommand{\cI}{\mathcal{I}}
\newcommand{\cT}{\mathcal{T}}
\newcommand{\cQ}{\mathcal{Q}}
\newcommand{\Ku}{\mathcal{K}u}
\newcommand{\cP}{\mathcal{P}}
\newcommand{\cD}{\mathcal{D}}

\newcommand{\cN}{\mathsf{N}}

\DeclareMathOperator{\cF}{\mathcal{F}}
\DeclareMathOperator{\cG}{\mathcal{G}}
\DeclareMathOperator{\E}{\mathcal{E}}
\DeclareMathOperator{\oh}{\mathcal{O}}
\numberwithin{equation}{section}

\usepackage{hyperref}
\hypersetup{
	colorlinks=true,
    linkcolor={blue},
    citecolor={blue},
	urlcolor={black}
}

\begin{document}

\title[BN reconstruction of index one prime Fano threefolds]{Brill--Noether reconstruction of index one \\ prime Fano threefolds}


\author{Augustinas Jacovskis, Zhiyu Liu and Shizhuo Zhang}

\address{Room B434, National Center for Physical and Technological Sciences, Sauletekio al. 3, Vilnius, 10257, Lithuania}
\email{augustinas.jacovskis@ff.vu.lt}

\address{School of Mathematical Sciences, Zhejiang University, Hangzhou, Zhejiang Province 310058, P. R. China}
\email{jasonlzy0617@gmail.com}

\address{School of Mathematics, Sun Yat-sen University, Guangzhou 510275, P. R. China}
\email{zhangshzh28@mail.sysu.edu.cn}

\begin{abstract}
We show by a uniform argument that every index one prime Fano threefold $X$ of genus $g\geq 6$ can be reconstructed as a Brill--Noether locus inside a Bridgeland moduli space of stable objects in the Kuznetsov component $\Ku(X)$. As an application, we verify Mukai's conjecture on the existence of dual embeddings of $X$. Moreover, we establish a refined categorical Torelli theorem for $X$ and classify autoequivalences of $\Ku(X)$. We also give an alternative disproof of Kuznetsov's Fano threefold conjecture. 

\end{abstract}

\vspace{-1em}
\maketitle


\section{Introduction}

In the landmark paper \cite{bondal2001reconstruction}, Bondal and Orlov showed that the bounded derived category of coherent sheaves $\Db(X)$ on a smooth projective variety $X$ determines the isomorphism class of $X$, provided that either the anticanonical divisor $-K_X$ or the canonical divisor $K_X$ is ample.

For smooth Fano threefolds, extensive work on semiorthogonal decompositions of $\Db(X)$ has made it possible to recover geometric information about $X$ from individual components of its derived category. This naturally raises the question of whether a smooth Fano threefold $X$ can be determined, up to isomorphism, by a distinguished component of $\Db(X)$. The most natural candidate is the \emph{Kuznetsov component} $\Ku(X)\subset \Db(X)$. This category has been studied extensively, for instance in \cite{kuznetsov2003derived,kuznetsov2010derivedcubicfourfolds,ku:fractional-cy,Kuznetsov_2016,kuznetsov:fano-threefolds,kuznetsov2018derived,brambilla2011moduli,brambilla2013rank,brambilla2014vector}, for many classes of Fano varieties, and it is widely expected to encode the birational geometry of $X$. For example, it is shown in \cite{bernardara2012categorical,PY20,bayer2020desingularization} that two cubic threefolds are isomorphic if and only if their Kuznetsov components are equivalent. We refer to \cite{pertusi2022categorical} for further related results.

In our previous work \cite{jacovskis2024categorical}, we showed that the Kuznetsov component $\Ku(X)$ of certain special prime Fano threefolds of genus $6$ determines their isomorphism class, while, for a general such threefold, $\Ku(X)$ only determines its birational class. Moreover, in the general case, we proved that $\Ku(X)$ together with certain additional data determines $X$ up to isomorphism. We call such a statement a \emph{refined categorical Torelli theorem}. It is therefore natural to ask whether refined categorical Torelli theorems of this kind hold for all index one prime Fano threefolds.

A complementary source of motivation comes from Mukai's Brill--Noether reconstructions. In \cite{mukai02,mukai:nonabelian}, Mukai showed that every index one prime Fano threefold $X$ of genus $g=7$ or $9$ can be reconstructed as a Brill--Noether locus in the moduli space $M_C(2,d)$ of stable vector bundles of rank $2$ and degree $d$ on a curve $C$, where $C$ has genus $7$ when $g=7$ and genus $3$ when $g=9$. In precisely these two cases, one has $\Ku(X)\simeq \Db(C)$; see \cite{brambilla2013rank,brambilla2014vector}. Thus the moduli spaces $M_C(2,d)$ may be viewed as Bridgeland moduli spaces of stable objects in the corresponding Kuznetsov components.

Moreover, in his ICM report \cite{mukai:icm}, Mukai conjectured the existence of a ``dual'' embedding for every index-one prime Fano threefold $X$ of genus $g\geq 6$, whose normal bundle is the twisted dual of the normal bundle of the natural embedding into a Grassmannian. When $g=7$ and $g=9$, these dual embeddings are precisely the embeddings into $M_C(2,d)$ as Brill--Noether loci. This suggests that a Brill--Noether type reconstruction should hold uniformly for all $g\geq 6$, and should provide the corresponding dual embedding of $X$.

In this paper, we prove the following theorem, which simultaneously generalizes the refined categorical Torelli theorem for genus $6$ threefolds \cite[Theorem 10.2]{jacovskis2024categorical}, Mukai's Brill--Noether reconstruction in genera $7$ and $9$ \cite{mukai02,mukai:nonabelian}, and verifies Mukai's conjecture.

Let $X$ be a prime Fano threefold of index one and genus $g\geq 6$, and consider the semiorthogonal decomposition
\[
\Db(X)=\langle\Ku(X),\cE_X,\oh_X\rangle,
\]
where $\cE_X$ is the pullback of the tautological subbundle on a Grassmannian $\mathbf{G}_X$ associated with $X$; see Section~\ref{sec:derivedfano}. Stability conditions on $\Ku(X)$ were constructed in \cite{bayer2017stability}, and were shown in \cite{Pertusi2021serreinv} to be Serre-invariant; see Definition~\ref{Serre_invariant_stab}.

Let $i^!$ and $i^*$ denote, respectively, the right and left adjoints of $i\colon \Ku(X)\hookrightarrow \Db(X)$. Set
\[
n_g\coloneqq \frac{g}{2}
\]
when $g$ is even, and set $n_7=5$ and $n_9=3$. Our main theorem is the following.

\begin{Thm}[{Theorem~\ref{thm: BN = X} and~\ref{thm:normal-bundle}}]
\label{main_theorem_Brill_Nother_reconstruction}
Let $X$ be a prime Fano threefold of index one and genus $g\geq 6$, and let $\sigma$ be a Serre-invariant stability condition on $\Ku(X)$.

\begin{enumerate}
    \item We have
\[
X\cong \mathsf{BN}_g \coloneqq
\left\{
F\in M_\sigma([i^*\oh_x])
\ \middle|\ 
\dim \Ext^k(F,i^!\cE_X)=n_g
\text{ for some } k\in\ZZ
\right\}.
\]
Here, $\oh_x$ is the skyscraper sheaf supported at a point $x\in X$, and $M_\sigma([i^*\oh_x])$ denotes the moduli space of $\sigma$-stable objects with numerical class $[i^*\oh_x]$.

\item When $g\neq 6$, or $g=6$ and $X$ is general, we have an isomorphism
\[N^{\vee}_{X/\mathbf{G}_X}(H)\cong N_{X/M_\sigma([i^*\oh_x])},\]
where $H$ is the restriction of the Pl\"ucker class on $\mathbf{G}_X$.
\end{enumerate}
\end{Thm}

A priori, this Brill--Noether reconstruction of $X$ depends on the choice of a stability condition $\sigma$ on $\Ku(X)$. However, the uniqueness of Serre-invariant stability conditions on $\Ku(X)$ was proved in \cite{jacovskis2024categorical,FeyzbakhshPertusi2021stab}. Hence, if $E\in \Ku(X)$ is stable with respect to a Serre-invariant stability condition and $\Phi\colon \Ku(X)\simeq \Ku(X')$ is an equivalence of Kuznetsov components, then $\Phi(E)$ is stable with respect to any Serre-invariant stability condition on $\Ku(X')$. As a consequence, we prove refined categorical Torelli theorems for all prime Fano threefolds of index one and genus $g\geq 6$.

\begin{Thm}
\label{Refined_Torelli}
Let $X$ and $X'$ be prime Fano threefolds of index one and genus $g\geq 6$. Suppose that there is an equivalence
\[
\Phi\colon\Ku(X)\simeq \Ku(X') \quad \text{ such that }\quad \Phi(i^!\cE_X)\cong i'{}^!\cE_{X'},
\]
where $i'$ denotes the inclusion $\Ku(X')\hookrightarrow \Db(X')$. Then $X\cong X'$.
\end{Thm}

Theorem~\ref{Refined_Torelli} has an interesting consequence for autoequivalence groups. In \cite{feyzbakhsh2023new}, the authors proved categorical Torelli theorems for del Pezzo threefolds of degree $d\geq 2$ and used them to compute the autoequivalence groups of their Kuznetsov components. Inspired by their method, we use Theorem~\ref{Refined_Torelli} to compute the autoequivalence groups of Kuznetsov components in the genus $6$ and genus $8$ cases.

\begin{Thm}[{Corollary \ref{cor_aut_genus_8} and \ref{cor_aut_genus_6}}]
\label{theorem_group_auto_equi_index_one}
Let $X$ be a prime Fano threefold of index one and genus $g=6$ or $8$. If $g=6$, assume further that $X$ is general. Then
\[
\mathrm{Aut}(\Ku(X))
=
\langle \mathrm{Aut}(X),S_{\Ku(X)},[1]\rangle,
\]
where $S_{\Ku(X)}$ denotes the Serre functor of $\Ku(X)$.
\end{Thm}

Combining this with the corresponding results for del Pezzo threefolds in \cite[Corollary 8.4]{feyzbakhsh2023new}, we obtain the following geometric consequence, which generalizes the classical results in \cite[Corollaries 4.2.3 and 4.3.5]{KPS}.

\begin{Cor}[{Corollary \ref{cor_aut_variety}}]
\label{cor_isom_genus8_cubic}
Let $X$ be a prime Fano threefold of index one and genus $8$, and let $Y$ be the Pfaffian cubic threefold associated with $X$. Then
\[
\mathrm{Aut}(X)\cong \mathrm{Aut}(Y).
\]
\end{Cor}

Kuznetsov conjectured in \cite[Conjecture 3.7]{kuznetsov:fano-threefolds} that certain pairs of prime Fano threefolds should have equivalent Kuznetsov components. For genus $6$ threefolds, this conjecture was disproved in \cite{zhang2020bridgeland,bayer2022kuznetsov}. Using Theorem~\ref{theorem_group_auto_equi_index_one}, we also obtain a new disproof of Kuznetsov's conjecture for Fano threefolds.

\begin{Cor}[{Corollary \ref{cor_ku_conj}}]
Let $X$ be a prime Fano threefold of index one and genus $6$, and let $Y$ be a prime Fano threefold of index two and degree $2$. Then $\Ku(X)$ is not equivalent to $\Ku(Y)$.
\end{Cor}

\subsection{Organization}

In Section~\ref{sec:derived categories}, we recall several useful results on semiorthogonal decompositions. In Section~\ref{sec:derivedfano}, we review semiorthogonal decompositions for index one prime Fano threefolds and define the corresponding \emph{gluing objects} for their Kuznetsov components. In Section~\ref{sec:stability conditions}, we recall weak stability conditions and state the algorithm used throughout the paper to find, and subsequently eliminate, possible numerical walls in tilt stability. We also summarize results on the existence of Bridgeland stability conditions on Kuznetsov components of index one Fano threefolds. In Section~\ref{sec:embedding}, we embed our Fano threefolds into moduli spaces of stable objects in $\Ku(X)$ and show that these are the desired dual embeddings. In Section~\ref{sec:BN}, we realize these embedded Fano threefolds as Brill--Noether loci and state the refined categorical Torelli theorem as a corollary. Finally, in Section~\ref{Section_auto-equi}, we characterize the groups of exact autoequivalences of the Kuznetsov components of genus $6$ and genus $8$ prime Fano threefolds, and provide a new disproof of Kuznetsov's Fano threefold conjecture. Some computational results are proved in the appendix.

\subsection*{Notation and conventions} \leavevmode
\begin{itemize}
    \item Throughout this article, we work over an algebraically closed field $k$ of characteristic zero. The integers $g$ and $d$ denote, respectively, the genus and degree of a Fano threefold. Whenever relevant, we specify the genus or degree of the threefold under consideration. Unless otherwise stated, $X$ denotes an index one Picard rank one Fano threefold of genus $g\geq 6$. When $X$ is written with a subscript, the subscript indicates the degree $2g-2$ of $X$.
  
    \item Let $\sigma$ be a weak stability condition. We denote its central charge and heart by $Z_{\sigma}$ and $\cA_{\sigma}$, respectively. We denote the phase and slope with respect to $\sigma$ by $\phi_{\sigma}$ and $\mu_{\sigma}$, respectively. The maximal and minimal slopes of the Harder--Narasimhan factors of an object are denoted by $\mu_{\sigma}^+$ and $\mu_{\sigma}^-$, respectively.
    
    \item We denote by $\cH_{\sigma}^i$ the $i$-th cohomology object with respect to the heart $\cA_{\sigma}$. When $\cA=\mathrm{Coh}(X)$, we write simply $\cH^i$ for the cohomology objects. All functors in this paper are derived.

    \item For any two objects $A,B$ in a $k$-linear triangulated category, we set $\hom(A,B)\coloneqq \dim_k \Hom(A,B)$ and $\ext^i(A,B)\coloneqq \dim_k \Ext^i(A,B)$.
    
    \item For any object $E$, we denote its numerical class in the numerical Grothendieck group by $[E]$. In our setting, specifying a numerical class is equivalent to specifying a Chern character.
\end{itemize}

\subsection*{Acknowledgements} 
We would like to thank Arend Bayer for suggesting this problem and for many very useful discussions throughout the project. We thank Daniele Faenzi for stimulating conversations on Torelli-type questions for Fano threefolds. We are also grateful to Sasha Kuznetsov, Laurent Manivel, Chunyi Li, Naoki Koseki, Laura Pertusi, and Sebastian Schlegel Mejia for useful discussions. Part of this work was completed while the third author was visiting the Tianyuan Mathematical Center in Southwest China (TMCSC) at Sichuan University. The third author thanks Professor Xiaojun Chen for the invitation and TMCSC for its hospitality. The first and third authors were supported by the ERC Consolidator Grant WallCrossAG, grant no.~819864.

\section{Derived categories} \label{sec:derived categories}

In this section, we collect several useful facts and results concerning semiorthogonal decompositions. Background on triangulated categories and derived categories of coherent sheaves can be found, for example, in \cite{huyb-book-FM}. From now on, $\Db(X)$ denotes the bounded derived category of coherent sheaves on a variety $X$. For $E,F\in\Db(X)$, we set
\[
\RHom^\bullet(E,F)=\bigoplus_{i\in\ZZ}\Ext^i(E,F)[-i].
\]
All triangulated categories in this paper are assumed to be $k$-linear.

\subsection{Exceptional collections and semiorthogonal decompositions}

\begin{Def}
Let $\cD$ be a triangulated category and $E\in\cD$. We say that $E$ is an \emph{exceptional object} if $\RHom^\bullet(E,E)\cong k.$ Let $(E_1,\dots,E_m)$ be a collection of exceptional objects in $\cD$. We say that it is an \emph{exceptional collection} if $\RHom^\bullet(E_i,E_j)=0$ for all $i>j$.
\end{Def}

\begin{Def}
Let $\cD$ be a triangulated category and let $\cC\subset\cD$ be a triangulated subcategory. The \emph{right orthogonal complement} of $\cC$ in $\cD$ is the full triangulated subcategory
\[
\cC^\bot
=
\{X\in\cD \mid \Hom(Y,X)=0 \text{ for all } Y\in\cC\}.
\]
Similarly, the \emph{left orthogonal complement} of $\cC$ in $\cD$ is
\[
{}^\bot\cC
=
\{X\in\cD \mid \Hom(X,Y)=0 \text{ for all } Y\in\cC\}.
\]
\end{Def}

Recall that a triangulated subcategory of $\cD$ is called \emph{admissible} if its inclusion functor admits both left and right adjoints.

\begin{Def}
Let $\cD$ be a triangulated category, and let $(\cC_1,\dots,\cC_m)$ be a collection of full admissible subcategories of $\cD$. We say that
\[
\cD=\langle \cC_1,\dots,\cC_m\rangle
\]
is a \emph{semiorthogonal decomposition} of $\cD$ if $\cC_j\subset \cC_i^\bot$ for all $i>j$, and the subcategories $\cC_1,\dots,\cC_m$ generate $\cD$, i.e. the smallest triangulated subcategory of $\cD$ containing all objects of $\cC_1,\dots,\cC_m$ is equivalent to $\cD$.
\end{Def}

By definition of a Serre functor, we have:

\begin{Lem} \label{serre s.o.d. proposition}
Let $\cD=\langle \cD_1,\cD_2\rangle$ be a semiorthogonal decomposition, and let $S_{\cD}$ be the Serre functor of $\cD$. Then $\cD= \langle S_{\cD}(\cD_2),\cD_1\rangle$ and $\cD= \langle \cD_2,S_{\cD}^{-1}(\cD_1)\rangle$
are also semiorthogonal decompositions.
\end{Lem}

\subsection{Mutations} \label{mutations subsection}

Let $\cC\subset\cD$ be an admissible triangulated subcategory. The inclusion functor
\[
i\colon \cC\hookrightarrow \cD
\]
admits both a left adjoint and a right adjoint, denoted by $i^*$ and $i^!$, respectively. The \emph{left mutation functor} $\bL_{\cC}$ through $\cC$ is defined by the functorial exact triangle
\[
i i^! \ra \identity \ra \bL_{\cC}.
\]
Similarly, the \emph{right mutation functor} $\bR_{\cC}$ through $\cC$ is defined by the functorial exact triangle
\[
\bR_{\cC} \ra \identity \ra i i^*.
\]

When $E\in\Db(X)$ is an exceptional object and $F\in\Db(X)$ is any object, the left mutation $\bL_E F$ fits into the exact triangle
\[
E\otimes\RHom^\bullet(E,F) \ra F \ra \bL_E F,
\]
and the right mutation $\bR_E F$ fits into the exact triangle
\[
\bR_E F \ra F \ra E\otimes\RHom^\bullet(F,E)^\vee.
\]
Moreover, if $(E_1,\dots,E_m)$ is an exceptional collection, then, for $i=1,\dots,m-1$, the collections
\[
(E_1,\dots,E_{i-1},\bL_{E_i}E_{i+1},E_i,E_{i+2},\dots,E_m)
\]
and
\[
(E_1,\dots,E_{i-1},E_{i+1},\bR_{E_{i+1}}E_i,E_{i+2},\dots,E_m)
\]
are also exceptional collections (cf.~\cite[Corollary 2.9]{kuznetsov2010derivedcubicfourfolds}).

\begin{Lem}[{\cite[Lemma 2.6]{ku:fractional-cy}}]\label{lem:serre-sod}
Let $\cD=\langle \cA,\cB\rangle$ be a semiorthogonal decomposition. Then
\[
S_{\cB}=\bR_{\cA}\circ S_{\cD}
\qquad\text{and}\qquad
S_{\cA}^{-1}=\bL_{\cB}\circ S_{\cD}^{-1}.
\]
\end{Lem}

\section{index one prime Fano threefolds and their Kuznetsov components}
\label{sec:derivedfano}

A smooth projective variety with ample anticanonical bundle is called \emph{Fano}. A Fano variety is called \emph{prime} if it has Picard number one. Thus, for a prime Fano variety $X$, there is a unique ample generator $H$ such that $\Pic(X)\cong \ZZ H.$ The \emph{index} of $X$ is the positive integer $i_X$ determined by $-K_X=i_XH.$ The \emph{degree} of a prime Fano threefold is $d_X\coloneqq H^3.$

In this paper, we mainly consider prime Fano threefolds of index one. These threefolds are classified in \cite{shafa}; we recall some of their properties below, following also \cite{KPS}. For an index one prime Fano threefold, the genus $g$ is determined by $H^3=2g-2$. We have:

\begin{itemize}
    \item $2\leq g\leq 12$ and $g\neq 11$;
    
    \item $X^O_{10}$, $g=6$: a section of $\mathbf{G}_X\coloneqq \Gr(2,5)\subset \PP^9$ by a linear space and a quadric;
    
    \item $X^S_{10}$, $g=6$: a double cover of a codimension-$3$ linear section of $\mathbf{G}_X\coloneqq \Gr(2,5)$, branched along an anticanonical divisor;
    
    \item $X_{12}$, $g=7$: a linear section of a connected component of the orthogonal Grassmannian $\mathbf{G}_X\coloneqq \mathrm{OGr}_+(5,10)\subset \PP^{15}$;
    
    \item $X_{14}$, $g=8$: a linear section of $\mathbf{G}_X\coloneqq \Gr(2,6)\subset \PP^{14}$;
    
    \item $X_{16}$, $g=9$: a linear section of the Lagrangian Grassmannian $\mathbf{G}_X\coloneqq \mathrm{LGr}(3,6)\subset \PP^{13}$;
    
    \item $X_{18}$, $g=10$: a linear section of the homogeneous space $\mathbf{G}_X\coloneqq G_2/P\subset \PP^{13}$;
    
    \item $X_{22}$, $g=12$: the zero locus of three sections of the bundle $\wedge^2\cU^\vee$, where $\cU$ is the universal subbundle on $\mathbf{G}_X\coloneqq \Gr(3,7)$.
\end{itemize}

A prime Fano threefold of type $X^S_{10}$ is called a \emph{special Gushel--Mukai (GM)} threefold, while a threefold of type $X^O_{10}$ is called an \emph{ordinary Gushel--Mukai (GM)} threefold.

We denote by $\cE_X$ the pullback of the tautological subbundle on $\mathbf{G}_X$. The corresponding quotient bundle is denoted by $\cQ_X$. When $g\geq 6$ is even, the bundle $\cE_X$ is the unique rank $2$ stable bundle on $X$ with certain constraints (cf.~\cite{bayer2024mukai}).

\subsection{Semiorthogonal decompositions}
\label{SOD}

A series of works of Bondal--Orlov and Kuznetsov shows that index one prime Fano threefolds admit the following semiorthogonal decompositions.

\begin{Prop}
Let $X$ be a prime Fano threefold of index one and genus $g$, and let $\Gamma_{g'}$ be a smooth curve of genus $g'$. Then:
\begin{itemize}
    \item for $g<6$,
    \[
    \Db(X)=
    \left\langle
    \Ku(X),\oh_{X}
    \right\rangle;
    \]
    
    \item for $g=6$,
    \[
    \Db(X)=
    \left\langle
    \Ku(X),\cE_X,\oh_X
    \right\rangle,
    \qquad
    S_{\Ku(X)}=\tau[2],
    \]
    where $\tau$ is an involution of $\Ku(X)$;
    
    \item for $g=7$,
    \[
    \Db(X)=
    \left\langle
    \Ku(X),\cE_X,\oh_X
    \right\rangle,
    \qquad
    \Ku(X)\simeq \Db(\Gamma_7);
    \]
    
    \item for $g=8$,
    \[
    \Db(X)=
    \left\langle
    \Ku(X),\cE_X,\oh_X
    \right\rangle,
    \qquad
    S_{\Ku(X)}^3=[5];
    \]
    
    \item for $g=9$,
    \[
    \Db(X)=
    \left\langle
    \Ku(X),\cE_X,\oh_X
    \right\rangle,
    \qquad
    \Ku(X)\simeq \Db(\Gamma_3);
    \]
    
    \item for $g=10$,
    \[
    \Db(X)=
    \left\langle
    \Ku(X),\cE_{X},\oh_X
    \right\rangle,
    \qquad
    \Ku(X)\simeq \Db(\Gamma_2);
    \]
    
    \item for $g=12$,
    \[
    \Db(X)=
    \left\langle
    \Ku(X),\cE_{X},\oh_X
    \right\rangle,
    \qquad
    \Ku(X)\simeq \Db(\mathrm{Rep}(K(3))),
    \]
    where $K(3)$ is the $3$-Kronecker quiver.
\end{itemize}
\end{Prop}

When $X_{10}$ is special, it is the double cover of a prime Fano threefold of index two and degree $5$. In this case, $\tau$ is induced by the geometric involution of $X_{10}$, which we also denote by $\tau$.

We call the subcategories $\Ku(X)$ defined above the \emph{Kuznetsov components} of $X$. The left adjoint of the inclusion $i\colon \Ku(X)\hookrightarrow \Db(X)$ is called the \emph{projection functor} and is denoted by $i^*$.

The Chern characters of $\cE_X$ are as follows:
\[
\ch(\cE_X)=
\begin{cases}
(2,-H,L,\frac{1}{3}P), & g=6, \\[2pt]
(5,-2H,0,P), & g=7, \\[2pt]
(2,-H,2L,\frac{1}{6}P), & g=8, \\[2pt]
(3,-H,0,\frac{1}{3}P), & g=9, \\[2pt]
(2,-H,3L,0), & g=10, \\[2pt]
(2,-H,4L,-\frac{1}{6}P), & g=12,
\end{cases}
\]
where $L$ and $P$ denote the classes of a line and a point on $X$, respectively. When the genus of $X$ is clear from context, we write simply $\cE\coloneqq\cE_X$.

When $g\geq 6$ is even, by \cite[Proposition 3.9]{kuznetsov2003derived}, the numerical Grothendieck group $\cN(\Ku(X))$ is a rank-two integral lattice generated by
\begin{equation}
\label{v_w}
\cN(\Ku(X))
=
\left\langle
v\coloneqq 1-\frac{g}{2}L+\frac{g-4}{4}P,\,
w\coloneqq H-\frac{3g-6}{2}L+\frac{7g-40}{12}P
\right\rangle .
\end{equation}
With respect to the ordered basis $(v,w)$, the Euler form is given by
\begin{equation}
\begin{pmatrix}
1-\frac{g}{2} & -\frac{g}{2} \\
3-g & 1-g
\end{pmatrix}.
\end{equation}

When $g=7$, using the Hirzebruch--Riemann--Roch theorem, one verifies by direct computation that $\cN(\Ku(X))$ is a rank-two integral lattice generated by
\[
\cN(\Ku(X))
=
\left\langle
v\coloneqq 2-5L+\frac{1}{2}P,\,
w\coloneqq H-6L
\right\rangle .
\]
With respect to the ordered basis $(v,w)$, the Euler form is
\begin{equation}
\begin{pmatrix}
-6 & -5 \\
-7 & -6
\end{pmatrix}.
\end{equation}

When $g=9$, the numerical Grothendieck group $\cN(\Ku(X))$ is a rank-two integral lattice generated by
\[
\cN(\Ku(X))
=
\left\langle
v\coloneqq1-3L+\frac{1}{2}P,\,
w\coloneqq H-8L+\frac{2}{3}P
\right\rangle .
\]
With respect to the ordered basis $(v,w)$, the Euler form is
\begin{equation}
\begin{pmatrix}
-2 & -3 \\
-5 & -8
\end{pmatrix}.
\end{equation}

When $g\geq 6$ is even, we also use the alternative semiorthogonal decomposition
\[
\Db(X)
=
\left\langle
\cA_{X},\oh_{X},\cE^\vee
\right\rangle .
\]
We call $\cA_{X}$ the \emph{alternative Kuznetsov component}. Its numerical Grothendieck group $\cN(\cA_{X})$ is a rank-two integral lattice generated by
\[
\cN(\cA_{X})
=
\left\langle
s\coloneqq 1-2L,\,
t\coloneqq H-\left(\frac{g}{2}+1\right)L-\frac{16-g}{12}P
\right\rangle .
\]
With respect to the ordered basis $(s,t)$, the Euler form is
\begin{equation}
\begin{pmatrix}
-1 & -2 \\
-\frac{g}{2}+1 & -g+1
\end{pmatrix}.
\end{equation}

We call a class $u\in \cN(\Ku(X))$, respectively $u\in \cN(\cA_{X})$, a \emph{$(-r)$-class} if
\[
\chi(u,u)=-r.
\]

\begin{Rem}
\label{rem:equiv kuz components}
When $g\geq 6$ is even, by \cite[Lemma 3.7]{jacovskis2024categorical}, there is an exact equivalence
\[
\Xi\colon \Ku(X)\xrightarrow{\sim}\cA_X,
\qquad
E\longmapsto \mathbf{L}_{\oh_X}(E\otimes \oh_X(H)).
\]
\end{Rem}

\subsection{Gluing objects in Kuznetsov components of even-genus prime Fano threefolds}

In this subsection, we define \emph{gluing objects} associated with the Kuznetsov components $\Ku(X)$ and $\cA_X$ arising from different semiorthogonal decompositions of the derived category of an even-genus prime Fano threefold $X$, where $g\geq 6$. We also compare the resulting gluing objects.

Let $\cD\coloneqq\langle \Ku(X),\cE\rangle$ and $\cD'\coloneqq\langle \cA_X,\cQ^\vee\rangle$, where $\cQ$ is the pullback of the tautological quotient bundle. Consider the inclusions
\[
i_{\cD}\colon \Ku(X)\hookrightarrow \cD,
\qquad
i_{\cD'}\colon \cA_X\hookrightarrow \cD'.
\]
Their left adjoints are $i_{\cD}^*=\bL_{\cE}$ and $i_{\cD'}^*=\bL_{\cQ^\vee}$, and their right adjoints are denoted by $i_{\cD}^!$ and $i_{\cD'}^!$, respectively.

Let $j_*\colon \langle \cE\rangle\hookrightarrow \cD$ denote the inclusion functor. In the sense of \cite[Definition 2.4]{kuznetsov2015categorical}, the functor $i_{\cD}^!\circ j_*$ is called the \emph{gluing functor}. Accordingly, we define the \emph{gluing objects} for $\cD$ and $\cD'$ to be
\[
i_{\cD}^!\cE
\qquad\text{and}\qquad
i_{\cD'}^!\cQ^\vee,
\]
respectively.

Recall that $i^!$ and $i'{}^!$ denote the right adjoints of the inclusions
\[
i\colon \Ku(X)\hookrightarrow \Db(X),
\qquad
i'\colon \cA_X\hookrightarrow \Db(X),
\]
respectively. Since $\cE\in \cD$ and $\cQ^\vee\in \cD'$, we have $i^!\cE=i_{\cD}^!(\cE)$ and $i'{}^!\cQ^\vee=i_{\cD'}^!(\cQ^\vee).$

\begin{Lem}
\label{cohomology objects of pi lemma}
The object $i^!\cE=i_{\cD}^!(\cE)$ is given by
\[
i^!\cE\cong \bL_{\cE}\cQ(-H)[1].
\]
It fits into a triangle
\[\cQ(-H)[1]\to i^!\cE\to \cE.\]
\end{Lem}

\begin{proof}
By the definition of left mutation, see Section~\ref{mutations subsection}, we have an exact triangle
\[
i_{\cD}i_{\cD}^!(\cE)\longrightarrow \cE\longrightarrow \bL_{\Ku(X)}\cE.
\]
By Lemma~\ref{serre s.o.d. proposition}, we have
\[
\langle \Ku(X),\cE\rangle
=
\langle S_{\cD}(\cE),\Ku(X)\rangle
=
\langle \bL_{\Ku(X)}\cE,\Ku(X)\rangle.
\]
Thus, the above triangle can be written as
\[
i_{\cD}i_{\cD}^!(\cE)\longrightarrow \cE\longrightarrow S_{\cD}(\cE).
\]
To compute $S_{\cD}(\cE)$ explicitly, we use Lemma~\ref{lem:serre-sod}, which gives
\[
S_{\cD}\cong \bR_{\oh_X(-H)}\circ S_{\Db(X)}.
\]
Since $\bR_{\oh_X(-H)}\cE(-H)\cong \cQ(-H)[-1],$ we obtain $S_{\cD}(\cE)\cong \cQ(-H)[2].$ Hence, the triangle above becomes
\[
i_{\cD}i_{\cD}^!(\cE)\longrightarrow \cE\longrightarrow \cQ(-H)[2].
\]
Applying $i_{\cD}^*=\bL_{\cE}$ to this triangle, we obtain $i_{\cD}^!(\cE)\cong \bL_{\cE}\cQ(-H)[1].$
\end{proof}

\begin{Lem}[{\cite[Lemma 5.4]{jacovskis2024categorical}}]
\label{gluing-data-equi-alter}
The object $i'{}^!(\cQ^\vee)=i_{\cD'}^!(\cQ^\vee)$ is given by
\[
i'{}^!(\cQ^\vee)\cong \bL_{\cQ^\vee}\cE[1].
\]
It fits into a triangle
\[\cE[1]\to i'{}^!(\cQ^\vee)\to \cQ^{\vee}.\]
\end{Lem}

\begin{proof}
The proof is analogous to that of Lemma~\ref{cohomology objects of pi lemma}. We have an exact triangle
\[
i_{\cD'}i_{\cD'}^!(\cQ^\vee)
\longrightarrow
\cQ^\vee
\longrightarrow
\bL_{\cA_X}\cQ^\vee.
\]
Moreover,
\[
\langle \cA_X,\cQ^\vee\rangle
\simeq
\langle S_{\cD'}(\cQ^\vee),\cA_X\rangle
\simeq
\langle \bL_{\cA_X}\cQ^\vee,\cA_X\rangle.
\]
Therefore the triangle above becomes
\[
i_{\cD'}i_{\cD'}^!(\cQ^\vee)
\longrightarrow
\cQ^\vee
\longrightarrow
S_{\cD'}(\cQ^\vee).
\]
To compute $S_{\cD'}(\cQ^\vee)$, we use the description of the Serre functor from Lemma~\ref{lem:serre-sod}. The standard mutation computation gives $S_{\cD'}(\cQ^\vee)\cong \cE[2].$ Thus, we obtain a triangle
\[
i_{\cD'}i_{\cD'}^!(\cQ^\vee)
\longrightarrow
\cQ^\vee
\longrightarrow
\cE[2].
\]
Applying $i_{\cD'}^*=\bL_{\cQ^\vee}$ and using $i_{\cD'}^*\cQ^\vee=0$, we obtain $i_{\cD'}^!(\cQ^\vee)\cong \bL_{\cQ^\vee}\cE[1].$ As $\RHom^\bullet(\cQ^\vee,\cE)= k[-2],$ the defining triangle for the left mutation is
\[
\cQ^\vee[-2]\longrightarrow \cE\longrightarrow \bL_{\cQ^\vee}\cE
\]
and the result follows.
\end{proof}

\begin{Rem}
\label{gluing_data_equal}
One checks that
\[
\Xi(i^!\cE)\cong i'{}^!(\cQ^\vee)[1],
\]
where $\Xi$ is the equivalence $\Ku(X)\simeq \cA_X$ from Remark~\ref{rem:equiv kuz components}. Indeed, applying $\Xi$ to the triangle defining $i^!\cE$ yields the triangle defining $i'{}^!(\cQ^\vee)$.
\end{Rem}

\begin{Rem}
The gluing object in the Kuznetsov component of an odd-genus prime Fano threefold is defined similarly. These objects were already considered in \cite[Lemma 3.5]{brambilla2013rank} for genus $9$ prime Fano threefolds, and in \cite[Lemma 2.9]{brambilla2014vector}. In those references, the Kuznetsov component is defined as $\mathcal{B}_X\coloneqq{}^{\perp}\langle \oh_X,\cE^\vee\rangle,$ and one has $\mathcal{B}_X\simeq \cA_X\otimes \oh_X(H).$
\end{Rem}

\section{Bridgeland stability conditions} \label{sec:stability conditions}

In this section, we recall the construction of (weak) Bridgeland stability conditions on $\Db(X)$, and the notions of tilt stability, double-tilt stability, and stability conditions induced on Kuznetsov components from weak stability conditions on $\Db(X)$. We follow \cite[\S~2]{bayer2017stability}.

\subsection{Weak stability conditions}

Let $\cD$ be a triangulated category, and $K_0(\cD)$ its Grothendieck group. Fix a surjective morphism $v \colon K_0(\cD) \ra \Lambda$ to a finite rank lattice. 

\begin{Def}
The \emph{heart of a bounded t-structure} on $\cD$ is an abelian subcategory $\cA \sst \cD$ such that the following conditions are satisfied:
\begin{enumerate}[(i)]
    \item for any $E, F \in \cA$ and $n <0$, we have $\Hom(E, F[n])=0$;
    \item for any object $E \in \cD$ there exist objects $E_i \in \cA$ and maps
    \[ 0=E_0 \xrightarrow{s_1} E_1 \xrightarrow{s_2} \cdots \xra{s_m} E_m=E \]
    such that $\cone(s_i) = A_i[k_i]$ where $A_i \in \cA$ and the $k_i$ are integers such that $k_1 > k_2 > \cdots > k_m$.
\end{enumerate}
\end{Def}

\begin{Def}
Let $\cA$ be an abelian category and $Z \colon K_0(\cA) \ra \CC$ be a group homomorphism such that for any $E \in \cD$ we have $\Im Z(E) \geq 0$ and if $\Im Z(E) = 0$ then $\Re Z(E) \leq 0$. Then we call $Z$ a \emph{weak stability function} on $\cA$. If furthermore we have for $0\neq E\in \cA$ that $\Im Z(E) = 0$ implies that $\Re Z(E) < 0$, then we call $Z$ a \emph{stability function} on $\cA$.
\end{Def}

\begin{Def}
A \emph{weak stability condition} on $\cD$ is a pair $\sigma = (\cA, Z)$ where $\cA$ is the heart of a bounded t-structure on $\cD$, and $Z \colon \Lambda \ra \CC$ is a group homomorphism such that 
\begin{enumerate}[(i)]
    \item the composition $Z \circ v \colon K_0(\cA) \cong K_0(\cD) \ra \CC$ is a weak stability function on $\cA$. From now on, we write $Z(E)$ rather than $Z(v(E))$.
\end{enumerate}
Much like the slope in classical $\mu$-stability, we can define a \emph{slope} $\mu_\sigma$ for $\sigma$ using $Z$. For any $E \in \cA$, set
\[
\mu_\sigma(E) \coloneqq \begin{cases}  - \frac{\Re Z(E)}{\Im Z(E)}, & \text{if} ~ \Im Z(E) > 0 \\
+ \infty , & \text{otherwise}.
\end{cases}
\]
We say an object $0 \neq E \in \cA$ is $\sigma$-(semi)stable if $\mu_\sigma(F) < \mu_\sigma(E/F)$ (respectively $\mu_\sigma(F) \leq \mu_\sigma(E/F)$) for all proper subobjects $F \sst E$. 
\begin{enumerate}[(i), resume]
    \item Any object $E \in \cA$ has a Harder--Narasimhan filtration in terms of $\sigma$-semistability defined above.
    \item There exists a quadratic form $Q$ on $\Lambda \otimes \mathbb{R}$ such that $Q|_{\ker Z}$ is negative definite, and $Q(E) \geq 0$ for all $\sigma$-semistable objects $E \in \cA$. This is known as the \emph{support property}.
\end{enumerate}
If the composition $Z \circ v$ is a stability function, then $\sigma$ is a \emph{stability condition} on $\cD$.
\end{Def}

\begin{Def}
Let $\sigma=(\cA, Z)$ be a weak stability condition on $\cD$. The \emph{phase} of a $\sigma$-semistable object $E\in \cA$ is
\[\phi(E) \coloneqq \frac{1}{\pi} \mathrm{arg}(Z(E)) \in (0,1].\]
In particular, if $Z(E) = 0$ then $\phi(E) = 1$. If $F = E[n]$, then we define 
\[\phi(F) \coloneqq \phi(E) + n.\]

We also define a collection of full additive subcategories $\cP(\phi) \subset \cD$ for each $\phi \in \mathbb{R}$ by
\begin{enumerate}[resume]
\item for $\phi\in (0,1]$, the subcategory $\cP(\phi)$ is given by the zero object and all $\sigma$-semistable objects whose phase is $\phi$;
\item for $\phi + n$ with $\phi\in (0,1]$ and $n\in \mathbb{Z}$, we set $\cP(\phi + n) \coloneqq \cP(\phi)[n]$.
\end{enumerate}
Then by the existence of Harder--Narasimhan filtration, for any non-zero $E\in \cD$, there exists a sequence of morphisms 
\[ 0=E_0 \xrightarrow{s_1} E_1 \xrightarrow{s_2} \cdots \xra{s_m} E_m=E \]
    such that $\cone(s_i) \in \cP(\phi_i)$ for a sequence of real numbers $\phi_1 > \phi_2 > \cdots > \phi_m$. We set $\phi_{\sigma}^+(E)\coloneqq \phi_1$ and $\phi_{\sigma}^-(E)\coloneqq \phi_m$.
\end{Def}

For a stability condition $\sigma=(\cA, Z)$ on $\cD$, we define the \emph{homological dimension} of $\cA$ as the smallest integer $\mathrm{homdim}(\cA)$ such that $\Hom(A,B[n])=0$ for any $n>\mathrm{homdim}(\cA)$. The homological dimension of a stability condition is defined as the homological dimension of its heart.

\subsection{Tilt stability}

Let $\sigma = (\cA, Z)$ be a weak stability condition on a triangulated category $\cD$. Now consider the following subcategories of $\cA$, where $\langle - \rangle$ denotes the extension closure:
\begin{align*} 
\cT_\sigma^\mu &= \langle  E \in \cA \mid E \text{ is } \sigma\text{-semistable with } \mu_\sigma(E) > \mu \rangle \\
\cF_\sigma^\mu &= \langle E \in \cA \mid E \text{ is } \sigma\text{-semistable with } \mu_\sigma(E) \leq \mu  \rangle.
\end{align*}
Then it is a result of \cite{happel1996tilting} that

\begin{Prop}
The abelian category $\cA^\mu_\sigma \coloneqq \langle \cT_\sigma^\mu, \cF_\sigma^\mu[1] \rangle$ is the heart of a bounded t-structure on $\cD$.
\end{Prop}

We call $\cA^\mu_\sigma$ the \emph{tilt} of $\cA$ around the torsion pair $(\cT_\sigma^\mu, \cF_\sigma^\mu)$. Let $X$ be an $n$-dimensional smooth projective complex variety. Tilting can be applied to the weak stability condition $(\Coh(X) , Z_H)$ to form the once-tilted heart $\Coh^\beta(X)$, where  $Z_H(E)\coloneqq-c_1(E)H^{n-1}+\mathfrak{i} \rk(E)H^n$ for any $E\in \Coh(X)$. The associated slope function is denoted by $\mu_H(-)$. For $E \in \Coh^\beta(X)$, we define
\[ Z_{\alpha, \beta}(E) \coloneqq \frac{1}{2} \alpha^2 H^n \ch_0^\beta(E) - H^{n-2} \ch_2^\beta(E) + \mathfrak{i} H^{n-1} \ch_1^\beta(E). \]

\begin{Prop}[{\cite{bayer2011bridgeland, bayer2016space}}]
Let $\alpha>0$ and $\beta \in \mathbb{R}$. Then the pair $\sigma_{\alpha, \beta} = (\Coh^\beta(X) , Z_{\alpha, \beta})$ defines a weak stability condition on $\Db(X)$. We have
\[ \Delta_H(E) = (H^{n-1} \ch_1(E))^2 - 2 H^n \ch_0(E) H^{n-2} \ch_2(E) \geq 0  \]
for any $\sigma_{\alpha,\beta}$-semistable object.
\end{Prop}

We denote by $\mu_{\alpha, \beta}(-)$ the slope function associated with $\sigma_{\alpha,\beta}$.

We now state a useful lemma which relates 2-Gieseker-stability and tilt stability. 

\begin{Lem} [{\cite[Lemma 2.7]{bayer2016space}, \cite[Proposition 4.8]{bayer2020desingularization}}] \label{bms lemma 2.7}
Let $E \in \Db(X)$ be an object.
If $H^2 \ch_1^\beta(E) > 0$, then $E\in \Coh^{\beta}(X)$ and $E$ is $\sigma_{\alpha, \beta}$-(semi)stable for $\alpha \gg 0$ if and only if $E$ is a 2-Gieseker-(semi)stable sheaf.

\end{Lem}

\subsection{Finding solutions for walls in tilt stability} \label{sec:tilt stability wall algorithm}

In this subsection, we describe a way of finding (potential) walls in tilt stability with respect to objects in the derived category with a given truncated Chern character. This is similar to the method used, e.g., in \cite{PY20} to find walls for certain objects. 

Let $X$ be a prime Fano threefold and $M \in \Coh^{\beta}(X)$ be an object. Assume there is a short exact sequence $0 \ra E \ra M \ra F \ra 0$ which makes $M$ strictly $\sigma_{\alpha,\beta}$-semistable. We can assume that $E$ and $F$ are tilt-semistable using the existence of Harder--Narasimhan or Jordan--H\"older filtrations. Then the following conditions must be satisfied:
\begin{enumerate}
    \item $ \ch_{\leq 2}(M) = \ch_{\leq 2}(E) + \ch_{\leq 2}(F)$;
    \item $\mu_{\alpha, \beta}(E) = \mu_{\alpha, \beta}(M) = \mu_{\alpha, \beta}(F)$;
    \item $\Delta_H(E) \geq 0$ and $\Delta_H(F) \geq 0$;
    \item $\Delta_H(E) \leq \Delta_H(M)$ and $\Delta_H(F) \leq \Delta_H(M)$.
\end{enumerate}

Since $E, F \in \Coh^\beta(X)$, we also have $\ch_1^\beta(E) \geq 0$ and $\ch_1^\beta(F) \geq 0$. Solving the system of inequalities above gives an even number of solutions of $(m_0, m_1, m_2)\in \mathbb{Z}^{\oplus 3}$, where
\[\ch_{\leq 2}(E)=\left(m_0,m_1H, m_2\frac{H^2}{H^3}\right).\]
Half of them are solutions for the destabilizing subobject $E$, and the other half are the corresponding quotients $F$.

\subsection{Stability conditions on Kuznetsov components}

\subsubsection{Double-tilted stability conditions}

Now as in \cite{bayer2017stability}, we pick a weak stability condition $\sigma_{\alpha, \beta}$ and tilt the once-tilted heart $\Coh^\beta(X)$ with respect to the tilt slope $\mu_{\alpha, \beta}$ and some second tilt parameter $\mu$. One gets a torsion pair $(\cT_{\alpha, \beta}^\mu, \cF_{\alpha, \beta}^\mu)$ and another heart $\Coh^\mu_{\alpha, \beta}(X)$ of $\Db(X)$. Now ``rotate" the stability function $Z_{\alpha, \beta}$ by setting 
\[ Z_{\alpha, \beta}^\mu \coloneqq \frac{1}{u} Z_{\alpha, \beta} \]
where $u \in \CC$ such that $|u|=1$ and $\mu =-\frac{\Re u}{\Im u}$.

\begin{Prop}[{\cite[Proposition 2.15]{bayer2017stability}}]
The pair $\sigma^{\mu}_{\alpha,\beta}\coloneqq (\Coh^\mu_{\alpha, \beta}(X), Z_{\alpha, \beta}^\mu)$ defines a weak stability condition on $\Db(X)$.
\end{Prop}

For example, if we choose $\mu=0$, we have
\[ Z^0_{\alpha, \beta}(E) = H^{n-1} \ch_1^\beta(E) + \mathfrak{i} \left(H^{n-2} \ch_2^\beta(E)-\frac{1}{2} \alpha^2 H^n \ch_0^\beta(E)\right). \]
The associated slope function is denoted by $\mu^0_{\alpha,\beta}(-)$.

\subsubsection{Stability conditions on Kuznetsov components}

More precisely, let $\cA(\alpha, \beta) = \Coh^\mu_{\alpha, \beta}(X) \cap \Ku(X)$ and $Z(\alpha, \beta) = Z_{\alpha, \beta}^\mu|_{\Ku(X)}$. Furthermore, let $0 < \epsilon \ll 1$, $\beta =  -1+\epsilon$ and $0 < \alpha < \epsilon$. Also impose the following condition on the second tilt parameter $\mu$:
\begin{equation}
    \mu_{\alpha, \beta}(\cE(-H)[1]) < \mu_{\alpha, \beta}(\oh_X(-H)[1]) < \mu < \mu_{\alpha, \beta}(\cE) < \mu_{\alpha, \beta}(\oh_X) .
\end{equation}
Then we get the following theorem.

\begin{Thm}[{\cite[Theorem 6.9]{bayer2017stability}}]
Let $X$ be a Fano threefold of genus $6, 8, 10$ or $12$, and let $\epsilon, \alpha, \beta$ and $\mu$ be as above. Then the pair $\sigma(\alpha, \beta)=(\cA(\alpha, \beta), Z(\alpha, \beta))$ defines a Bridgeland stability condition on $\Ku(X)$.
\end{Thm}

In our paper, we fix $\mu=0$, i.e. $\sigma(\alpha, \beta)\coloneqq\sigma^0_{\alpha, \beta}|_{\Ku(X)}$.

\begin{Prop}

\label{induce_range}
Let $X$ be a prime Fano threefold of index one and genus $g\geq 6$. Then $\sigma(\alpha, \beta)\coloneqq(\cA^0_{\alpha, \beta}|_{\Ku(X)}, Z^0_{\alpha, \beta}|_{\Ku(X)})$ is a stability condition for $(\alpha, \beta)$ listed below:

\begin{itemize}
    \item $g=6$: $\beta=-\frac{9}{10}, 0<\alpha<1+\beta$,
    
    \item  $g=7$: $\beta=-\frac{5}{6}$ or $-\frac{71}{84}$,  $0<\alpha<1+\beta$,
    
    \item $g=8$: $\beta=-\frac{22}{25}$ or $-\frac{122}{125}$,  $0<\alpha<1+\beta$,
    
    \item  $g=9$: $\beta=-\frac{3}{4}$ or $-\frac{31}{40}$,  $0<\alpha<1+\beta$,
    
    \item  $g=10$: $\beta=-\frac{22}{25}$ or $-\frac{10}{11}$,  $0<\alpha<1+\beta$,
    
    \item  $g=12$: $\beta=-\frac{21}{25}$ or $-\frac{19}{22}$,  $0<\alpha<1+\beta$. 
\end{itemize}

Moreover, $\sigma(\alpha, \beta)$ is Serre-invariant for these $(\alpha, \beta)$.
\end{Prop}

\begin{proof}
It is not hard to see that $\cE, \cE(-H)[1], \oh_X, \oh_X(-H)[1]\in \Coh^{\beta}(X)$, and that they satisfy
\[\mu_{\alpha, \beta}(\cE(-H)[1]) < \mu_{\alpha, \beta}(\oh_X(-H)[1]) < 0 < \mu_{\alpha, \beta}(\cE) < \mu_{\alpha, \beta}(\oh_X) \]
for each $(\alpha, \beta)$ listed above.

First we assume that $g\geq 6$ is even. Then from \cite[Proposition 3.2]{Pertusi2021serreinv} we know that $\sigma(\alpha, \beta)$ is a stability condition for $(\alpha, \beta)$ as above. The Serre-invariance follows from \cite[Theorem 3.18]{Pertusi2021serreinv}.

When $g=7$ or $9$, we know that $\Ku(X)$ is equivalent to the derived category of a certain curve with positive genus. Thus if one proves that $\sigma(\alpha, \beta)$ is a stability condition, the Serre-invariance follows from \cite{macri2007stability}. To this end, it is sufficient to show that $\cE$ and $\cE(-H)[1]$ are $\sigma_{\alpha, \beta}$-semistable for the $(\alpha, \beta)$ listed above. From Lemma \ref{bms lemma 2.7} we know that $\cE$ and $\cE(-H)[1]$ are both $\sigma_{\alpha, \beta}$-semistable for the $\beta$ listed above and for $\alpha \gg 0$. The result then follows from Lemma \ref{g=7wall1}, Lemma \ref{g=7wall2}, Lemma \ref{g=9wall1} and Lemma \ref{g=9wall2}.
\end{proof}

\subsection{Serre-invariance of stability conditions on Kuznetsov components}

Recall the universal covering $\widetilde{\mathrm{GL}}^+(2,\mathbb{R})$ of $\mathrm{GL}^+(2, \mathbb{R})$ acts on the space of stability conditions, see \cite[Lemma 8.2]{bridgeland}.

\begin{Def}\label{Serre_invariant_stab}
Let $\sigma$ be a stability condition on a triangulated category $\cD$. It is called \emph{Serre-invariant} if $S_{\cD} \cdot \sigma=\sigma \cdot g$ for some $g\in\widetilde{\mathrm{GL}}^+(2,\mathbb{R})$.
\end{Def}

By \cite{PY20,Pertusi2021serreinv}, the stability conditions considered in Proposition \ref{induce_range} are all Serre-invariant. In the following, we recall several properties of Serre-invariant stability conditions.

\begin{Prop}\label{ogm homo dim 2}
Let $\sigma$ be a Serre-invariant stability condition on $\Ku(X)$, where $X$ is a prime Fano threefold of index one and genus $g\geq 6$. Then:

\begin{enumerate}[\normalfont(i)]
    \item the heart of $\sigma$ has homological dimension $\leq 2$, and
    \item let $A$ be a non-trivial object in the heart of $\sigma$ and $g\geq 6$. Then \begin{itemize}
        \item $\mathrm{ext}^1(A,A)\geq 2$ if $g=6,8,12$ and
        \item $\mathrm{ext}^1(A,A)\geq 1$ if $g=7,9,10$.
    \end{itemize}
\end{enumerate}
\end{Prop}

\begin{proof}
We treat the proof case by case.
\begin{enumerate}[(i)]

    \item If $g=7,9,10$, we have $\Ku(X)\simeq \Db(C)$ for certain curves $C$ of positive genus (see \cite{kuznetsov2005derived,kuznetsov:fano-threefolds}). Then by \cite{macri2007stability}, the stability condition is given by slope stability on the curve $C$ up to some action of $\widetilde{\mathrm{GL}}^+(2,\mathbb{R})$, whose heart is $\mathrm{Coh}(C)$. So the homological dimension is $1$. If $g=12$, then by \cite[Section 7]{dimitrov2019bridgeland} the heart of a stability condition on $\Ku(X_{22})\simeq \Db(\mathrm{Rep}(K(3)))$ is generated by two exceptional objects, so the homological dimension is $1$. If $g=6$ or $8$, the result follows from \cite[Proposition 4.13]{zhang2020bridgeland} and \cite[Lemma 5.10]{PY20}.
    
    \item Since $\chi(A,A)\geq -1$ for every non-trivial object $A$ in $\Ku(X_{10})$, $\Ku(X_{14})$ and $\Ku(X_{22})$, by $(i)$ we have $\mathrm{ext}^1(A,A)\geq 2$. For $g=7,9$ and $10$, we get $\mathrm{ext}^1(A,A)\geq 1$ since $\chi(A,A)\geq 0$ for every non-trivial object $A$ in $\Ku(X_{2g-2})$.
\end{enumerate}
\end{proof}

\begin{Lem}[{Weak Mukai Lemma, \cite[Lemma 5.12]{PY20}, \cite[Lemmas 3.15 and 3.16]{liu2021note}}] \label{weak-mukai-lemma}
 Let $X$ be a prime Fano threefold of index one and genus $g$, and $\sigma=(\cA_{\sigma}, Z_{\sigma})$ be a Serre-invariant stability condition on $\Ku(X)$. Let $A\rightarrow E\rightarrow B$ be an exact triangle in $\Ku(X)$. Under one of the following assumptions, we have
\[\mathrm{ext}^1(A,A)+\mathrm{ext}^1(B,B)\leq\mathrm{ext}^1(E,E).\]
 
\leavevmode\begin{enumerate}
   \item $g=8$, $\mathrm{Hom}(A,B)=0$, and $\phi^-_{\sigma}(A)\geq \phi^+_{\sigma}(B)$.
   \item $g=7,9,10,12$, $A,B\in \cA_{\sigma}$, and $\mathrm{Hom}(A,B)=0$.
   
    \item $g=6$ and $\mathrm{Hom}(A,B)\cong\mathrm{Hom}(A,\tau(B))=0$.
\end{enumerate}
\end{Lem}

By the following result, all Serre-invariant stability conditions on Kuznetsov components of certain Fano threefolds are unique up to a group action.

\begin{Thm}[{\cite[Theorem A.10]{jacovskis2024categorical}, \cite[Theorem 3.1]{FeyzbakhshPertusi2021stab}}]
\label{all_in_one_orbit}

Let $X$ be a prime Fano threefold of genus $g\geq 6$. Then all Serre-invariant stability conditions on $\mathcal{K}u(X)$ are in the same $\widetilde{\mathrm{GL}}^+(2,\mathbb{R})$-orbit.
\end{Thm}

\section{Embedding Fano threefolds into Bridgeland moduli spaces} \label{sec:embedding}

In this section, we embed our Fano threefolds $X$ into certain moduli spaces of stable objects in $\Ku(X)$. In the first subsection, we give an explicit description of projections of skyscrapers into $\Ku(X)$. In the second subsection, we deal with the stability of these projections with respect to stability conditions on $\Ku(X)$. In the last subsection, we construct an embedding from $X$ to moduli spaces of stable objects and verify Mukai's conjecture.

\subsection{Projections of skyscraper sheaves}

Recall that $n_g\coloneqq\frac{g}{2}$ if $g$ is even, $n_7=5$, and $n_9=3$.

\begin{Prop} \label{projection}
Let $X$ be a prime Fano threefold of genus $g\geq 6$ and $x\in X$. Then we have $i^*(\oh_x)=\bL_{\cE}(I_x)[1]$. Moreover, if we set 
\[K_x\coloneqq \ker(\cE^{\oplus n_g} \to I_x),\]
where $h\colon \cE^{\oplus n_g} \to I_x$ is the natural evaluation map, then
\begin{enumerate}[\normalfont(i)]
    \item $\bL_{\cE}(I_x)\cong K_x[1]$ when $X$ is not a special GM threefold, and
     \item we have $K_x=\ker(\cE^{\oplus 3}\twoheadrightarrow I_{\pi^{-1}(\pi(x))})$ and a triangle
     \[K_{x}[1]\to \bL_{\cE}(I_x) \to \oh_{\tau(x)}\]
     when $X$ is a special GM threefold, $\pi$ is the natural degree $2$ map $\pi\colon X\to \Gr(2,5)$, and $\tau$ is the geometric involution on $X$ induced by $\pi$.
\end{enumerate}

\end{Prop}

\begin{proof}
It is clear that $\bL_{\oh_X} \oh_x \cong I_x[1]$, so $i^*(\oh_x)=\bL_{\cE}(I_x)[1]$. We divide the remaining proof into several cases.

\begin{itemize}
    \item Suppose $g\geq 6$ is even and $X$ is not a special GM threefold. Then we have an embedding $X \hookrightarrow \Gr(2, \frac{g}{2}+2)$. Applying $\Hom(\cE, -)$ to the exact sequence $0 \ra I_x \ra \oh_X \ra \oh_x \ra 0$, we get an exact sequence
\begin{equation*}
    0 \ra \Hom(\cE, I_x) \ra \Hom(\cE, \oh_X) \ra \Hom(\cE, \oh_x) \ra \Ext^1(\cE, I_x) \ra 0 .
\end{equation*}
Note that $\RHom^\bullet(\cE, \oh_X) = k^{\frac{g}{2}+2}$ and $\RHom^\bullet(\cE, \oh_x)=k^2$. 
Then $\mathrm{hom}(\cE,I_x)\geq\frac{g}{2}$. If $\mathrm{hom}(\cE,I_x)\geq\frac{g}{2}+1$, then $x$ would be contained in the zero locus of at least $\frac{g}{2}+1$ linearly independent sections of $\cE^{\vee}$, which is an empty set. Thus $\mathrm{hom}(\cE,I_x)=\frac{g}{2}$. Therefore $\RHom^\bullet(\cE, I_x) = k^{\frac{g}{2}}$ and we have a triangle
\begin{equation} \label{eq:even genus projection original triangle}
    \cE^{\oplus \frac{g}{2}} \xra{h} I_x \ra \bL_{\cE} I_x .
\end{equation}
We claim that the map $h$ is surjective. Indeed, $\im(h) = I_D$ where $D$ is the zero locus of $\frac{g}{2}$ linearly independent sections of $\cE^{\vee}$ containing the point $x$. For all cases of $g$ that we consider, this zero locus is $\Gr(2,2)$, which is just a point. Since $X \hookrightarrow \Gr(2, \frac{g}{2}+2)$ is an embedding, it follows that $\im(h) = I_x$. So $h$ is surjective and $\bL_{\cE} I_x \cong \ker(h)[1]$.
\item If $g=7$, $i^*(\oh_x)\cong\bL_{\cE}I_x[1]$. By similar computations, $\mathrm{RHom}^\bullet(\cE,I_x)\cong k^5$, so that we have
$$\cE^{\oplus 5}\xra{h} I_x\rightarrow\bL_{\cE}I_x.$$
If the map $h$ is not surjective, then $\mathrm{im}(h)$ would be an ideal sheaf $I_D$, where $D$ is the zero locus of five linearly independent sections of $\cE^{\vee}$. But in this case $D=\mathrm{Gr}(5,5)\cap X$, and since $X\hookrightarrow \mathrm{OGr}_+(5,10)$ is cut out by a linear section, we know that $D=\{x\}$. Hence the map is surjective as well. 
\item If $g=9$, then $X$ is a linear section of the Lagrangian Grassmannian $\mathrm{LGr}(3,6)\subset \mathbb{P}^{13}$. The map $h\colon \cE^{\oplus 3}\rightarrow I_x$ is surjective by \cite[Lemma 3.6]{brambilla2013rank}.

\item If $X$ is a special GM threefold, since $X$ does not embed into $\Gr(2,5)$, the map $h$ 
is no longer a surjective map. Instead, its image is $I_{\pi^{-1}(\pi(x))}$. Moreover, we have $\pi^{-1}(\pi(x))=x\cup\tau(x)$ or $\pi^{-1}(\pi(x))$ is a non-reduced point with multiplicity two. Then the desired result follows.
\end{itemize}
\end{proof}

\begin{Rem} \label{rem:K shorthand}
It is an easy computation to show that if $g\geq 6$ and $g$ is even, then 
\[\ch(i^*\oh_x)=(g-1)v-\frac{g}{2}w =\left(  g-1, -\frac{g}{2}H, \frac{g(g-4)}{4}L, -\frac{1}{24}(g+2)(g-12)P \right). \]
If $g=7$, then $\ch(i^*\oh_x)=12v-10w = (24,-10H, 0, 6P)$ and if $g=9$, then $\ch(i^*\oh_x)=8v-3w = (8,-3H,0,2P)$.
\end{Rem}

\subsection{Stability of projections of skyscraper sheaves}

Now, we show that $i^*\oh_x$ is stable with respect to every Serre-invariant stability condition on $\Ku(X)$. We begin with the following stability result for $K_x$.

\begin{Lem} \label{Kstable}
Let $X$ be a prime Fano threefold of even genus $g\geq 6$. Then $K_x$ is a $\mu_H$-stable reflexive sheaf for any $x\in X$. 
\end{Lem}

\begin{proof}
It is clear that $K_x$ is reflexive. Assume that $K_x$ is not $\mu_H$-semistable, and let $K'\subset K_x$ be the maximal destabilizing subsheaf. Since $K_x\subset \cE^{\oplus \frac{g}{2}}$ and $\cE$ is $\mu_H$-stable, we know $$-\frac{g}{2g-2}<\mu_H(K')\leq -\frac{1}{2}.$$ It is not hard to check that the only possible case is  $\mu_H(K')=-\frac{1}{2}$. As $K'$ is reflexive and $\cE$ is locally free and $\mu_H$-stable, we know that $K'$ is a direct sum of $\cE$, which contradicts $\Hom(\cE, K_x)=0$.
\end{proof}

\begin{Rem}
The stability of $K_x$ in odd genus can be proved by a similar argument as above, but we will not need this result.
\end{Rem}

To apply Lemma \ref{weak-mukai-lemma}, we need the following computation of Ext groups.

\begin{Lem} \label{ext of K}
Let $X$ be a prime Fano threefold of genus $g\geq 6$ and $x\in X$. We have
\[ 
\RHom^\bullet(i^*\oh_x, i^*\oh_x) = \begin{cases}
k\oplus k^6[-1], & \text{if } g=6 \text{ and ordinary} \\
k\oplus k^6[-1] ~\text{or}~k\oplus k^7[-1]\oplus k[-2] , & \text{if } g=6 \text{ and special} \\
k\oplus k^{25}[-1], & \text{if } g=7 \\
k\oplus k^9[-1], & \text{if } g=9 \\
k\oplus k^g[-1], & \text{if } g\geq 8 \text{ and } g \text{ is even}. 
\end{cases}
\]
\end{Lem}

\begin{proof}

First, we assume that $g\geq 6$ is even. Then we have a triangle $\cE^{\oplus \frac{g}{2}}\to I_x\to i^*\oh_x[-1]$. We use the standard spectral sequence, see e.g. \cite[Lemma 2.27]{pirozhkov2020admissible}, to do the computations. The first page is
	
$$ E^{p,q}_1= \begin{array}{c|cc}
	     \Ext^3(I_x, \cE^{\oplus \frac{g}{2}}) &  \Ext^3(\cE^{\oplus \frac{g}{2}}, \cE^{\oplus \frac{g}{2}})\oplus \Ext^3(I_x, I_x) & \Ext^3(\cE^{\oplus \frac{g}{2}},I_x)\\
	     \Ext^2(I_x, \cE^{\oplus \frac{g}{2}}) &  \Ext^2(\cE^{\oplus \frac{g}{2}}, \cE^{\oplus \frac{g}{2}})\oplus \Ext^2(I_x, I_x) & \Ext^2(\cE^{\oplus \frac{g}{2}},I_x)\\
	     \Ext^1(I_x, \cE^{\oplus \frac{g}{2}}) & \Ext^1(\cE^{\oplus \frac{g}{2}}, \cE^{\oplus \frac{g}{2}})\oplus \Ext^1(I_x, I_x) & \Ext^1(\cE^{\oplus \frac{g}{2}},I_x) \\
	     \Hom(I_x, \cE^{\oplus \frac{g}{2}}) & \Hom(\cE^{\oplus \frac{g}{2}}, \cE^{\oplus \frac{g}{2}})\oplus \Hom(I_x, I_x) & \Hom(\cE^{\oplus \frac{g}{2}},I_x) \\ \hline
	     0 &  0 & 0
	\end{array}$$
	
It is not hard to show that the dimensions appearing in the first page are of the form
	
	$$\dim E^{p,q}_1= \begin{array}{c|cc}
	     0 &  0 & 0\\
	     g &  3 & 0\\
	     0 & 3 & 0 \\
	     0 & \frac{g^2}{4}+1 & \frac{g^2}{4} \\ \hline
	     0 &  0 & 0
	\end{array}$$

From $i^*\oh_x\in \Ku(X)$, if we apply $\Hom(\cE^{\oplus \frac{g}{2}},-)$ to the triangle $\cE^{\oplus \frac{g}{2}}\to I_x\to i^*\oh_x[-1]$, we have a natural isomorphism $\Hom(\cE^{\oplus \frac{g}{2}}, \cE^{\oplus \frac{g}{2}}) \cong \Hom(\cE^{\oplus \frac{g}{2}}, I_x)$. This implies that the differential $$E^{0,0}_1=k^{\frac{g^2}{4}+1}\to E^{1,0}_1=k^{\frac{g^2}{4}}$$ is surjective. Thus, we get $E_2^{0,0}=k$ and $E^{1,0}_2=0$.

Next, we compute the differential $d\colon E^{-1,2}_1\to E^{0,2}_1$. Since $\Ext^2(\cE, \cE)=0$, the map $$d\colon E^{-1,2}_1=\Ext^2(I_x, \cE^{\oplus \frac{g}{2}})\to E^{0,2}_1=\Ext^2(I_x, I_x)$$ is the natural map induced by applying $\Hom(I_x,-)$ to the triangle $\cE^{\oplus \frac{g}{2}}\to I_x\to i^*\oh_x[-1]$. When $X$ is not a special GM threefold, we have $\Ext^2(I_x,i^*\oh_x[-1])=\Ext^3(I_x, K_x)$. By Serre duality and Lemma \ref{Kstable}, we have \[ \Ext^2(I_x,i^*\oh_x[-1])=\Ext^3(I_x, K_x) \cong \Hom(K_x, I_x(-H))=0 . \] This means the natural map $d\colon E^{-1,2}_1\to E^{0,2}_1$ is surjective. Thus, the dimensions of the second page are of the form

	$$\dim E^{p,q}_2= \begin{array}{c|cc}
	     0 &  0 & 0\\
	     g-3 &  0 & 0\\
	     0 & 3 & 0 \\
	     0 & 1 & 0 \\ \hline
	     0 &  0 & 0
	\end{array}$$
Therefore, in this case, the spectral sequence degenerates at $E_2$ and the desired result follows.

When $X$ is a special GM threefold and $\pi\colon X\to \Gr(2,5)$ is the natural map, we first assume that $\pi(x)$ is in the branch locus of $\pi$. Then Proposition \ref{projection} gives a triangle $K_x[1]\to i^*\oh_x[-1]\to \oh_x$. In this case, we have $\Ext^2(I_x,i^*\oh_x[-1])=\Ext^2(I_x, \oh_x)=k$, thus the dimensions of the second page are of the form

	$$\dim E^{p,q}_2= \begin{array}{c|cc}
	     0 &  0 & 0\\
	     4 &  1 & 0\\
	     0 & 3 & 0 \\
	     0 & 1 & 0 \\ \hline
	     0 &  0 & 0
	\end{array}$$ 
Therefore, the spectral sequence degenerates at $E_2$ as well and we have $$\RHom^\bullet(i^*\oh_x, i^*\oh_x)=k\oplus k^7[-1]\oplus k[-2].$$ When $\pi(x)$ is not in the branch locus, we know that $$\Ext^2(I_x,i^*\oh_x[-1])=\Ext^2(I_x, \oh_{\tau(x)})=0.$$ Then from the above argument, we obtain $\RHom^\bullet(i^*\oh_x, i^*\oh_x)=k\oplus k^6[-1]$.

When $g$ is odd, a similar argument shows that $\Hom(i^*\oh_x, i^*\oh_x)=1$. Since the homological dimension of $\cA_{\sigma}$ is one, $i^*\oh_x$ is the direct sum of shifts of its cohomology objects with respect to $\cA_{\sigma}$. Then from $\mathrm{hom}(i^*\oh_x,i^*\oh_x)=1$, we know that $i^*\oh_x$ is in the heart up to shifts. Hence we have $\Ext^i(i^*\oh_x, i^*\oh_x)=0$ for $i\geq 2$. So the dimension of $\Ext^1$ follows from the Euler characteristics.
\end{proof}

\begin{Lem} \label{coho two}
We have $\cH^j_{\mathrm{Coh}^0_{\alpha, \beta}(X)}(i^*\oh_x)=0$ for $j\neq -2,-1$ and $(\alpha, \beta)$ as in Proposition \ref{induce_range}.
\end{Lem}

\begin{proof}
By definition of $i^*\oh_x[-1]$, we have a triangle $\cE^{\oplus n_g}\to I_x\to i^*\oh_x[-1]$. Since $\cE, I_x\in \mathrm{Coh}^0_{\alpha, \beta}(X)$, the result follows from the long exact sequence of cohomology objects.
\end{proof}

For simplicity, in each case of $g\geq 6$, we define 

\begin{itemize}
    \item $g=6$: $(\alpha_0, \beta_0)=(\frac{1}{20}, -\frac{9}{10})$,
    
    \item $g=7$: $(\alpha_0, \beta_0)=(\frac{1}{12}, -\frac{5}{6})$,
    
    \item $g=8$: $(\alpha_0, \beta_0)=(\frac{1}{25}, -\frac{22}{25})$,
    
    \item $g=9$: $(\alpha_0, \beta_0)=(\frac{1}{8}, -\frac{3}{4})$,
    
    \item $g=10$: $(\alpha_0, \beta_0)=(\frac{1}{25}, -\frac{22}{25})$,
    
    \item $g=12$: $(\alpha_0, \beta_0)=(\frac{1}{25}, -\frac{21}{25})$.
\end{itemize}

\begin{Lem} \label{K in heart}
Let $X$ be a prime Fano threefold of genus $g\geq 6$ and $x\in X$. Then we have $i^*\oh_x[-1]\in \cA(\alpha_0, \beta_0)$.
\end{Lem}

\begin{proof}
Let $\sigma\coloneqq\sigma(\alpha_0, \beta_0)$. If $g\neq 6, 8$, the homological dimension of $\cA_{\sigma}$ is one by Proposition \ref{ogm homo dim 2}. Thus, $i^*\oh_x$ is the direct sum of shifts of its cohomology objects with respect to $\cA_{\sigma}$. Since $\mathrm{hom}(i^*\oh_x,i^*\oh_x)=1$ by Lemma \ref{ext of K}, we know that $i^*\oh_x$ has only one cohomology object, i.e. $i^*\oh_x$ is in the heart $\cA_{\sigma}$ up to a shift.

When $g=6,8$, as in \cite[Lemma 4.5]{bernardara2012categorical}, we consider the spectral sequence for objects in $\mathcal{K}u(X)$ whose second page is given by 
	\[E^{p,q}_2=\bigoplus_i \mathrm{Hom}^p(\mathcal{H}^i_{\sigma}(i^*\oh_x), \mathcal{H}^{i+q}_{\sigma}(i^*\oh_x))\Rightarrow \mathrm{Hom}^{p+q}(i^*\oh_x,i^*\oh_x)\]
where the cohomology is taken with respect to the heart $\cA(\alpha_0, \beta_0)$. Let $r$ be the number of non-zero cohomology objects of $i^*\oh_x$ with respect to the heart. Note that since the homological dimension of $\cA_{\sigma}$ is at most $2$ by Proposition \ref{ogm homo dim 2}, we have $E_2^{1,q}=E^{1,q}_{\infty}$. Therefore, if we take $q=0$, we see that
\[\mathrm{ext}^1(i^*\oh_x,i^*\oh_x)\geq \sum_i \mathrm{ext}^1(\mathcal{H}_{\sigma}^i(i^*\oh_x),\mathcal{H}_{\sigma}^i(i^*\oh_x))\geq 2r.\]
If $r\geq 2$, by Lemma \ref{coho two} we know that $r=2$. Together with Lemma \ref{ext of K}, we obtain
\begin{equation}\label{eq:heart-1}
1-\chi(\cH^{-1}_{\sigma}(i^*\oh_x),\cH^{-1}_{\sigma}(i^*\oh_x))+1-\chi(\cH^{-2}_{\sigma}(i^*\oh_x),\cH^{-2}_{\sigma}(i^*\oh_x)) \leq \delta_g,
\end{equation}
where $\delta_6=7$ and $\delta_8=8$.

Recall that we have  $\cN(\Ku(X))=\langle v, w \rangle$. If we assume that  $[\cH^{-1}_{\sigma}(i^*\oh_x)[1]]=av+bw$ and  $[\cH^{-2}_{\sigma}(i^*\oh_x)[2]]=cv+dw$ for some $a,b,c,d\in \mathbb{Z}$, then from $$[\cH^{-1}_{\sigma}(i^*\oh_x)[1]]+[\cH^{-2}_{\sigma}(i^*\oh_x)[2]]=[i^*\oh_x],$$ we have
\begin{equation}\label{eq:heart-2}
a+c=g-1, \quad b+d=-\frac{g}{2}.
\end{equation}
Moreover, from $\cH^{l}_{\sigma}(i^*\oh_x)\in \cA_{\sigma}$ for any $l$, we see that 
\begin{equation}\label{eq:heart-3}
\Im Z^0_{\alpha_0, \beta_0}(av+bw)\leq 0, \quad \Im Z^0_{\alpha_0, \beta_0}(cv+dw)\geq 0.
\end{equation}
A direct computation shows that there is no $(a,b,c,d)\in \ZZ^4$ satisfying \eqref{eq:heart-1}, \eqref{eq:heart-2}, and \eqref{eq:heart-3} at the same time, and the result follows.
\end{proof}

Now, we prove the main result of this section. We need the following lemma.

\begin{Lem} \label{chi-1}
Let $\cT\subset \Db(X)$ be a triangulated subcategory and $\cA$ be a heart of $\cT$ with homological dimension $\leq 2$. Let $0\to A\to E\to B\to 0$ be an exact sequence in $\cA$. Suppose that $\Hom(B, A[2])=0$ and $\Hom(E,E)=k$. Then we have $\chi(B,A)\leq -1$.
\end{Lem}

\begin{proof}
First note that $\Ext^1(B,A)\neq 0$, otherwise $E\cong A\oplus B$ which contradicts $\Hom(E,E)=k$. 
By assumption, we have $\Hom(B, A[i])=0$ for $i\geq 2$. We also have $\Hom(B,A)=0$, otherwise we have a non-zero composition of maps $E\twoheadrightarrow B\to A\hookrightarrow E$, which contradicts $\Hom(E,E)=k$. Therefore we obtain $\chi(B,A)=-\ext^1(B,A)\leq -1$.
\end{proof}

\begin{Prop} \label{prop:skyscraper projection stability}
Let $X$ be a prime Fano threefold of genus $g\geq 6$ and $x \in X$. Then $i^*\oh_x$ is $\sigma$-stable for every Serre-invariant stability condition $\sigma$ on $\Ku(X)$.
\end{Prop}

\begin{proof} 
By Proposition \ref{induce_range} and Theorem \ref{all_in_one_orbit}, we may assume that $\sigma=\sigma(\alpha_0, \beta_0)$.

When $g\geq 6$, we know that for $(\alpha_0, \beta_0)$ as above, $i^*\oh_x[-1]$ is in the heart $\cA(\alpha_0, \beta_0)$ by Lemma \ref{K in heart}. If $i^*\oh_x$ is not $\sigma$-semistable, we can find an exact sequence
\[0\to A\to i^*\oh_x[-1]\to B\to 0\]
in $\cA(\alpha_0, \beta_0)$ with $\mu_{\sigma}(A)> \mu_{\sigma}(B)$. By Lemma \ref{weak-mukai-lemma}, we have
\[\ext^1(A,A)+\ext^1(B,B)\leq \ext^1(i^*\oh_x,i^*\oh_x) . \]
Note that $1-\chi(A,A)\leq \ext^1(A,A)$. Therefore, we have 
\[[A]+[B]=[i^*\oh_x[-1]],\]
\[\Im Z^0_{\alpha_0, \beta_0}(A)\geq 0,\quad \Im Z^0_{\alpha_0, \beta_0}(B)\geq 0,\]
\[\mu^0_{\alpha_0, \beta_0}(A)>\mu^0_{\alpha_0, \beta_0}(B),\]
and
\[1-\chi(A,A)+1-\chi(B,B)\leq \ext^1(i^*\oh_x,i^*\oh_x).\]

Thus, if we assume $[A]=av+bw$ and $[B]=cv+dw$, by the results in Section \ref{5.9computation}, we have the following integer solutions for $(a,b,c,d)$:

\begin{itemize}
    \item $g=6$ and ordinary:  $(a,b,c,d)=(-2,1,-3,2)$;
    
    \item $g=6$ and special:  $(a,b,c,d)=(-2,1,-3,2)$ or $(a,b,c,d)=(-4,2,-1,1)$;
    
    \item $g=7$: 
    there are no solutions;

    \item $g=8$: $(a,b,c,d)=(-2,1,-5,3)$ or $(a,b,c,d)=(-4,2,-3,2)$;
    
    \item $g=9$: 
    there are no solutions;
    
    \item $g=10$: 
    there are no solutions;
    
    \item $g=12$: 
    there are no solutions.
    
\end{itemize}

Therefore, if $g\neq 6,8$, we know that $i^*\oh_x[-1]$ is $\sigma(\alpha_0, \beta_0)$-semistable. In the following, we consider the cases $g=6$ and $8$ separately.

\begin{itemize}
    \item Assume that $X$ is an ordinary GM threefold. The solution $(a,b)=(-2,1)$ gives a potential destabilizing subobject $A$ of $i^*\oh_x[-1]$.  It is known that $A\cong i^*(E)[1]$ for some $\mu_H$-stable sheaf $E$ with $\mu_H(E)=-\frac{1}{2}$ by \cite[Theorem 8.9]{jacovskis2024categorical}. By Lemma \ref{Kstable} and $\mu_H$-stability, we have $\Hom(E, K_x)=0$. Since $K_x=i^*\oh_x[-2]$, the adjointness implies $$\Hom(i^*(E)[1], i^*\oh_x[-1])=\Hom(E[1],i^*\oh_x[-1])=\Hom(E, K_x)=0,$$ which gives a contradiction. Therefore, $i^*\oh_x[-1]$ is $\sigma(\alpha_0, \beta_0)$-semistable.

    \item Assume that $X$ is a special GM threefold. If $(a,b)=(-2,1)$, then we also have $A\cong i^*(E)[1]$ for some $\mu_H$-stable sheaf $E$ with $\mu_H(E)=-\frac{1}{2}$, hence $\Hom(E, K_x)=0$. If we apply $\Hom(E,-)$ to the triangle $K_x[1]\to i^*\oh_x[-1]\to \oh_{\tau(x)}$, from $$\Hom(E,K_x)=\Hom(E, \oh_{\tau(x)}[-1])=0,$$ we obtain $\Hom(E,i^*\oh_x[-2])=\Hom(A, i^*(\oh_{x})[-1])=0$, which is a contradiction. If $(a,b)=(-4,2)$, then we have $\chi(B,A)=0$, which contradicts Lemma \ref{chi-1}. Therefore, $i^*\oh_x[-1]$ is $\sigma(\alpha_0, \beta_0)$-semistable.

    \item Assume that $g=8$. If $(a,b)=(-2,1)$, then by \cite[Theorem 1.1]{liu2021note} we know that $A\cong i^*(E)[1]$ for some $\mu_H$-stable sheaf $E$ with $\mu_H(E)=-\frac{1}{2}$. By Lemma \ref{Kstable} and stability, we have $\Hom(E, K_x)=0$. Since $K_x=i^*\oh_x[-2]$, by adjointness, we have $$\Hom(i^*(E)[1], i^*\oh_x[-1])=\Hom(E[1], i^*\oh_x[-1])=\Hom(E, K_x)=0,$$ which gives a contradiction. If $(a,b)=(-4,2)$, then we have $\chi(B,A)=0$, which contradicts  Lemma \ref{chi-1}. Therefore, $i^*\oh_x[-1]$ is $\sigma(\alpha_0, \beta_0)$-semistable.
\end{itemize}

It remains to show the $\sigma(\alpha_0, \beta_0)$-stability of $i^*\oh_x[-1]$. When $g\neq 7$, $[i^*\oh_x[-1]]$ is a primitive class, hence $i^*\oh_x[-1]$ is $\sigma$-stable. 

When $g=7$, the only possible Jordan--H\"older factors of $i^*\oh_x[-1]$ have numerical class $-6v+5w$. By \cite{kuznetsov2005derived}, there is an equivalence $\Theta\colon \Ku(X)\to \Db(\Gamma_7)$, where $\Gamma_7$ is a smooth projective curve of genus $7$. Moreover, up to some autoequivalences, we have $\Theta(i^*\oh_x[-1])\cong E$, where $E$ is a rank $2$ semistable bundle on $\Gamma_7$ of degree $12$ and $h^0(E)=5$ (see e.g. \cite[Theorem 2.3]{brambilla2014vector}). Therefore, if $i^*\oh_x[-1]$ is strictly $\sigma$-semistable, $\Theta$ maps the Jordan--H\"older filtration of $i^*\oh_x[-1]$ to an exact sequence
$0\to L_1\to E\to L_2\to 0$ on $\Gamma_7$, where $L_i$ are line bundles of degree $6$. By Riemann--Roch, we know that $h^0(L_i)=h^1(L_i)$ for each $i$. Therefore, there exists $j\in \{1,2\}$ such that $h^0(L_j)=4$ or $3$. By \cite[Table 1]{mukai95}, we know that the curve $\Gamma_7$ is not trigonal nor tetragonal, hence the only possibility is $h^0(L_j)=4$. However, by Clifford's theorem, this implies $\Gamma_7$ is hyperelliptic, contradicting \cite[Theorem 8.1]{mukai:nonabelian}.

Thus, the above argument shows that for every $g\geq 6$, the object $i^*\oh_x[-1]$ is $\sigma$-stable.
\end{proof}

\subsection{Embedding to a Bridgeland moduli space} \label{subsec:embedding}

Finally, we can embed $X$ in a Bridgeland moduli space.

Let $\sigma$ be a Serre-invariant stability condition on $\Ku(X)$. By \cite[Corollary 26.2]{bayer2021stability}, for any $0\neq \xi\in \cN(\Ku(X))$, there is a proper moduli space $M_{\sigma}(\xi)$, whose $k$-points one-to-one correspond to $\sigma$-stable objects of class $\xi$ in $\cA_{\sigma}$.

\begin{Prop}\label{closed_embedding_into_modulispace}
Let $X$ be a prime Fano threefold of genus $g\geq 6$ and $\sigma$ be a Serre-invariant stability condition on $\Ku(X)$. Then up to a unique shift of $\sigma$, the functor $i^*$ induces a closed embedding
\[p\colon X\hookrightarrow M_{\sigma}([i^*\oh_x]), \quad x\mapsto [i^*\oh_x].\]
\end{Prop}

\begin{proof}
By Proposition \ref{prop:skyscraper projection stability}, we may assume that $i^*\oh_x\in \cA_{\sigma}$ for any $x\in X$. Then we get a morphism $p\colon X\to M_{\sigma}([i^*\oh_x])$.

We first show that $p$ is injective. Indeed, if $g\geq 7$ or if $X$ is an ordinary GM threefold, $i^*\oh_x[-2]$ is a reflexive sheaf that is not locally free only at the point $x$ by Proposition \ref{projection}. Therefore, we see that $i^*\oh_x\cong i^*\oh_y$ if and only if $x=y$. If $X$ is a special GM threefold, then $\cH^0(i^*\oh_x[-1])\cong \cH^0(i^*\oh_y[-1])$ if and only if $\tau(x)=\tau(y)$ by Proposition \ref{projection}, which is also equivalent to $x=y$.

Since $M_{\sigma}([i^*\oh_x])$ is proper, to show $p$ is a closed embedding, it remains to prove that the tangent map of $p$ at each point is injective. For any $x\in X$, note that the tangent map of $p$ at $x$ is given by the natural map
\[\Ext^1(\oh_x, \oh_x)\to \Ext^1(i^*\oh_x,i^*\oh_x).\]
Moreover, the functor $\bL_{\oh_X}$ induces an isomorphism $\Ext^1(\oh_x, \oh_x)\cong \Ext^1(I_x, I_x)$. Thus, we only need to show that $d_2\colon \Ext^1(I_x, I_x)\to \Ext^1(\bL_{\cE} I_x, \bL_{\cE} I_x)$ is injective.

We consider the defining triangle of $\bL_{\cE} I_x$ and apply $\Hom(I_x, -)$ to it. Inside the long exact sequence, we have the sequence
\[ \cdots \ra \Ext^1(I_x, \cE^{\oplus n_g}) \ra \Ext^1(I_x, I_x) \xra{d_2} \Ext^1(I_x, \bL_{\cE} I_x)=\Ext^1(\bL_{\cE} I_x, \bL_{\cE} I_x) \ra \cdots . \]
Since $\Ext^1(I_x, \cE)=0$, we can see from the sequence above that the map $d_2$ is indeed injective.
\end{proof}

Using the above description, we can verify Mukai's conjecture on dual embeddings.

\begin{Thm}\label{thm:normal-bundle}
Let $X$ be a prime Fano threefold of index one and genus $g\geq 6$. Assume that either $g\geq 7$, or $g=6$ and $X$ is ordinary. Then we have an isomorphism
\[N^{\vee}_{X/\mathbf{G}_X}(H)\cong N_{X/M_\sigma([i^*\oh_x])}.\]
\end{Thm}

\begin{proof}
By our assumption, we have a closed embedding $\gamma\colon X\hookrightarrow \mathbf{G}_X$. Set $M\coloneqq M_\sigma([i^*\oh_x]).$ Let
\[
\pi_1,\pi_2\colon X\times X\to X
\]
be the two projections, and let $\Delta\subset X\times X$ be the diagonal. We write \(\cI_\Delta\) for its ideal sheaf. Consider the vector bundle $\cF\coloneqq\pi_1^*\cQ^\vee\otimes \pi_2^*\cE$ on \(X\times X\). By \cite{kuznetsov:base-change}, we have a relative semiorthogonal decomposition
\[\Db(X\times X)=\langle \Db(X)\boxtimes \Ku(X), \Db(X)\boxtimes \cE, \Db(X)\boxtimes \oh_X \rangle \]
over the first factor. The relative left mutation of $\cI_\Delta$ through \(\cE\) is
defined by the functorial triangle
\[
\pi_1^*\pi_{1*}\mathcal Hom(\pi_2^*\cE,\cI_\Delta)\otimes \pi_2^*\cE
\longrightarrow
\cI_\Delta
\longrightarrow
\bL_{\cE}^{\mathrm{rel}}(\cI_\Delta).\]
Using the exact sequence $0\to \cI_{\Delta}\to \oh_{X\times X}\to \oh_{\Delta}\to 0$ and the tautological sequence $$0\to \cQ^{\vee}\to H^0(\cE^{\vee})\otimes \oh_X\to \cE^{\vee}\to 0,$$ we obtain a canonical isomorphism $\pi_{1*}(\pi_{2}^*\cE^{\vee}\otimes \cI_{\Delta})\cong \cQ^{\vee}$. Therefore, the relative mutation triangle becomes
\begin{equation}\label{eq:relative-mutation}
\cF\xra{\epsilon} \cI_{\Delta}\to \bL_{\cE}^{\mathrm{rel}}(\cI_\Delta)=:\cP.
\end{equation}
By Proposition \ref{projection} and our assumption, we see that $\epsilon$ is surjective and $\cP[-1]$ is a sheaf with $(\cP[-1])|_{\{x\}\times X}\cong K_x$ for each $x\in X$.

We write
\[
\mathcal Ext^k_{\pi_1}(-,-)
\coloneqq
\mathcal H^k\!\left(\pi_{1*}\mathcal Hom(-,-)\right).
\]
By Lemma \ref{ext of K} and our assumption, we see that $M$ is smooth along $p(X)$, and
\[p^*T_M
\simeq
\mathcal Ext^1_{\pi_1}(\cP,\cP).\]
On the other hand, the universal ideal sheaf of the diagonal gives
\[
T_X
\simeq
\mathcal Ext^1_{\pi_1}(\cI_\Delta,\cI_\Delta).\]
Note that $\pi_{1*}\mathcal Hom(\cF,\cP)\cong 0$. Therefore, under the above identifications and the isomorphism $\pi_{1*}\mathcal Hom(\cI_{\Delta},\cP)\cong \pi_{1*}\mathcal Hom(\cP,\cP)$, the tangent map 
\[dp\colon T_X\to p^*T_M\]
is induced by applying $\pi_{1*}\mathcal Hom(\cI_{\Delta},-)$ to the triangle \eqref{eq:relative-mutation}. Moreover, since $\Ext^1(I_x, \cE)=0$ for any $x\in X$, we have $\mathcal Ext^1_{\pi_1}(\cI_\Delta,\cF)=0$. Thus, the long exact sequence gives
\[
0
\longrightarrow
\mathcal Ext^1_{\pi_1}(\cI_\Delta,\cI_\Delta)
\longrightarrow
\mathcal Ext^1_{\pi_1}(\cI_\Delta,\cP)
\longrightarrow
\ker(\rho)
\longrightarrow 0,
\]
where
\[
\rho\colon
\mathcal Ext^2_{\pi_1}(\cI_\Delta,\cF)
\longrightarrow
\mathcal Ext^2_{\pi_1}(\cI_\Delta,\cI_\Delta)
\]
is induced by $\epsilon$. In particular, we get $N_{X/M}\cong \ker(\rho)$.

It remains to identify \(\ker(\rho)\). To this end, by definition, we have
\[\pi_{1*}\cH om (\cI_{\Delta}, \cF)\cong \cQ^{\vee}\otimes \pi_{1*}\cH om (\cI_{\Delta}, \pi_2^*\cE).\]
Using the exact sequence
\[
0\to \cI_\Delta\to \oh_{X\times X}\to \oh_\Delta\to 0
\]
and relative duality for the diagonal, we get
\[
\pi_{1*}\mathcal Hom(\cI_\Delta,\pi_2^*\cE)
\cong
\cE(H)[-2]
\]
and hence
\[\pi_{1*}\cH om (\cI_{\Delta}, \cF)[2]\cong \mathcal Ext^2_{\pi_1}(\mathcal I_\Delta,\cF)
\cong
\Omega^1_{\mathbf{G}_X}|_X(H).\]
Moreover, we have a natural isomorphism $\mathcal Ext^2_{\pi_1}(\cI_\Delta,\cI_\Delta)\cong \wedge^2 T_X\cong \Omega_X(H)$. So under the above natural identifications, the map $\rho$ becomes
\[\rho\colon \Omega^1_{\mathbf{G}_X}|_X(H)\to \Omega_X(H).\]
Therefore, it is straightforward to check that $\rho=(d\gamma)^{\vee}\otimes \oh_X(H)$, and the result follows.
\end{proof}

\section{Fano threefolds as Brill--Noether loci} \label{sec:BN}

In this section, we exhibit $X$ as a Brill--Noether locus inside the moduli space $M_{\sigma}([i^*\oh_x])$. We first fix some notation. Recall that  $n_g\coloneqq\frac{g}{2}$ if $g$ is even, $n_7=5$ and $n_9=3$. Let $x \in X$ be a point, and let $$i_{\cD} \colon \Ku(X) \hookrightarrow \cD=\langle \Ku(X), \cE \rangle$$ be the inclusion; this has a left adjoint $i_{\cD}^*$ and a right adjoint $i_{\cD}^!$. In particular, $i_{\cD}^* = \bL_{\cE}$. Therefore, we can rewrite the triangle \ref{eq:even genus projection original triangle} as 
\begin{equation*}
    \cE^{\oplus n_g} \ra I_x \ra i_{\cD}^* I_x \ra \cE^{\oplus n_g}[1] .
\end{equation*}

\begin{Rem}\label{inside_BN_locus}
We have 
\begin{equation*}
    \Ext^1_{\cD}(i_{\cD}^*I_x, \cE) = \Ext^1_{\cD}(i_{\cD}i_{\cD}^*I_x, \cE) = \Ext^1_{\Ku(X)}(i_{\cD}^*I_x, i^!\cE)=\Hom(i^*\oh_x[-1],i^!\cE[1])
\end{equation*}
because $i_{\cD}^*$ and $i_{\cD}^!$ are left and right adjoints to $i_{\cD}$, respectively.
\end{Rem}

We define 

\begin{itemize}
    \item $(\alpha_6,\beta_6)=(\frac{1}{20},-\frac{9}{10})$, 
    
    \item $(\alpha_7, \beta_7)=(\frac{\sqrt{71}}{84},-\frac{71}{84})$,
    
    \item $(\alpha_8, \beta_8)=(\frac{2\sqrt{79}}{875},-\frac{122}{125})$, 
    
    \item $(\alpha_9, \beta_9)=(\frac{\sqrt{31}}{40},-\frac{31}{40})$, 
    
    \item $(\alpha_{10},\beta_{10})=(\frac{1}{33}\sqrt{\frac{5}{3}},-\frac{10}{11})$, and
    
    \item $(\alpha_{12}, \beta_{12})=(\frac{1}{22},-\frac{19}{22})$.
\end{itemize}

Note that $\mu^0_{\alpha_g,\beta_g}([i^*\oh_x])=+\infty$ when $g\neq 6$.

From Proposition \ref{induce_range} we know that $\sigma(\alpha_g, \beta_g)$ is a Serre-invariant stability condition on $\Ku(X)$. In this section, we fix $\sigma\coloneqq\sigma(\alpha_g, \beta_g)$. Then $i^*\oh_x[-1]\in \cA_{\sigma}$.

\begin{Def}
We define the \emph{Brill--Noether locus} in $M_{\sigma}(-[i^*\oh_x])$ with respect to the object $i^!(\cE)$ as
\begin{equation*}
    \mathsf{BN}_g \coloneqq \{ F \mid \ext^1(F, i^!(\cE)) = n_g \} \sst M_{\sigma}(-[i^*\oh_x]) .
\end{equation*}
\end{Def}

The main theorem we prove in this section is

\begin{Thm} \label{thm: BN = X}
The natural morphism $p\colon X\hookrightarrow M_{\sigma}(-[i^*\oh_x])$ induces an isomorphism
\[X\cong \mathsf{BN}_g. \]
\end{Thm}

\subsection{Proof of Theorem \ref{thm: BN = X}} \label{sec_proof}

We split the proof of the theorem into a series of propositions and lemmas. First, we show that $X$ is embedded in the Brill--Noether locus $\mathsf{BN}_g$. Then we show that each object in $\mathsf{BN}_g$ is of the form $i^*(\oh_x[-1])$. The strategy will be to take an object $F \in \mathsf{BN}_g \setminus p(X)$, and show that it cannot exist.

Let $X$ be a prime Fano threefold of index one and genus $g\geq 6$.

\begin{Lem}
\label{X_embedded_BN}
There is a closed embedding $p\colon X\hookrightarrow\mathsf{BN}_g$ induced by $i^*$. 
\end{Lem}

\begin{proof}
By Proposition~\ref{closed_embedding_into_modulispace}, the projection functor $i^*$ already induces a closed embedding 
\[ p\colon X\hookrightarrow M_{\sigma}(-[i^*\oh_x]) . \] 
It suffices to check that $i^*\oh_x[-1]$ satisfies the Brill--Noether condition, i.e. that $$\mathrm{ext}^1(i^*\oh_x[-1], i^!\cE)=n_g$$ for $g\geq 6$. By Proposition~\ref{projection}, the triangle above becomes
$$\cE^{\oplus n_g}\rightarrow I_x\rightarrow\bL_{\cE}I_x.$$
By Remark~\ref{inside_BN_locus}, $\mathrm{Ext}^1(i^*\oh_x[-1],i^!\cE)\cong\mathrm{Ext}^1(\bL_{\cE}I_x,\cE)$. Since $i_{\cD}^*I_x=\bL_{\cE}I_x$, applying $\mathrm{Hom}(-,\cE)$ to the triangle above, we get a long exact sequence
$$0\rightarrow\mathrm{Hom}(i_{\cD}^*I_x,\cE)\rightarrow\mathrm{Hom}(I_x,\cE)\rightarrow\mathrm{Hom}(\cE,\cE^{\oplus n_g})\rightarrow\mathrm{Ext}^1(i_{\cD}^*I_x,\cE)\rightarrow\mathrm{Ext}^1(I_x,\cE)\rightarrow\cdots$$
It is easy to check $\mathrm{Hom}(I_x,\cE)=\mathrm{Ext}^1(I_x,\cE)=0$, so $\mathrm{Ext}^1(i^*I_x,\cE)\cong\mathrm{Hom}(\cE,\cE^{\oplus n_g})\cong k^{n_g}$. 
Then the desired result follows. 
\end{proof}

\begin{Lem} \label{prop:F is shift of vb}
Let $F \in \mathsf{BN}_g \setminus p(X)$. Then $F[-1]$ is a vector bundle.
\end{Lem}

\begin{proof}
By Serre duality, for each $j\in \ZZ$, we have 
\begin{align*}
    \Ext^j(F, \oh_x) \cong \Ext^{2-j}(i^*\oh_x[-1],F)^\vee \cong \Ext^{-j}(F[2], S_{\Ku(X)}(i^*\oh_x[-1])) .
\end{align*}
Since $F$ and $i^*\oh_x[-1]$ are in the same heart with homological dimension at most two by Lemma \ref{ogm homo dim 2}, we have $\Ext^i(F, \oh_x)=0$ for $2<i$ and $i<0$. Since $F\notin p(X)$ by assumption, using $\sigma$-stability of $F$ and $i^*\oh_x[-1]$, we get $\Ext^2(F, \oh_x)=\Hom(i^*\oh_x[-1],F)^\vee=0$. Thus, $\Ext^j(F, \oh_x)=0$ for $j\notin \{0,1\}$.

First, we assume that $g\neq 6$.  Note that $S_{\Ku(X)}(i^*\oh_x[-1])$ is stable in the Kuznetsov component by Serre-invariance of $\sigma$. For $g=8$, we have \[ \Hom(F, \oh_x)=\Ext^2(i^*\oh_x[-1], F) = 0 \] by \cite[Corollary 5.5, Lemma 5.9]{PY20}, and  \[ \Hom(F, \oh_x)=\Ext^2(i^*\oh_x[-1], F) = 0 \] for $g=7,9,10, 12$ since the homological dimensions of $\cA_{\sigma}$ in these cases are one. This means $\RHom^\bullet(F, \oh_x)=\Ext^1(F, \oh_x)[-1]$.  By \cite[Proposition 5.4]{bm99}, $F[-1]$ is a vector bundle.

Now assume that $g=6$. If $X$ is a special GM threefold, then 
\begin{align*}
&\mathrm{Hom}(F,\oh_x) \cong \mathrm{Ext}^2(i^*\oh_x[-1],F)^{\vee} \\ \cong&\mathrm{Hom}(F,\tau(i^*\oh_x[-1])) \cong\mathrm{Hom}(F,i^*(\oh_{\tau(x)})[-1]).
\end{align*} 
By our assumption, we see that $\mathrm{Hom}_{\Db(X)}(F,\oh_x)=0$ and the desired result again follows from \cite[Proposition 5.4]{bm99}.

Finally, if $X$ is an ordinary GM threefold, then \begin{align*}\mathrm{Hom}(F,\oh_x) \cong\mathrm{Ext}^2(i^*\oh_x[-1],F)^{\vee} \cong\mathrm{Hom}(F,\tau(i^*\oh_x[-1])) . \end{align*} We claim that $\mathrm{Hom}(F,\tau(i^*\oh_x[-1]))=0$. Indeed, if $\mathrm{Hom}(F,\tau(i^*\oh_x[-1]))\neq 0$, then we have $F\cong\tau(i^*\oh_x[-1])$ since $F$ and $\tau(i^*\oh_x[-1])$ are both $\sigma$-stable objects of the same phase in $\Ku(X)$. But then this means that $\tau(i^*\oh_x[-1])$ is also in the Brill--Noether locus $\mathsf{BN}_g$, i.e. $\mathrm{Ext}^1(\tau(i^*\oh_x[-1]),i^!\cE)=k^3$, which is impossible by Lemma~\ref{Lemma_not_in_BN} below. This completes the proof.
\end{proof}

\begin{Lem}
\label{Lemma_not_in_BN}
We have $\mathrm{Ext}^1(\tau(i^*\oh_x[-1]),i^!\cE)\neq k^3$. 
\end{Lem}

\begin{proof}
By Remark~\ref{gluing_data_equal}, the gluing object $i^!\cE$ is given by the exact triangle $$\mathcal{Q}(-H)[1]\rightarrow  i^!\cE\rightarrow\cE.$$
Since $\tau^{-1}\cong i^*\circ(-\otimes\oh_X(H))[-1]$ by Lemma \ref{lem:serre-sod}, we have $\tau^{-1}(i^!\cE)\cong\bL_{\cE}\cQ^{\vee}$. Thus $$\mathrm{Ext}^1(\tau(i^*\oh_x[-1]),i^!\cE)\cong\mathrm{Ext}^2(i^*\oh_x, \bL_{\cE}\cQ^{\vee})\cong\mathrm{Ext}^2(\oh_x,\bL_{\cE}\cQ^{\vee}).$$ On the other hand, by \cite[Lemma 5.1]{jacovskis2024categorical}, there is a triangle
$$\cE^{\oplus 2}\rightarrow\cQ^{\vee}\rightarrow\bL_{\cE}\cQ^{\vee}.$$
Applying $\mathrm{Hom}(\oh_x,-)$ to it, we get the following long exact sequence:
$$0\rightarrow\mathrm{Ext}^2(\oh_x,\bL_{\cE}\cQ^{\vee})\rightarrow k^4\rightarrow k^3\rightarrow\mathrm{Ext}^3(\oh_x,\bL_{\cE}\cQ^{\vee})\rightarrow 0.$$
Thus $\mathrm{Ext}^1(\tau(i^*\oh_x[-1]),i^!\cE)=k^3$ if and only if $\mathrm{Ext}^3(\oh_x,\bL_{\cE}\cQ^{\vee})=k^2$. By Serre duality, we get $\mathrm{Ext}^3(\oh_x,\bL_{\cE}\cQ^{\vee})\cong\mathrm{Hom}(\bL_{\cE}\cQ^{\vee},\oh_x)$.  Since $X$ is an ordinary GM threefold, by \cite[Lemma 5.1]{jacovskis2024categorical}, we have $\mathrm{Hom}(\cE,\cQ^{\vee})=k^2$ and there is a short exact sequence 
$$0\rightarrow\cE\rightarrow\cQ^{\vee}\rightarrow I_C\rightarrow 0$$ for some conic $C\subset X$. Then \cite[Lemma 7.16]{jacovskis2024categorical} gives $$\bL_{\cE}\cQ^{\vee}\cong\bL_{\cE}I_C\cong\mathbb{D}(I_C)\otimes\oh_X(-H)[1],$$ where $\mathbb{D}(-)$ is the derived dual functor.
Therefore, we see that $$\mathrm{Hom}(\bL_{\cE}\cQ^{\vee},\oh_x)\cong\mathrm{Hom}(\mathbb{D}(I_C)\otimes\oh_X(-H)[1],\oh_x)\cong\mathrm{Hom}(\mathbb{D}(I_C)[1],\oh_x).$$ 
Note that $\mathbb{D}(I_C)[1]$ is given by the exact triangle
$$\oh_X[1]\rightarrow\mathbb{D}(I_C)[1]\rightarrow\oh_C.$$
Applying $\mathrm{Hom}(-,\oh_x)$, we get $\mathrm{Hom}(\oh_C,\oh_x)\cong\mathrm{Hom}(\mathbb{D}(I_C)[1],\oh_x)$. An easy computation shows that $\mathrm{hom}(\oh_C,\oh_x)$ is either $0$ or $1$, and the result follows.
\end{proof}

Let $F\in \mathsf{BN}_g\setminus p(X)$. By assumption, we have a natural map $F\to \cE^{\oplus n_g}[1]$. Since $F$ and $\cE\in \Coh^0_{\alpha_g, \beta_g}(X)$, the extension $C$ of $F$ and $\cE^{\oplus n_g}$
\[\cE^{\oplus n_g}\to C\to F\]
is also in the heart $\Coh^0_{\alpha_g, \beta_g}(X)$. Note that $\ch(C)=\ch(I_x)$. If we apply $\Hom(-, \cE)$ to this exact sequence, the natural map $\Hom(\cE^{\oplus n_g}, \cE)\to \Ext^1(F, \cE)$ is bijective by construction. Since $\Ext^1(\cE^{\oplus n_g}, \cE)=0$, we have
\begin{equation} \label{CE=0}
    \Hom(C, \cE)=\Ext^1(C, \cE)=0.
\end{equation}

\begin{Prop} \label{gneq6cone}
The object $C\in \Coh^0_{\alpha_g, \beta_g}(X)$ is $\sigma^0_{\alpha_g, \beta_g}$-semistable.
\end{Prop}

\begin{proof}

We assume that $C\in \Coh^0_{\alpha_g, \beta_g}(X)$ is not $\sigma^0_{\alpha_g, \beta_g}$-semistable. Let $B$ be the minimal destabilizing quotient object of $C$. Then we have an exact sequence in $\Coh^0_{\alpha_g, \beta_g}(X)$
\[0\to A\to C\to B\to 0\]
where
\begin{equation}\label{eq:abc-slope}
\mu^{0 -}_{\alpha_g, \beta_g}(A)>\mu^0_{\alpha_g, \beta_g}(C)>\mu^0_{\alpha_g, \beta_g}(B)
\end{equation}
and $B$ is $\sigma^0_{\alpha_g, \beta_g}$-semistable. Therefore, we have inequalities
\[\Im(Z^0_{\alpha_g, \beta_g}(A))\geq 0,\quad \Im(Z^0_{\alpha_g, \beta_g}(B))> 0\]
and
\[\mu^0_{\alpha_g, \beta_g}(C)>\mu^0_{\alpha_g, \beta_g}(B).\]

We have $\mu^0_{\alpha_g, \beta_g}(C)=\mu^0_{\alpha_g, \beta_g}(\oh_X)<0$. Since $B\in \Coh^0_{\alpha_g, \beta_g}(X)$ is $\sigma^0_{\alpha_g, \beta_g}$-semistable with slope $\mu^0_{\alpha_g, \beta_g}(B)<\mu^0_{\alpha_g, \beta_g}(C)<0$ and $F[-1]$ is a bundle by Lemma \ref{prop:F is shift of vb}, we know that $B\in \Coh^{\beta_g}(X)$ is $\sigma_{\alpha_g, \beta_g}$-semistable by \cite[Proposition 4.1]{FeyzbakhshPertusi2021stab}.

Note that $\oh_X, \oh_X(-H)[2], \cE, \cE(-H)[2]\in \Coh^0_{\alpha_g, \beta_g}(X)$ and they are $\sigma^0_{\alpha_g, \beta_g}$-semistable. Then by Serre duality and the definition of heart, we know $\Hom(\oh_X, A[n])=\Hom(\cE, A[n])=0$ for every $n<0$ and $n\geq 2$ and the same holds for $B$. Also, from $\mu^0_{\alpha_g, \beta_g}(C)=\mu^0_{\alpha_g, \beta_g}(\oh_X)>\mu^0_{\alpha_g, \beta_g}(B)$ we have $\Hom(\oh_X, B)=0$. Therefore, if we apply $\Hom(\oh_X, -)$ and $\Hom(\cE, -)$ to the triangle $A\to C\to B$, from $\RHom^\bullet(\oh_X, C)=0$ and $\RHom^\bullet(\cE, C)=k^{n_g}$ we obtain $A,B\in \oh^{\perp}_X$, $$\RHom^\bullet(\cE, B)=\Hom(\cE, B),\quad \RHom^\bullet(\cE, A)=\Hom(\cE, A)\oplus \Hom(\cE, A[1])[-1],$$ and a long exact sequence
\[0\to \Hom(\cE, A)\to \Hom(\cE, C)=k^{n_g}\to \Hom(\cE, B)\to \Ext^1(\cE, A)\to 0.\]

Recall that in \eqref{v_w} we have  $\cN(\Ku(X))=\langle v, w \rangle$. If we assume that $[B]=av+bw+c[\cE]$ for $a,b,c\in \mathbb{Z}$, from $\chi(\cE, B)=\hom(\cE, B)\geq 0$ we have $c=\chi(\cE, B)\geq 0$. 

If we apply the projection functor $i^*$ to the triangle, since $i^*(C)\cong F$, we obtain a triangle
\begin{equation} \label{prABC}
    i^*(B)[-1]\to i^*(A)\to F.
\end{equation}

By the definition of projection functor, we get a triangle:

\begin{equation} \label{prB}
    i^*(B)[-1]\to \cE^{\oplus c}\to B.
\end{equation}

\medskip

\noindent\textbf{Claim 1}: $i^*(A)\in \Coh^0_{\alpha_g, \beta_g}(X)$.

Let $t\coloneqq\ext^1(\cE, A)$. Then $\hom(\cE, A)=n_g-c+t$. Therefore we have a triangle $$\cE^{\oplus n_g-c+t}\oplus \cE^{\oplus t}[-1]\to A\to i^*(A)$$ and a long exact sequence of cohomologies in the heart $\Coh^0_{\alpha_g, \beta_g}(X)$:

\begin{equation} 
    0\to \cH^{-1}_{\Coh^0_{\alpha_g, \beta_g}(X)}(i^*(A))\to \cE^{\oplus n_g-c+t}\to A\to \cH^0_{\Coh^0_{\alpha_g, \beta_g}(X)}(i^*(A)) \to \cE^{\oplus t}\to 0.
\end{equation}

Since $\RHom^{\bullet}(\cE, C)=k^{n_g}$ and the natural morphism $\cE^{\oplus n_g}\to C$ is injective in $\Coh^0_{\alpha_g, \beta_g}(X)$, from $\Hom(\cE, A)\subset \Hom(\cE, C)$, we know that the natural map $\cE^{\oplus n_g-c+t}\to A$ is injective as well. This implies $\cH^0_{\Coh^0_{\alpha_g, \beta_g}(X)}(i^*(A))\cong i^*(A)\in \Coh^0_{\alpha_g, \beta_g}(X)$ and we have an exact sequence in $\Coh^0_{\alpha_g, \beta_g}(X)$:

\begin{equation} \label{prA}
    0\to \cE^{\oplus n_g-c+t}\to A\to i^*(A) \to \cE^{\oplus t}\to 0.
\end{equation}

\medskip

\noindent\textbf{Claim 2}: $i^*(B)[-1]\in \Coh(X)$ is torsion-free. Thus $a<0$ and $\mu_H^+(i^*(B)[-1])< \mu_H(\cE)$.

From Proposition \ref{prop:F is shift of vb} and the construction of $C$, we have $\cH^i(C)=0$ for $i\neq -1,0$. Since $B\in \Coh_{\alpha_g,\beta_g}(X)$, we know that $\cH^i(B)=0$ for $i\neq -1,0$.
Thus, if we take the cohomology long exact sequence associated to \eqref{prB}
with respect to the standard heart, we have $\cH^i(i^*(B)[-1])=0$ for $i\notin \{0,1\}$.

Next, as $i^*(A)\in \Coh^0_{\alpha_g, \beta_g}(X)$ by Claim 1, we obtain $\cH^{j}(i^*(A))=0$ for any $j\geq 1$. Hence, if we take the cohomology long exact sequence associated to (\ref{prABC}) 
with respect to the standard heart, since $F[-1]\in \Coh(X)$, we have $\cH^{1}(i^*(B)[-1])=0$, i.e. $i^*(B)[-1]\cong \cH^{0}(i^*(B)[-1])\in \Coh(X)$. Therefore, we have a long exact sequence in $\Coh(X)$:
\[0\to \cH^{-1}(B)\to i^*(B)[-1] \xra{\theta} \cE^{\oplus c} \to \cH^0(B)\to 0\]
and $i^*(B)[-1]\in \Coh(X)$ is a torsion-free sheaf. Moreover, by $[B]=av+bw+c[\cE]$, we have $[i^*(B)]=av+bw$. So $a\leq 0$. And if $a=0$, then $i^*(B)=0$ from torsion-freeness, which implies $B\cong \cE^{\oplus c}$ and contradicts \eqref{CE=0}. Thus, we always have $a<0$. By $\mu_H(\cE)>\beta_g\geq \mu_H(\cH^{-1}(B))$, we also get $\mu_H^+(i^*(B)[-1])\leq \mu_H(\cE)$. Note that when $\mu_H^+(i^*(B)[-1])=\mu_H(\cE)$, $i^*(B)[-1]$ contains a $\mu_H$-stable subsheaf of slope $\mu_H(\cE)$, which is isomorphic to $\cE$ by the uniqueness of Jordan--H\"older factors and contradicts $i^*(B)\in \Ku(X)$. Therefore, we obtain $\mu_H^+(i^*(B)[-1])< \mu_H(\cE)$.

\medskip

\noindent\textbf{Claim 3}: $c>0$ and $\mu^0_{\alpha_g, \beta_g}(B)\geq \mu^0_{\alpha_g, \beta_g}(\cE)$.

We claim that $c=\hom(\E, B)\neq 0$, which implies $\mu^0_{\alpha_g, \beta_g}(B)\geq \mu^0_{\alpha_g, \beta_g}(\cE)$ by $\sigma^0_{\alpha_g, \beta_g}$-stability of $B$ and $\cE$. Indeed, if $c=0$, then we have $B\in \Ku(X)$ and we obtain an exact sequence in the heart $\cA(\alpha_g, \beta_g)$
\[i^*(A)\to F\to B .\]
Since $F$ is $\sigma(\alpha_g, \beta_g)$-stable, we have $\mu^0_{\alpha_g, \beta_g}(B)>\mu^0_{\alpha_g, \beta_g}(F)$, which gives a contradiction since $\mu^0_{\alpha_g, \beta_g}(F)>\mu^0_{\alpha_g, \beta_g}(\oh_X)>\mu^0_{\alpha_g, \beta_g}(B)$.

\bigskip

Now, we are ready to prove our main statement. We divide the argument into several cases.

\medskip

\noindent\textbf{Case 1}: $g\neq 6$.

By Claim 3, we have a triangle $\cE^{\oplus c}\xra{\lambda} B\to i^*(B)$. Therefore, $\cH^{-1}_{\Coh^0_{\alpha_g, \beta_g}(X)}(i^*(B))\cong \ker(\lambda)$ and $\cH^0_{\Coh^0_{\alpha_g, \beta_g}(X)}(i^*(B))\cong \cok(\lambda)$. Note that $$\cH^{j}_{\Coh^0_{\alpha_g, \beta_g}(X)}(i^*(B))\in \cA(\alpha_g, \beta_g),\quad \forall j\in \ZZ.$$ Taking the cohomology long exact sequence of \eqref{prABC} with respect to the heart $\cA(\alpha_g, \beta_g)$, we have an exact sequence in $\cA(\alpha_g, \beta_g)$
\[0\to \ker(\lambda)\to i^*(A)\to F\to \cok(\lambda)\to 0.\]

From the $\sigma(\alpha_g, \beta_g)$-stability of $F$ with $\mu^0_{\alpha_g, \beta_g}(F)=+\infty$, we know that either $\cok(\lambda)\cong F$ or $\cok(\lambda)=0$.

\medskip

\textbf{Case 1.1}: $\cok(\lambda)\cong F$.

In this case, we have $i^*(A)\cong \ker(\lambda)$ and hence we obtain a triangle
\[i^*(A)[1]\to i^*(B)\to F.\]
Using $A,B\in \Coh^0_{\alpha_g, \beta_g}(X)$ and \eqref{eq:abc-slope}, we have $$\Hom(A, \cE)=\Hom(A[1], \cE)=\Hom(A[1], B)=0.$$ Therefore, applying $\Hom(A[1],-)$ to the triangle \eqref{prB}, we obtain $$\Hom(A[1], i^*(B))=\Hom(i^*(A)[1], i^*(B))=0.$$ This implies $i^*(B)\cong \cok(\lambda) \cong  F$ and $\ker(\lambda)\cong i^*(A) \cong 0$. However, by \eqref{prA}, we know that $t=\Ext^1(\cE, A)=0$ and $A\cong \cE^{\oplus n_g-c}$, contradicting \eqref{eq:abc-slope}.

\medskip

\textbf{Case 1.2}: $\cok(\lambda)\cong 0$.

In this case, we have $i^*(B)\cong \ker(\lambda)[1]$ and two exact sequences in $\Coh^0_{\alpha_g, \beta_g}(X)$
\[0\to \ker(\lambda)\to i^*(A)\to F\to 0\]
and
\begin{equation} \label{EB}
    0\to \ker(\lambda)\xra{\theta} \cE^{\oplus c}\xra{\lambda} B\to 0.
\end{equation}

Thus, if we apply $\Hom(-, \cE)$ to the sequence \eqref{EB} and use \eqref{CE=0}, we obtain that $$\Hom(\ker(\lambda), \cE)\cong \Hom(\cE^{\oplus c}, \cE) =k^c.$$ Then $\theta$ coincides with the natural map $\ker(\lambda)\to \cE\otimes \Hom(\ker(\lambda), \cE) \cong\cE^{\oplus c}$. Hence, we know that  $\pr_j\circ \theta\neq 0$ for any $1\leq j\leq c$, where $\pr_j\colon \cE^{\oplus c}\to \cE$ is the projection map of the $j$-th component.

By Claim 3, we have $\mu^0_{\alpha_g, \beta_g}(B)\geq \mu^0_{\alpha_g, \beta_g}(\cE)$. If $\mu^0_{\alpha_g, \beta_g}(B)= \mu^0_{\alpha_g, \beta_g}(\cE)$, since $\cE$ is $\sigma^0_{\alpha_g, \beta_g}$-stable and $B$ is the quotient of $\cE^{\oplus c}$, we know that $B\cong \cE^{\oplus \frac{\rk B}{2}}$, and this contradicts \eqref{CE=0}. Hence we have $\mu^0_{\alpha_g, \beta_g}(B)> \mu^0_{\alpha_g, \beta_g}(\cE)$. From this and the $\sigma^0_{\alpha_g, \beta_g}$-stability of $\cE$, we also know that $\mu^0_{\alpha_g, \beta_g}(\ker(\lambda))< \mu^0_{\alpha_g, \beta_g}(\cE)$ and hence $\Im(Z^0_{\alpha_g, \beta_g}(\ker(\lambda)))>0$.

To summarize, in this case, we have:

\begin{itemize}
    \item $\Im(Z^0_{\alpha_g, \beta_g}(A))\geq 0, \Im(Z^0_{\alpha_g, \beta_g}(B))> 0$,
    
    \item $\Im(Z^0_{\alpha_g, \beta_g}(\ker(\lambda)))>0$, 
    
    \item $\mu^0_{\alpha_g, \beta_g}(\ker(\lambda))<\mu^0_{\alpha_g, \beta_g}(\cE)$,

    \item $\mu^0_{\alpha_g, \beta_g}(C)>\mu^0_{\alpha_g, \beta_g}(B)> \mu^0_{\alpha_g, \beta_g}(\cE)$,
    
    \item $c>0, a<0$.

\end{itemize}

From Claim 2, we know that $\ker(\lambda)=i^*(B)[-1]$ is a torsion-free sheaf. By Lemma \ref{inequality}, these inequalities imply $\mu_H(\ker(\lambda))\geq \mu_H(\cE)$. By Claim 2, the only possible situation is when $\mu_H(\ker(\lambda))=\mu_H(\cE)$ and $\ker(\lambda)$ is a $\mu_H$-semistable sheaf. In this case, by the $\mu_H$-stability of $\cE$, we know that $\im(\theta)$ is contained in $\cE^{\oplus c'}$ for $c'=\frac{\rk(\im(\theta))}{2}$. But from $\pr_j\circ \theta\neq 0$ for any $1\leq j\leq c$, the only possible case is $c=c'$, and either $\cH^{-1}(B)\cong 0$ and $B\cong \cH^{0}(B)$ is a torsion sheaf supported in codimension $\geq 2$, or $\cH^{-1}(B)\neq 0$ is a $\mu_H$-semistable sheaf with $\mu_H(\cH^{-1}(B))=\mu_H(\cE)$. But the first case contradicts $\mu^0_{\alpha_g, \beta_g}(B)<0$ and the second case contradicts $\mu_H^+(\cH^{-1}(B))\leq \beta_g$. Therefore, such a minimal destabilizing quotient object $B$ cannot exist, and we can conclude that $C$ is $\sigma^0_{\alpha_g, \beta_g}$-semistable when $g\neq 6$.

\bigskip

\noindent\textbf{Case 2}: $g=6$.

In this case, by the claims above, we also have a system of inequalities:

\begin{itemize}
    \item $\Im(Z^0_{\alpha_g, \beta_g}(A))\geq 0, \Im(Z^0_{\alpha_g, \beta_g}(B))> 0$,
    
    \item $\Im(Z^0_{\alpha_g, \beta_g}(i^*(A)))\geq 0$,

    \item $\mu^0_{\alpha_g, \beta_g}(C)>\mu^0_{\alpha_g, \beta_g}(B)\geq  \mu^0_{\alpha_g, \beta_g}(\cE)$,
    
    \item $c>0,a<0$.

\end{itemize}

By Claim 2, we have $\mu_H(i^*(B)[-1])\leq \mu_H^+(i^*(B)[-1])< \mu_H(\cE).$ Therefore, Lemma \ref{g=6inequality} implies that the only possible classes of $i^*(B)[-1]$ are $[i^*(B)[-1]]=v-w$ or $3v-2w$.

\medskip

\textbf{Case 2.1}: $[i^*(B)[-1]]=v-w$.

By Claim 2, we know that $i^*(B)[-1]$ is a torsion-free sheaf with $\ch(i^*(B)[-1])=\ch(I_C(-H))$, where $C$ is a conic on $X$. Thus $i^*(B)[-1]\cong I_C(-H)$ for some conic $C$ on $X$. But then $\Ext^3(\oh_X, I_C(-H))=\Hom(I_C, \oh_X)\neq 0$ by Serre duality, contradicts $i^*(B)\in \Ku(X)$.

\medskip

\textbf{Case 2.2}: $[i^*(B)[-1]]=3v-2w$.

In this case, $i^*(B)[-1]$ is a torsion-free sheaf with rank $3$ and $\mu_H(i^*(B)[-1])=-\frac{2}{3}$. We first show that $i^*(B)[-1]$ is $\mu_H$-stable. Indeed, if the minimal destabilizing quotient sheaf of $i^*(B)[-1]$ has rank one, then we have $$\Hom(i^*(B)[-1], \oh_X(-nH))\neq 0$$ for some $n\geq 1$. As $\oh_X(-nH)[1]\in \Coh^{\beta_g}(X)$ is $\sigma_{\alpha_g, \beta_g}$-semistable, from Claim 3, we know that $\Hom(B, \oh_X(-nH)[1])=0$ for every $n\geq 1$. Now, applying $\Hom(-, \oh_X(-nH))$ to the triangle \eqref{prB} and using $\cE, B\in \Coh^{\beta_g}(X)$, we get $\Hom(i^*(B)[-1], \oh_X(-nH))=0$, a contradiction. If the maximal destabilizing subsheaf of $i^*(B)[-1]$ has rank one, then we know that 
$$\Hom(\oh_X(nH),i^*(B)[-1])\neq 0$$ for some $n\geq 0$. Since $\Hom(\oh_X(nH), \cE)=0$ and $\oh_X(nH)\in \Coh^{\beta_g}(X)$, if we apply the functor $\Hom(\oh_X(nH),-)$ to the triangle \eqref{prB}, we obtain $\Hom(\oh_X(nH), i^*(B)[-1])=0$, again a contradiction. Therefore, $i^*(B)[-1]$ is $\mu_H$-stable.

Now, by Lemma \ref{bms lemma 2.7} and Lemma \ref{wall_3v-2w}, we know that $i^*(B)[-1]\in \Coh^{\beta_g}(X)$ is $\sigma_{\alpha, \beta_g}$-semistable for all $\alpha>0$. In particular, $i^*(B)[-1]$ is $\sigma_{\alpha_g, \beta_g}$-semistable, hence $i^*(B)\in \Coh^0_{\alpha_g, \beta_g}(X)$ is $\sigma^0_{\alpha_g, \beta_g}$-semistable and is also $\sigma(\alpha_g, \beta_g)$-semistable. But by $\sigma(\alpha_g, \beta_g)$-stability of $F$ and $$\mu^0_{\alpha_g, \beta_g}(F)>\mu^0_{\alpha_g, \beta_g}(i^*(B)),$$ we get $\Hom(F, i^*(B))=\Hom(F, B)=0$, a contradiction.

Therefore, such a minimal destabilizing quotient object $B$ cannot exist, and we conclude that $C$ is $\sigma^0_{\alpha_g, \beta_g}$-semistable when $g=6$.
\end{proof}

\begin{proof}[Proof of Theorem \ref{thm: BN = X}]
Take $F\in \mathsf{BN}_g \setminus p(X)$. By Proposition \ref{gneq6cone}, the cone $C$ of the natural map $F[-1] \ra \cE^{\oplus n_g}$ is in $\Coh^0_{\alpha_g, \beta_g}(X)$ and is $\sigma^0_{\alpha_g, \beta_g}$-semistable. From $\mu^0_{\alpha_g, \beta_g}(C)<0$, we know that $C\in \Coh_{\alpha_g, \beta_g}(X)$ is $\sigma_{\alpha_g, \beta_g}$-semistable by \cite[Proposition 4.1]{FeyzbakhshPertusi2021stab}. Moreover, $C$ is $\sigma_{\alpha, \beta_g}$-semistable for every $\alpha>0$ by $\Delta_H(C)=0$ and \cite[Corollary 3.10]{bayer2016space}. Thus, $C$ is a torsion-free $\mu_H$-semistable sheaf with $\ch(C)=\ch(I_x)$. From $\Pic(X)\cong \mathbb{Z}H$, we conclude that $C\cong I_x$ for a point $x\in X$.

Combining with Lemma \ref{prop:F is shift of vb}, we have a short exact sequence $$0 \ra F[-1] \ra \cE^{\oplus n_g} \ra I_x \ra 0$$ of sheaves. Applying $\Hom(\oh_x, -)$ to this gives a long exact sequence containing
\[ \cdots \ra \Ext^1(\oh_x, \cE^{\oplus n_g}) \ra \Ext^1(\oh_x, I_x) \ra \Ext^2(\oh_x, F[-1]) \ra \cdots . \]
Since $\cE$ and $F[-1]$ are vector bundles, the first and last terms of the above vanish. However, $\Ext^1(\oh_x, I_x) \neq 0$, and we have a contradiction. This gives $\mathsf{BN}_g= p(X)$ as required.
\end{proof}

\subsection{Refined categorical Torelli}

The above theorem can be used to prove refined categorical Torelli theorems for all index one Fano threefolds of genus $g \geq 6$. 

We first introduce a functor on $\Ku(X)$, defined by $$T(-)\coloneqq i^*\circ \mathcal{H}om(-, \oh_X(-H)[1]).$$ By \cite[Proposition 3.8]{zhang2020bridgeland}, this functor $T\colon \Ku(X)\to \Ku(X)$ is an antiequivalence with the property $T\circ T\cong \mathrm{id}_{\Ku(X)}$. Moreover, it induces a linear isometry on $\cN(\Ku(X))$, which we also denote by $T$. When $g=6$, we have $T(v)=-3v+2w$ and $T(w)=-4v+3w$. Note that when $g=6$, there also exists an involution $\Ku(X)\rightarrow\Ku(X)$, which is $\tau=S_{\Ku(X)}[-2]$. The action of $\tau$ on $\cN(\Ku(X))$ is trivial.

\begin{Lem} \label{isometry}
Let $g\geq 6$ and $$\Phi\colon \cN(\Ku(X))=\langle v, w\rangle \to \cN(\Ku(X'))=\langle v', w'\rangle$$ be a linear isometry with respect to Euler forms, such that $\Phi([i^!\cE_X])=[i'^!\cE_{X'}]$.

\begin{itemize}
    \item When $g\geq 7$, we have $\Phi(v)=v'$ and $\Phi(w)=w'$.

    \item When $g=6$, we have $\Phi(v)=v'$ and $\Phi(w)=w'$ or $T\circ \Phi(v)=v'$ and $T\circ \Phi(w)=w'$.
\end{itemize}
\end{Lem}

\begin{proof}
Since $\Phi$ preserves the Euler form and $\Phi([i^!\cE_X])=[i'^!\cE_{X'}]$, an elementary computation shows that $\Phi(v)=v'$ and $\Phi(w)=w'$ if $g\geq 7$. If $g=6$, we have either $\Phi(v)=v'$ and $\Phi(w)=w'$ or $\Phi(v)=-3v'+2w'$ and $\Phi(w)=-4v'+3w'$. Then the result follows from $T(v')=-3v'+2w'$, $T(w')=-4v'+3w'$ and $T\circ T=\mathrm{id}$.
\end{proof}

\begin{Lem} \label{inv_on_data}
When $g=6$, we have $T(i^{!}\cE)\cong \tau(i^{!}\cE)$.
\end{Lem}

\begin{proof}
By Lemma \ref{cohomology objects of pi lemma}, we have a triangle 
\[\cQ(-H)[1]\to i^!\cE\to \cE.\]
Note that $T(\cQ(-H)[1])\cong i^*(\cQ^{\vee})$ and $T(\cE)\cong i^*(\cE[1])=0$, so we have $T(i^!\cE)\cong i^*(\cQ^{\vee})$. Then the result follows from $\tau\cong \tau^{-1}$ and $i^*(\cQ^{\vee})=\bL_{\cE}\cQ^{\vee}\cong \tau^{-1}(i^!\cE)$.
\end{proof}

\begin{proof}[{Proof of Theorem \ref{Refined_Torelli}}]
First, we assume that $g\neq 6$. By Lemma \ref{isometry} and Theorem \ref{all_in_one_orbit}, the functor $\Phi$ induces an isomorphism $M_{\sigma}(-[i^*\oh_x])\cong M_{\sigma'}(-[i'^*\oh_{x'}])$. Now, the result follows from $\Phi(i^!\cE)\cong i'^!\cE_{X'}$ and Theorem~\ref{thm: BN = X}. When $g=6$, by Lemma \ref{isometry}, after replacing $\Phi$ with $\tau\circ T\circ \Phi$ if necessary,  we can assume that $\Phi([i^*\oh_x])=[i'^*\oh_{x'}]$ and $\Phi(i^!\cE_X)\cong i'^!\cE_{X'}$. Therefore, the result in this case also follows from Theorem \ref{all_in_one_orbit} and Theorem~\ref{thm: BN = X}.
\end{proof}

\section{Autoequivalences of Kuznetsov components of index one prime Fano threefolds}\label{Section_auto-equi}
Our goal in this section is to describe exact autoequivalences of Kuznetsov components of index one prime Fano threefolds of genus $6$ and $8$ as in Theorem~\ref{theorem_group_auto_equi_index_one}.

Let $X$ be an index one prime Fano threefold of even genus $g\geq 6$. Recall that by Remark \ref{rem:equiv kuz components}, there is an equivalence $\Xi\colon \Ku(X)\xra{\simeq} \cA_X$. We begin with a lemma, which is a generalization of the argument in \cite[Section 4]{liu2023autoeq}. 

\begin{Lem} \label{lem_kernel_trivial}
    Let $X$ be a smooth projective variety and $i\colon \cA\hookrightarrow \Db(X)$ be an admissible subcategory with a left adjoint $i^*$. Assume furthermore that there exists $n\in \ZZ$ such that $i^*\oh_x[-n]=K^x$ is a Gieseker-stable sheaf for any $x\in X$ and $i^*$ induces a closed immersion
    \[\pi\colon X\hookrightarrow M_X([K^x]),\quad x\mapsto i^*\oh_x[-n]=K^x\]
    where $M_X([K^x])$ is a fine moduli space of Gieseker-semistable sheaves of class $[K^x]$. Assume furthermore that for any non-trivial line bundle $L$ on $X$, there is a point $x\in X$ with
    \begin{equation}\label{assumption}
        K^x\neq i^*(K^x\otimes L).
    \end{equation}
    Then a Fourier--Mukai type functor $\Phi\colon \cA\to \cA$ is isomorphic to $\identity_{\cA}$ if and only if $\Phi(K^x)\cong K^x$ for any point $x\in X$.
\end{Lem}

\begin{proof}
The 'only if' part is obvious. We now prove the 'if' part.

By \cite[Lemma 3.31]{huyb-book-FM}, the kernel of $\Phi$ is isomorphic to a coherent sheaf $\mathcal{K}[n]$ on $X\times X$ flat over the first factor and $i_x^*\cK\cong K^x$ for any $x\in X$ and the closed immersion $i_x\colon \{x\}\times X \hookrightarrow X\times X$. Since $M_X([K^x])$ is fine, we have an induced morphism $\pi'\colon X\to M_X([K^x])$. By \cite[Theorem 7.1]{kuznetsov:base-change}, the functor $i\circ i^*\colon \Db(X)\to \Db(X)$ is of Fourier--Mukai type, and let $\cP[n]$ be the kernel. Note that $\identity_{\cA}=i\circ i^*|_{\cA}$. Then $\cP$ also induces a morphism $X\to M_X([K^x])$, which is exactly $\pi$.

It is clear that $\pi$ and $\pi'$ have the same image as continuous maps. Then by taking the scheme-theoretic image of $\pi'$ and the uniqueness of the reduced scheme structure of a closed subset, we see that $\pi'$ can be factored as $\pi\circ f$, where $f\colon X\to X$ is a bijective morphism. Since $X$ is smooth and connected, we deduce that $f=\identity_X$, hence $\pi'$ is isomorphic to $\pi$. By the definition of the moduli functor, this implies that $\cK$ and $\cP$ differ only by a line bundle $L$ on $X$ (cf.~\cite[Section 4.1]{huybrechts:geometry-of-moduli-space-of-sheaves}).

We claim that $L\cong \oh_X$ and hence $\Phi\cong \identity_{\cA}$. Indeed, if we denote $q_i$ the projection of $X\times X$ to the $i$-th factor, then we have
\[K^x\cong \Phi(K^x)\cong q_{2*}(\cP\otimes q_1^*(K^x\otimes L))\cong i^*(K^x\otimes L)\in \cA\]
for any point $x\in X$. Then by the assumption \eqref{assumption}, we see $L\cong \oh_X$ and the result follows.
\end{proof}

\begin{Ex}\label{example_index_one}
    Let $X$ be a prime Fano threefold of even genus $g\geq 6$. When $g=6$, we furthermore assume that $X$ is ordinary. If we take $\cA=\Ku(X)$, Lemma~\ref{Kstable} and the same argument in Proposition~\ref{closed_embedding_into_modulispace} show that $i^*$ induces the closed immersion $\pi$ in Lemma \ref{lem_kernel_trivial}. Moreover, the moduli space is fine since the corresponding Chern character is primitive (cf.~Remark~\ref{rem:K shorthand}). Finally, it is straightforward to check that
    \[\ch(K_x)\neq \ch(i^*(K_x\otimes L))\]
    for any $L\neq \oh_X\in \Pic(X)=\ZZ H$. Hence the assumption \eqref{assumption} is also satisfied in this case.
\end{Ex}

We also need the following two lemmas. Recall that for a prime Fano threefold $Y$ of index $2$, we define $\Ku(Y)\coloneqq \langle \oh_Y, \oh_Y(1) \rangle^{\perp}$. The inclusion $\Ku(Y)\hookrightarrow \Db(Y)$ is denoted by $l$.

\begin{Lem}\label{lem_bundle}
Let $Y$ and $Y'$ be two prime Fano threefolds of index $2$ and degree $2\leq d\leq 3$. Let $\Psi\colon \Ku(Y)\xra{\simeq} \Ku(Y')$ be an exact equivalence that maps classes $\bv$ and $\bw$ to $\bv'$ and $\bw'$, respectively. If $E\in \Ku(Y)$ is a vector bundle, $\Psi(E)$ is also a bundle up to shift.
\end{Lem}

\begin{proof}
    By \cite[Remark 7.3]{feyzbakhsh2023new} and \cite[Theorem 6.2]{feyzbakhsh2023new}, we know that for any $p'\in Y$, there is a point  $p\in Y$ such that $\Psi(l^*\oh_p)\cong l'^*\oh_{p'}$. Then we have
    \[\RHom^{\bullet}(\oh_{p'}, \Psi(E))=\RHom^{\bullet}(l'^*\oh_{p'}, \Psi(E))=\RHom^{\bullet}(l^*\oh_{p}, E)=\RHom^{\bullet}(\oh_p, E),\]
where the first and the last equalities follow from the fact that $l^*$ and $l'^*$ are left adjoint to $l$ and $l'$, respectively. Then the local freeness of $\Psi(E)$ follows from the local freeness of $E$ and \cite[Proposition 5.4]{bm99}.
\end{proof}

\begin{Lem}\label{unique_extend_index_one}
    Let $f,g\colon X\to X'$ be two isomorphisms between index one prime Fano threefolds of genus $g\geq 6$. If $f_{*}|_{\Ku(X)}=g_*|_{\Ku(X)}\colon \Ku(X)\to \Ku(X')$, then $f=g$. Thus the homomorphism 
    \[\Aut(X)\to \Aut(\Ku(X)), \quad f\mapsto f_*|_{\Ku(X)}\]
    is injective.
\end{Lem}
\begin{proof}
It is clear that
\begin{equation*}
f_*(i^*\oh_p)=i'^*\oh_{f(p)} \qquad \text{and } \qquad g_*(i^*\oh_p)=i'^*\oh_{g(p)}. 
\end{equation*}
Since $f_{*}|_{\Ku(X)}=g_*|_{\Ku(X)}$, we get $i'^*\oh_{f(p)} = i'^*\oh_{g(p)}$. Thus, the injectivity part in Proposition~\ref{closed_embedding_into_modulispace} implies that $f(p) = g(p)$ for any point $p \in X$. Since both $X$ and $X'$ are smooth, we get $f=g$.  
\end{proof}

\subsection{Genus 8 case}

We start with genus $8$ case.

\begin{Thm}\label{thm_genus_8}
    Let $X$ be an index one prime Fano threefold of genus $8$. Then for any exact autoequivalence $\Phi\colon \Ku(X)\xra{\simeq} \Ku(X)$, after composing with the Serre functor $S_{\Ku(X)}$ and shift functor, we have
    \[\Phi(i^!\cE)\cong i^!\cE.\]
\end{Thm}

\begin{proof}
By \cite[Theorem 1.3]{li2022derived}, $\Phi$ is of Fourier--Mukai type. We define an equivalence $$\Phi'\colon \cA_X\xra{\simeq} \cA_X$$ by $\Phi'\coloneqq\Xi\circ \Phi\circ \Xi^{-1}$. Let $Y$ be the Pfaffian cubic threefold associated with $X$. Then by \cite[Theorem 4.7]{kuznetsov:fano-threefolds} and Proposition~\ref{Prop_instanton_bundle}, there is an equivalence of Fourier--Mukai type $\Theta\colon \cA_X\xra{\simeq} \Ku(Y)$ which maps $\Xi(i^!\cE)$ to an instanton bundle of charge $2$ on $Y$ up to shift in the sense of \cite[Definition 2.4]{kuznetsov2003derived}; see also Definition \ref{def:instanton}. We set $\Psi\coloneqq\Theta\circ \Phi'\circ \Theta^{-1}$. 

By \cite[Lemma 7.4]{feyzbakhsh2023new}, after composing $\Phi$ with the Serre functor $S_{\Ku(X)}$ and shift functor, we can assume that $\Psi$ acts trivially on $\cN(\Ku(Y))$. Hence, using Lemma~\ref{lem_bundle}, $\Psi(\Theta(\Xi(i^!\cE)))$ is also a bundle up to shift. Moreover, $\Psi(\Theta(\Xi(i^!\cE)))$ is stable with respect to every Serre-invariant stability condition on $\Ku(Y)$ as $\Theta(\Xi(i^!\cE))$ is (cf.~\cite[Theorem 7.6]{liu2021note}). By \cite[Theorem 7.6]{liu2021note} and the local freeness of $\Psi(\Theta(\Xi(i^!\cE)))$, we can assume that $\Psi(\Theta(\Xi(i^!\cE)))$ is also an instanton bundle of charge $2$ on $Y$ after composing $\Phi$ with shift functor. Then from \cite[Theorem 4.7]{kuznetsov:fano-threefolds} and Proposition~\ref{Prop_instanton_bundle}, there is another index one genus $8$ prime Fano threefold $X'$ with an equivalence $\Theta'\colon \cA_{X'}\xra{\simeq} \Ku(Y)$ such that $\Theta'(\Xi'(i'^!\cE_{X'}))\cong \Psi(\Theta(\Xi(i^!\cE_X)))$, where $\Xi'\colon \Ku(X')\to \cA_{X'}$ is the equivalence in Remark \ref{rem:equiv kuz components}. 

By Theorem \ref{Refined_Torelli}, there is an isomorphism $g\colon X\to X'$ such that $$\Theta'^{-1}\circ \Theta(i^*\oh_x)=i'^*\oh_{g(x)}$$ for any $x\in X$. Using Lemma~\ref{lem_kernel_trivial} and Example \ref{example_index_one}, we obtain $\Theta'^{-1}\circ \Theta=g_*$. Thus, the isomorphism $\Theta'(\Xi'(i'^!\cE_{X'}))\cong \Psi(\Theta(\Xi(i^!\cE_X)))$ gives $$\Xi'(i'^!\cE_{X'})\cong g_*(\Phi'(\Xi(i^!\cE_X))).$$ Therefore, we get $g^*(\Xi'(i'^!\cE_{X'}))\cong \Phi'(\Xi(i^!\cE_X))$. By definition and the uniqueness of $\cE$ (cf.~\cite{bayer2024mukai}), we have $g^*(\Xi'(i'^!\cE_{X'}))\cong \Xi(i^!\cE_{X})$ and the result follows.
\end{proof}

\begin{Cor}\label{cor_aut_genus_8}
     Let $X$ be an index one prime Fano threefold of genus $8$. Then we have
     \[\Aut(\Ku(X))= \langle \Aut(X), S_{\Ku(X)}, [1] \rangle.\]
\end{Cor}

\begin{proof}
By Lemma~\ref{unique_extend_index_one}, we have an injection 
\[\Aut(X)\to \Aut(\Ku(X)), \quad f\mapsto f_*|_{\Ku(X)}.\]
Note that $\Aut(X)$ acts trivially on $\cN(\Ku(X))$, while the only elements in $\langle S_{\Ku(X)},[1]\rangle$ that act trivially on $\cN(\Ku(X))$ are of the form $[2m]$. Hence $\Aut(X)\cap \langle S_{\Ku(X)},[1]\rangle=\identity_{\Ku(X)}$ and we see that the induced homomorphism
\[\eta\colon \Aut(X)\to \frac{\Aut(\Ku(X))}{\langle S_{\Ku(X)},[1]\rangle}\]
is also injective. On the other hand, we have a homomorphism
\[\eta'\colon \frac{\Aut(\Ku(X))}{\langle S_{\Ku(X)},[1]\rangle}\to \Aut(X)\]
given by Theorem~\ref{thm_genus_8} and Theorem~\ref{main_theorem_Brill_Nother_reconstruction}. Then any element in the kernel of $\eta'$ can be represented by an autoequivalence $\Phi$ such that $\Phi(i^*\oh_p)\cong i^*\oh_p$ for any point $p\in X$. Then by Lemma~\ref{lem_kernel_trivial} and Example \ref{example_index_one}, we have $\Phi\cong\identity_{\Ku(X)}$ and $\eta'$ is injective as well. It is straightforward to check that $\eta'\circ \eta=\identity_{\Aut(X)}$, hence $\eta$ and $\eta'$ are inverse to each other and the result follows.
\end{proof}

Corollary~\ref{cor_aut_genus_8} has an immediate application to automorphism groups of Fano threefolds. 

\begin{Cor}\label{cor_aut_variety}
Let $X$ be an index one prime Fano threefold of genus $8$ and $Y$ be the Pfaffian cubic threefold associated with $X$. Then we have
    \[\Aut(X)\cong\Aut(Y).\]
\end{Cor}

\begin{proof}
    By \cite[Theorem 4.7]{kuznetsov:fano-threefolds}, we have an equivalence $\Ku(X)\simeq \Ku(Y)$. Hence we get an isomorphism 
    \[s\colon \Aut(\Ku(X))\cong \Aut(\Ku(Y)).\]
    Since an exact equivalence commutes with the Serre functor and shift functor, we see that $s(\langle S_{\Ku(X)}, [1] \rangle)=\langle S_{\Ku(Y)}, [1] \rangle$. Taking quotient on both sides and using Corollary \ref{cor_aut_genus_8} and \cite[Corollary 8.4]{feyzbakhsh2023new}, we get an induced isomorphism
    \[\Aut(X)\cong\Aut(Y).\]
\end{proof}

\subsection{Genus 6 case}

Using results of \cite{jacovskis2024categorical}, the case of genus $6$ can be treated analogously.

\begin{Thm} \label{thm_genus_6}
    Let $X$ be a general ordinary Gushel--Mukai threefold. Then for any exact autoequivalence $\Phi\colon \Ku(X)\xra{\simeq} \Ku(X)$, after composing with the Serre functor $S_{\Ku(X)}$ and shift functor, we have
    \[\Phi(i^!\cE_X)\cong i^!\cE_X.\]
\end{Thm}

\begin{proof}
By \cite[Theorem 1.3]{li2022derived}, $\Phi$ is of Fourier--Mukai type. Let $\Psi\coloneqq\Xi\circ \Phi\circ \Xi^{-1}\colon \cA_X\simeq \cA_X$ be the induced equivalence. Since $X$ is general, $X$ is not the period dual of itself since the involution on the period domain defined in  \cite[(1.0.14)]{ogrady:double-EPW-period} is non-trivial. Hence, by \cite[Theorem 10.3]{jacovskis2024categorical} and its proof, up to shift, $\Psi$ and $\Phi$ act trivially on $\cN(\cA_X)$ and $\cN(\Ku(X))$, respectively. Moreover, by \cite[Theorem 7.12]{jacovskis2024categorical},  $\Psi$ induces an automorphism $\cC_m(X)\cong \cC_m(X)$ of the minimal model of the Fano surface of conics on $X$, which is either the identity map or an involution (cf.~\cite[Corollary 9.2]{debarre2012period}). Then by \cite[Proposition 7.13]{jacovskis2024categorical}, up to composing with $S_{\cA_X}$ and shift functor, $\Psi$ induces an automorphism $\cC_m(X)\cong \cC_m(X)$ that maps the point $\Xi(i^!\cE_X)$ to $\Xi(i^!\cE_X)$, and the result follows.
\end{proof}

\begin{Cor} \label{cor_aut_genus_6}
     Let $X$ be a general ordinary Gushel--Mukai threefold. Then we have
     \[\Aut(\Ku(X))=   \langle \Aut(X), S_{\Ku(X)},  [1]\rangle.\]
\end{Cor}

\begin{proof}
By Lemma \ref{unique_extend_index_one}, we have an injection 
\[\Aut(X)\to \Aut(\Ku(X)), \quad f\mapsto f_*|_{\Ku(X)}.\]
Since $\Aut(X)\cap \langle S_{\Ku(X)},[1]\rangle=\identity_{\Ku(X)}$, we see that the induced homomorphism
\[\eta\colon \Aut(X)\to \frac{\Aut(\Ku(X))}{\langle S_{\Ku(X)},[1]\rangle}\]
is injective as well. On the other hand, we have a homomorphism
\[\eta'\colon \frac{\Aut(\Ku(X))}{\langle S_{\Ku(X)},[1]\rangle}\to \Aut(X)\]
given by Theorem \ref{thm_genus_6} and Theorem~\ref{main_theorem_Brill_Nother_reconstruction}. Then any element in the kernel of $\eta'$ can be represented by an autoequivalence $\Phi$ such that $\Phi(i^*\oh_p)\cong i^*\oh_p$ for any point $p\in X$. Then by Lemma \ref{lem_kernel_trivial} and Example \ref{example_index_one}, we have $\Phi\cong\identity_{\Ku(X)}$ and $\eta'$ is also injective. It is straightforward to check that $\eta'\circ \eta=\identity_{\Aut(X)}$, hence $\eta$ and $\eta'$ are inverse to each other and the result follows.
\end{proof}

Corollary~\ref{cor_aut_genus_6} has a nice application to Kuznetsov's Fano threefold conjecture \cite[Conjecture 3.7]{kuznetsov:fano-threefolds}. It was disproved in \cite{bayer2022kuznetsov} and \cite{zhang2020bridgeland} independently in its most general form. We present a simple disproof by assuming the Gushel--Mukai threefold is general.

\begin{Cor} \label{cor_ku_conj}
Let $X$ be a general Gushel--Mukai threefold and $Y$ be a prime Fano threefold of index $2$ and degree $2$. Then $\Ku(X)$ is not equivalent to $\Ku(Y)$. 
\end{Cor}

\begin{proof}
Assume that there is an exact equivalence $\Ku(X)\simeq \Ku(Y)$. Then it induces an isomorphism of the numerical Grothendieck groups and 
    \[\Aut(\Ku(X))\cong \Aut(\Ku(Y)).\]
However, Corollary \ref{cor_aut_genus_6} shows that any element in $\Aut(\Ku(X))$ acts on $\cN(\Ku(X))$ by identity up to sign, while the action of the rotation functor on $\cN(\Ku(Y))$ is non-trivial (cf.~\cite[Section 2.2]{feyzbakhsh2023new}). Thus, we get a contradiction.
\end{proof}

\begin{appendix} 

\section{A calculation in the genus 8 case}

Let $X$ be a prime Fano threefold of index $1$ and genus $8$. Following \cite{kuznetsov2003derived}, we write
\[
    X=\mathbb P(f(A)^{\perp})\cap \Gr(2,V)\subset \mathbb P(\wedge^2V),
\]
where $V$ is six-dimensional, $A\subset \wedge^2V^{\vee}$ is a
five-dimensional vector subspace, $f\colon A\to \wedge^2 V^{\vee}$ is a linear map, and $f(A)^{\perp}$ denotes the annihilator of $f(A)\subset \wedge^2 V^{\vee}$.

Let $Y=\mathsf{Pf}(f)\subset \mathbb P(A)$ be the corresponding Pfaffian cubic threefold. We denote by $H_Y$ the hyperplane class of $Y$.

\begin{Def}\label{def:instanton}
An \emph{instanton sheaf of charge} $2$ on $Y$ is a $\mu_{H_Y}$-stable sheaf $E$ such that $H^1(E(-H_Y))=0$ and $\ch(E)=2-2l$. Here, $l$ is the class of lines on $Y$.
\end{Def}

As in \cite[Section 2]{kuznetsov2003derived}, a natural instanton bundle $E_f(-H_Y)$ of charge $2$ on $Y$ can be constructed as follows. The map $f$ naturally corresponds to a morphism $$V\otimes \oh_{\PP(A)}(-1)\to V^{\vee}\otimes \oh_{\PP(A)},$$ which is also denoted by $f$. Such a morphism is injective, and its cokernel is schematically supported on $Y$. In particular, there is a unique sheaf $E_f\in \Coh(Y)$ so that $\mathrm{cok}(f)\cong \alpha_*E_f$, where $\alpha\colon Y\hookrightarrow \PP(A)$ is the inclusion. By \cite[Theorem 2.2]{kuznetsov2003derived}, $E_f(-H_Y)$ is an instanton bundle of charge $2$ on $Y$.

Recall that for $X$, we use the alternative Kuznetsov component $\cA_X$
and the equivalence
\[
    \Xi\colon \Ku(X)\xrightarrow{\sim}\cA_X .
\]
Set
\[
    \cG\coloneqq \Xi(i^!\cE_X)\in\cA_X.
\]
As recalled in Section~\ref{sec:derivedfano}, the object $\cG$ is given by the exact triangle
\begin{equation}\label{eq:G-triangle-appendix}
    \cE[2]\longrightarrow \cG\longrightarrow \cQ^{\vee}[1].
\end{equation}
Here, $\cE$ is the restriction of the tautological rank-two subbundle on
$\Gr(2,V)$, and $\cQ$ is the restriction of the tautological quotient bundle.

Let
\[
    \Theta\coloneqq \bL_{\oh_Y}\circ (-\otimes \oh_Y(H_Y))\circ \Phi_{W}\colon \Db(X)\longrightarrow \Db(Y),
\]
where $\Phi_{W}$ is the Fourier--Mukai functor with kernel $I_W\otimes p_Y^*\oh_Y(H_Y)$, and $W\subset X\times Y$ is the closed subscheme parameterizing all points $$(U,a)\in \Gr(2,V)\times \PP(A)$$ such that the kernel of the skew form $a\in A$ intersects the two-dimensional subspace $U\subset V$. Then by \cite{kuznetsov2003derived} and \cite[Remark B.6.5]{KPS}, $\Theta$ induces an equivalence between $\cA_X$ and $\Ku(Y)$.  We denote by $p_X,p_Y$ the two projections from $X\times Y$.

This appendix aims to prove the following.

\begin{Prop}\label{Prop_instanton_bundle}
We have $\Theta(\cG)\cong E_f(-H_Y)$.
\end{Prop}

\begin{proof}
Let $j\colon W\hookrightarrow Y\times X$ be the inclusion. We first compute $\Phi_W(\cE)$. From the exact sequence $0\to I_W\to \oh_{X\times Y}\to j_*\oh_W\to 0$, we get a triangle $$p_X^*\cE\otimes I_W(H_Y)\to p_X^*\cE \otimes p_Y^*\oh_Y(H_Y)\to p_X^*\cE \otimes j_*((p_Y^*\oh_Y(H_Y))|_W).$$ As $\RHom^{\bullet}(\oh_X,\cE)=0$, applying $p_{Y*}$ to this triangle, we get $$\Phi_W(\cE)\cong p_{Y*}(j_*\cE|_W)\otimes \oh_Y(H_Y)[-1].$$
Note that $W_y\coloneqq W\cap p_Y^{-1}(y)$ is a curve of degree $4$ and arithmetic genus $0$, so it is straightforward to check $\Phi_W(\cQ^{\vee})\cong 0$. Therefore, we obtain 
\[p_{Y*}(j_*\cE|_W)\otimes \oh_Y(H_Y)[1]\cong \Phi_W(\cE)[2]\cong \Phi_W(\cG).\]
Since $\bL_{\oh_Y}(E_f(-H_Y))=E_f(-H_Y)$, it suffices to show 
\[p_{Y*}(j_*\cE|_W)\cong E_f(-3H_Y)[-1].\]

Let \(Q\subset \PP(V)\) be the quartic hypersurface in \cite[Theorem 2.18]{kuznetsov2003derived}. Then we have \(W=\PP_Y(E_f^\vee)\times_Q \PP_X(\mathcal E)\). We write $\chi\colon W\to Q$ for the natural projection. Let $q_X\coloneqq p_X\circ j$ and $q_Y\coloneqq p_Y\circ j$.

By \cite[Proposition 2.20(i)]{kuznetsov2003derived}, there is an exact sequence on
\(Y\times X\)
\begin{equation}\label{A.1}
0\to
E_f^\vee(-H_Y)\boxtimes\mathcal O_X
\to
\mathcal O_Y(-H_Y)\boxtimes\mathcal Q
\to
\mathcal O_Y\boxtimes\mathcal E^\vee
\to
j_*\chi^*\mathcal O_Q(1)
\to0.
\end{equation}

On \(\pi\colon \PP_X(\mathcal E)\to X\), the tautological sequence is
\[
0\to \mathcal O_{\PP_X(\mathcal E)}(-e)
\to \pi^*\mathcal E
\to \mathcal O_{\PP_X(\mathcal E)}(e-H)
\to0,
\]
where \(\mathcal O_{\PP_X(\mathcal E)}(e)=\phi^*\mathcal O_Q(1)\) and $\phi\colon \PP_X(\mathcal E)\to Q$ is the morphism in \cite[Theorem 2.18]{kuznetsov2003derived}. Pulling this sequence back
to \(W\), we obtain
\[
0\to
\chi^*\mathcal O_Q(-1)
\to
\mathcal E|_W
\to
\chi^*\mathcal O_Q(1)\otimes q_X^*\mathcal O_X(-H_X)
\to0.
\]

We first claim that
\begin{equation}\label{A.3}
q_{Y*}\chi^*\mathcal O_Q(-1)=0.
\end{equation}
Indeed, for every \(y\in Y\), the fiber \(W_y\) is a rational quartic
curve. Moreover, \(\chi^*\mathcal O_Q(1)|_{W_y}\) has degree \(1\) on
the distinguished component mapping to the line
\(\PP(\ker f(y))\subset Q\), and degree \(0\) on any contracted
component. Hence
\[
H^0(W_y,\chi^*\mathcal O_Q(-1)|_{W_y})=0.
\]
Since \(p_a(W_y)=0\) and $\deg\bigl(\chi^*\mathcal O_Q(-1)|_{W_y}\bigr)=-1,$ we also get $H^1(W_y,\chi^*\mathcal O_Q(-1)|_{W_y})=0,$ and \eqref{A.3} follows by cohomology and base change. Therefore, we get
\begin{equation}\label{A.2}
q_{Y*}(\cE|_W)\cong q_{Y*}\bigl(\chi^*\mathcal O_Q(1)\otimes q_X^*\mathcal O_X(-H)\bigr).
\end{equation}

By the projection formula for the closed embedding \(j\), we have 
\[
j_*\bigl(\chi^*\mathcal O_Q(1)\otimes q_X^*\mathcal O_X(-H)\bigr)
\cong
j_*\chi^*\mathcal O_Q(1)\otimes p_X^*\mathcal O_X(-H).
\]
Tensoring the resolution \eqref{A.1} by \(p_X^*\mathcal O_X(-H)\) and pushing it forward to \(Y\), we then obtain
\[
q_{Y*}\bigl(\chi^*\mathcal O_Q(1)\otimes q_X^*\mathcal O_X(-H)\bigr)
\cong
E_f^\vee(-H_Y)[-1].
\]
Combining this with \eqref{A.2} and \eqref{A.3}, we obtain
\[
q_{Y*}(\mathcal E|_W)
\cong
E_f^\vee(-H_Y)[-1].
\]
Finally, since \(E_f\) is a rank $2$ bundle with $\det E_f\simeq \mathcal O_Y(2H_Y)$, we have $E_f^\vee\simeq E_f(-2H_Y).$
Thus, the desired isomorphism
\[
p_{Y*}(j_*\mathcal E|_W)
\cong
E_f(-3H_Y)[-1]
\]
follows.
\end{proof}

\section{Numerical computations} \label{appA}

In this appendix, we collect some numerical computations used in the main content. Readers can skip this section safely.

\subsection{Wall-crossing computations for genus 7 and 9 cases}

In this subsection, we compute potential walls for $\cE$ and $\cE(-H)[1]$. The lemmas here are used in the proof of  Proposition \ref{induce_range}.

\begin{Lem} [{\cite[Proposition 3.2]{li2018stability}}] \label{SBGg=7}
Let $X\coloneqq X_{12}$ and $F\in \Db(X)$ be a $\sigma_{\alpha, \beta}$-semistable object for some $\beta$ and $\alpha>0$.

\begin{enumerate}
    \item If $|\mu_H(F)|\leq \frac{1}{2\sqrt{2}}$, then $\frac{H\ch_2(F)}{H^3\ch_0(F)}\leq 0$,
    
    \item If $\frac{1}{2\sqrt{2}} \leq |\mu_H(F)|\leq 1-\frac{1}{2\sqrt{2}}$, then 
    $\frac{H\ch_2(F)}{H^3\ch_0(F)}\leq \frac{1}{2}|\mu_H(F)|^2-\frac{1}{16}$,
   
    \item If $1-\frac{1}{2\sqrt{2}} \leq |\mu_H(F)|\leq 1+\frac{1}{2\sqrt{2}}$, then $\frac{H\ch_2(F)}{H^3\ch_0(F)}\leq |\mu_H(F)|-\frac{1}{2}$,
    
    \item If $1+\frac{1}{2\sqrt{2}} \leq |\mu_H(F)|\leq 2-\frac{1}{2\sqrt{2}}$, then 
    $\frac{H\ch_2(F)}{H^3\ch_0(F)}\leq \frac{1}{2}|\mu_H(F)|^2-\frac{1}{16}$.

\end{enumerate}

\end{Lem}

\begin{Lem} \label{g=7wall1}
Let $X\coloneqq X_{12}$ and $\beta=-\frac{5}{6}$ or $-\frac{71}{84}$. Let $E\in \Coh^{\beta}(X)$ be a $\sigma_{\alpha, \beta}$-semistable object for some $\alpha>0$ with $\ch_{\leq 2}(E)=\ch_{\leq 2}(\cE_7)$. Then $E$ is $\sigma_{\alpha, \beta}$-semistable for all $\alpha>0$.
\end{Lem}

\begin{proof}
We only do computations for $\beta=-\frac{5}{6}$ here. Since $-\frac{71}{84}$ is very close to $-\frac{70}{84}=-\frac{5}{6}$, the argument is almost the same.

We are going to show that there is no wall for $E$ when $\beta=-\frac{5}{6}$. As in Section  \ref{sec:tilt stability wall algorithm}, a wall would be given by a short exact sequence
\[0\to A\to E\to B\to 0\]
in $\Coh^{-\frac{5}{6}}(X)$ such that following conditions hold:

\begin{enumerate}
    \item $\mu_{\alpha, -\frac{5}{6}}(A)=\mu_{\alpha, -\frac{5}{6}}(E)=\mu_{\alpha, -\frac{5}{6}}(B)$;
    
    \item $\Delta_H(A)\geq 0$ and $\Delta_H(B)\geq 0$;
    
    \item $\Delta_H(A)\leq \Delta_H(E)$ and $\Delta_H(B)\leq \Delta_H(E)$;
    
    \item $\ch_1^{-\frac{5}{6}}(A)\geq 0$ and $\ch_1^{-\frac{5}{6}}(B)=\ch_1^{-\frac{5}{6}}(E)-\ch_1^{-\frac{5}{6}}(A)\geq 0$.
\end{enumerate}

Note that $\ch_1^{-\frac{5}{6}}(E)>0$, thus the inequalities in $(d)$ are actually strict. We can assume that $\ch_{\leq 2}(A)=(a,bH,cL)$ and so $\ch_{\leq 2}(B)=(5-a, (-2-b)H, -cL)$
for some $a,b,c\in \mathbb{Z}$. If we divide the discriminant $\Delta_H(-)$ by $(H^3)^2$, the conditions above can be rewritten as

\begin{enumerate}
    \item $\frac{-36a\alpha^2+25a+60b+6c}{12(5a+6b)}=\frac{5-180\alpha^2}{156}$;
    
    \item $b^2-\frac{ac}{6}\geq 0$ and $(-2-b)^2-\frac{(5-a)(-c)}{6}\geq 0$;
    
    \item $b^2-\frac{ac}{6}\leq 4$ and $(-2-b)^2-\frac{(5-a)(-c)}{6}\leq 4$;
    
    \item $\frac{13}{6}> b+\frac{5}{6}a> 0$.
\end{enumerate}

Since $\alpha^2>0$ and $6b+5a>0$ by $(d)$, $(a)$ implies
\begin{equation} \label{g=7a1}
    (50a+125b+13c)(2a+5b)<0.
\end{equation}

Now $(a)$ and $(d)$ imply the following four cases:

\begin{enumerate}[(i)]
    \item $a>5$, $-\frac{5}{6}a<b<\frac{13-5a}{6}, c>-\frac{25}{13}(2a+5b)$;
    
    \item $0<a\leq 5$, $-\frac{5}{6}a<b<-\frac{2}{5}a, c>-\frac{25}{13}(2a+5b)$;
    
    \item $0<a\leq 5$, $-\frac{2}{5}a<b<\frac{13-5a}{6}, c<-\frac{25}{13}(2a+5b)$;
    
    \item $a\leq 0$, $-\frac{5}{6}a<b<\frac{13-5a}{6}, c<-\frac{25}{13}(2a+5b)$.
\end{enumerate}

Combined with $(b)$ and $(d)$, each case $(i)$ to $(iv)$ gives the following:

\begin{enumerate}[(i)]
    \item $5<a<\frac{169}{5}$, $-\frac{10a}{13}<b\leq -\frac{2}{5}(a+\sqrt{(a-5)a})$, $-\frac{25}{13}(2a+5b)<c\leq \frac{6b^2}{a}$; or $5<a<\frac{169}{5}$, $-\frac{2}{5}(a+\sqrt{(a-5)a})<b< \frac{13-5a}{6}$, $-\frac{25}{13}(2a+5b)<c\leq \frac{6(b+2)^2}{a-5}$
    
    \item $0<a\leq 5$, $-\frac{10a}{13}<b<-\frac{2a}{5}$, $-\frac{25}{13}(2a+5b)<c\leq \frac{6b^2}{a}$;
    
    \item $0<a<5$, $-\frac{2a}{5}<b<-\frac{2(5a-12)}{13}$, $\frac{6(b+2)^2}{a-5}\leq c<-\frac{25}{13}(2a+5b)$;
    
    \item $-\frac{144}{5}<a\leq 0$, $-\frac{5a}{6}<b\leq \frac{2}{5}(-a+\sqrt{(a-5)a})$, $\frac{6b^2}{a}\leq c<-\frac{25}{13}(2a+5b)$; or $-\frac{144}{5}<a\leq 0$, $\frac{2}{5}(-a+\sqrt{(a-5)a})<b<-\frac{2(5a-12)}{13}$, $\frac{6(b+2)^2}{a-5}\leq c<-\frac{25}{13}(2a+5b)$.
\end{enumerate}

Now by a careful computation for each case $(i)$ to $(iv)$, we obtain all possible truncated Chern characters of $A$ and $B$:

\begin{enumerate}[(1)]
    \item $(-11, 10H, -54L)$ and $(16, -12H, 54L)$;
    
    \item $(-5, 5H, -29L)$ and $(10, -7H, 29L)$;
    
    \item $(-4, 4H, -24L)$ and $(9, -6H, 24L)$;
    
    \item $(-3, 3H, -18L)$ and $(8, -5H, 18L)$;
    
    \item $(-3, 4H, -27L)$ and $(8, -6H, 27L)$;
    
    \item $(-2, 2H, -12L)$ and $(7, -4H, 12L)$;
    
    \item $(-1, H,-6L)$ and $(6,-3H, 6L)$;
    
    \item $(-1,2H,-16L)$ and $(6,-4H, 16L)$;
    
    \item $(0, H,-10L)$ and $(5,-3H, 10L)$;
    
    \item $(1,0,-6L)$ and $(4,-2H, 6L)$;
    
    \item $(1,0,-5L)$ and $(4,-2H, 5L)$;
    
    \item $(1,0,-4L)$ and $(4,-2H, 4L)$;
    
    \item $(2,-H, 2L)$ and $(3,-H,-2L)$;
    
    \item $(2,-H, 3L)$ and $(3,-H, -3L)$;
    
    \item $(2,0,-8L)$ and $(3,-2H, 8L)$.
\end{enumerate}

Since $A$ and $B$ are both $\sigma_{\alpha, -\frac{5}{6}}$-semistable for some $\alpha>0$, the cases $(13)$ and $(14)$ are ruled out by using Lemma \ref{SBGg=7} on the first character. The other cases are ruled out by using Lemma \ref{SBGg=7} on the second character. This implies that there are no walls when $\beta=-\frac{5}{6}$ for $E$, and $E$ is $\sigma_{\alpha, -\frac{5}{6}}$-semistable for every $\alpha>0$.

When $\beta=-\frac{71}{84}$, the computation argument is similar, and the solutions for $\ch_{\leq 2}(A)$ and $\ch_{\leq 2}(B)$ are the same as those when $\beta=-\frac{5}{6}$. Thus from the same argument using Lemma \ref{SBGg=7}, there are no walls when $\beta=-\frac{71}{84}$ for $E$, and $E$ is $\sigma_{\alpha, -\frac{71}{84}}$-semistable for every $\alpha>0$.
\end{proof}

\begin{Lem} \label{g=7wall2}
Let $X\coloneqq X_{12}$ and $\beta=-\frac{5}{6}$ or $-\frac{71}{84}$. Let $E\in \Coh^{\beta}(X)$ be a $\sigma_{\alpha, \beta}$-semistable object for some $\alpha>0$ with $\ch_{\leq 2}(E)=\ch_{\leq 2}(\cE(-H)[1])$. Then $E$ is $\sigma_{\alpha, \beta}$-semistable for all $\alpha>0$.
\end{Lem}

\begin{proof}
We assume that there is a wall when $\beta=-\frac{5}{6}$ or $-\frac{71}{84}$ for $E$, and that it is given by $A\to E\to B$. Then a similar computation as in Lemma \ref{g=7wall1} shows that all possible truncated Chern characters of $A$ and $B$ are:

\begin{enumerate}[(1)]
    \item $(-6,7H, -49L)$ and $(1,0,-5L)$;
    
    \item $(-5, 6H, -43L)$ and $(0,H,-11L)$;
    
    \item $(-4, 5H, -37L)$ and $(-1, 2H,-17L)$;
    
    \item $(-3, 4H, -32L)$ and $(-2, 3H, -22L)$;
    
    \item $(-3, 4H, -31L)$ and $(-2, 3H, -23L)$.
\end{enumerate}

Now using Lemma \ref{SBGg=7} on the first character in each case, all of the cases $(1)$ to $(5)$ are ruled out. This means that there are no walls for $E$ when $\beta=-\frac{5}{6}$ or $-\frac{71}{84}$, and hence $E$ is $\sigma_{\alpha, \beta}$-semistable for every $\alpha>0$.
\end{proof}

\begin{Lem} [{\cite[Proposition 3.2]{li2018stability}}] \label{SBGg=9}
Let $X\coloneqq X_{16}$ and $F\in \Db(X)$ be a $\sigma_{\alpha, \beta}$-semistable object for some $\beta$ and $\alpha>0$. If $\mu_H(F)=-\frac{1}{2}$, then $\frac{H\ch_2(F)}{H^3\ch_0(F)}\leq \frac{5}{64}$.

\end{Lem}

\begin{Lem} \label{g=9wall1}
Let $X\coloneqq X_{16}$ and $\beta=-\frac{3}{4}$ or $-\frac{31}{40}$. Let $E\in \Coh^{\beta}(X)$ be a $\sigma_{\alpha, \beta}$-semistable object for some $\alpha>0$ with $\ch_{\leq 2}(E)=\ch_{\leq 2}(\cE)$. Then $E$ is $\sigma_{\alpha, \beta}$-semistable for all $\alpha>0$.
\end{Lem}

\begin{proof}
We assume that there is a wall when $\beta=-\frac{3}{4}$ or $-\frac{31}{40}$ for $E$, and that it is given by $A\to E\to B$. Then a similar computation as in Lemma \ref{g=7wall1} shows that all of the possible truncated Chern characters of $A$ and $B$ are:

\begin{enumerate}[(1)]
    \item $(-1,H, -8L)$ and $(4,-2H,8L)$;
    
    \item $(1, 0, -4L)$ and $(2,-H,4L)$.
\end{enumerate}

Now using Lemma \ref{SBGg=9} on the second character in each case, both cases $(1)$ and $(2)$ are ruled out. This means that there are no walls for $E$ on $\beta=-\frac{3}{4}$ or $-\frac{31}{40}$, and hence $E$ is $\sigma_{\alpha, \beta}$-semistable for every $\alpha>0$.
\end{proof}

\begin{Lem} \label{g=9wall2}
Let $X\coloneqq X_{16}$ and $\beta=-\frac{3}{4}$ or $-\frac{31}{40}$. Let $E\in \Coh^{\beta}(X)$ be a $\sigma_{\alpha, \beta}$-semistable object for some $\alpha>0$ with $\ch_{\leq 2}(E)=\ch_{\leq 2}(\cE_9(-H)[1])$. Then $E$ is $\sigma_{\alpha, \beta}$-semistable for all $\alpha>0$.
\end{Lem}

\begin{proof}
We assume that there is a wall when $\beta=-\frac{3}{4}$ or $-\frac{31}{40}$ for $E$, and that it is given by $A\to E\to B$. Then a similar computation as in Lemma \ref{g=7wall1} shows that there are no such truncated Chern characters of $A$ and $B$. This means there are no walls for $E$ when $\beta=-\frac{3}{4}$ or $-\frac{31}{40}$, and hence $E$ is $\sigma_{\alpha, \beta}$-semistable for every $\alpha>0$.
\end{proof}

\subsection{Inequalities in Proposition \ref{prop:skyscraper projection stability}} \label{5.9computation}

In this subsection, we compute inequalities used in the proof of Proposition \ref{prop:skyscraper projection stability}.

Recall that
\begin{itemize}
    \item $g=6$: $(\alpha_0, \beta_0)=(\frac{1}{20}, -\frac{9}{10})$,
    
    \item $g=7$: $(\alpha_0, \beta_0)=(\frac{1}{12}, -\frac{5}{6})$,
    
    \item $g=8$: $(\alpha_0, \beta_0)=(\frac{1}{25}, -\frac{22}{25})$,
    
    \item $g=9$: $(\alpha_0, \beta_0)=(\frac{1}{8}, -\frac{3}{4})$,
    
    \item $g=10$: $(\alpha_0, \beta_0)=(\frac{1}{25}, -\frac{22}{25})$,
    
    \item $g=12$: $(\alpha_0, \beta_0)=(\frac{1}{25}, -\frac{21}{25})$.
\end{itemize}

\begin{Lem}
Let $X\coloneqq X_{2g-2}$. If we assume that $[A]=av+bw$ and $[B]=cv+dw$ for $a,b,c,d\in \mathbb{Z}$, then the solutions of inequalities

\begin{enumerate}[(1)]
    \item $[A]+[B]=[i^*\oh_x[-1]]$;
    
    \item $\Im Z^0_{\alpha_0, \beta_0}(A)\cdot \Im Z^0_{\alpha_0, \beta_0}(i^*\oh_x[-1])\geq 0$ and $\Im Z^0_{\alpha_0, \beta_0}(B)\cdot \Im Z^0_{\alpha_0, \beta_0}(i^*\oh_x[-1])\geq 0$;
    
    \item $\mu^0_{\alpha_0, \beta_0}(A)>\mu^0_{\alpha_0, \beta_0}(B)$;
    
    \item $1-\chi(A,A)+1-\chi(B,B)\leq \ext^1(i^*\oh_x,i^*\oh_x)$
\end{enumerate}

are listed below:

\begin{enumerate}[(i)]
    \item $g=6$ and ordinary:  $(a,b,c,d)=(-2,1,-3,2)$;
    
    \item $g=6$ and special:  $(a,b,c,d)=(-2,1,-3,2)$ or $(a,b,c,d)=(-4,2,-1,1)$;
    
    \item $g=7$: 
    there are no solutions;

    \item $g=8$: $(a,b,c,d)=(-2,1,-5,3)$ or $(a,b,c,d)=(-4,2,-3,2)$;
    
    \item $g=9$: 
    there are no solutions;
    
    \item $g=10$: 
    there are no solutions;
    
    \item $g=12$:
    there are no solutions.
    
\end{enumerate}
\end{Lem}

\begin{proof}
Note that for genus $g= 6, 8$ and $12$, there are only finitely many $(a,b,c,d)\in \mathbb{Z}^{\oplus 4}$ satisfying the conditions $(1)$, $(2)$, and $(4)$. Thus from a simple computation we obtain solutions to $(i), (ii), (iv)$ and $(vii)$.

\begin{itemize}

    \item $g=7$: $(1)$ and $(4)$ give

\begin{equation} \label{71}
    \frac{1}{6}(-6a-\sqrt{33}-6)\leq b\leq \frac{1}{6}(-6a+\sqrt{33}-6).
\end{equation}

From $(2)$ we obtain
\begin{equation}\label{72}
    \frac{1}{12}>\frac{13a+16b}{48}\geq 0.
\end{equation}
The combination of \eqref{71} and \eqref{72} implies $-11\leq a\leq -1$. Then it is not hard to check that the only possible solution $(a,b)$ for \eqref{71} and \eqref{72} is $(a,b)=(-6, 5)$. But this contradicts $(3)$.

    \item $g=9$:  $(1)$ and $(4)$ give

\begin{equation} \label{91}
    \frac{1}{4}(-2a-\sqrt{3}-2)\leq b\leq \frac{1}{4}(-2a+\sqrt{3}-2).
\end{equation}

From $(2)$ we obtain
\begin{equation}\label{92}
    \frac{1}{16}>\frac{11a+32b}{128}\geq 0.
\end{equation}
The combination of \eqref{91} and \eqref{92} implies $-7\leq a\leq -1$. Then it is not hard to check that there are no integer solutions $(a,b)$ for \eqref{91} and \eqref{92}.

    \item $g=10$: $(1)$ and $(4)$ give

\begin{equation} \label{101}
    \frac{1}{6}(-4a-\sqrt{7}-3)\leq b\leq \frac{1}{6}(-4a+\sqrt{7}-3).
\end{equation}

From $(2)$ we obtain
\begin{equation}\label{102}
    \frac{167}{1875}>\frac{611a+1200b}{5625}\geq 0.
\end{equation}
The combination of \eqref{101} and \eqref{102} implies $-8\leq a\leq -1$. Then it is not hard to check that there are no integer solutions $(a,b)$ for \eqref{101} and \eqref{102}.
\end{itemize}
\end{proof}

\subsection{Inequalities in Proposition \ref{gneq6cone}}

In this subsection, we compute the inequalities used in the proof of Proposition \ref{gneq6cone}.

\begin{Lem} \label{inequality}

Let the notation and assumptions be as in Case 1.2 of the proof of Proposition \ref{gneq6cone}. If we assume $[B]=av+bw+c[\cE]$, then the inequalities

\begin{itemize}
    \item $\Im(Z^0_{\alpha_g, \beta_g}(A))\geq 0, \Im(Z^0_{\alpha_g, \beta_g}(B))> 0$,
    
    \item $\Im(Z^0_{\alpha_g, \beta_g}(\ker(\lambda)))>0$,
    
    \item $\mu^0_{\alpha_g, \beta_g}(\ker(\lambda))<\mu^0_{\alpha_g, \beta_g}(\cE)$,

    \item $\mu^0_{\alpha_g, \beta_g}(C)>\mu^0_{\alpha_g, \beta_g}(B)> \mu^0_{\alpha_g, \beta_g}(\cE)$,
    
    \item $c>0,a<0$.

\end{itemize}

imply that $\frac{b}{a}\geq \mu_H(\cE)$ for $g\neq 7$ and $\frac{b}{2a}\geq \mu_H(\cE)$ when $g=7$.

\end{Lem}

\begin{proof}

We assume that $\frac{b}{a}< \mu_H(\cE)$ for $g\neq 7$ and $\frac{b}{2a}< \mu_H(\cE)$ when $g=7$,  and we will show that there are no such integers $a,b,c \in \mathbb{Z}$.
Recall that $(\alpha_7, \beta_7)=(\frac{\sqrt{71}}{84},-\frac{71}{84})$, $(\alpha_8, \beta_8)=(\frac{2\sqrt{79}}{875},-\frac{122}{125})$, $(\alpha_9, \beta_9)=(\frac{\sqrt{31}}{40},-\frac{31}{40})$, $(\alpha_{10},\beta_{10})=(\frac{\sqrt{\frac{5}{3}}}{33},-\frac{10}{11})$, and  $(\alpha_{12}, \beta_{12})=(\frac{1}{22},-\frac{19}{22})$.

First we assume that $g$ is even. 
In this case $\mu_H(\cE)=-\frac{1}{2}$ and  $\ch_{\leq 2}(B)=(a+2c,(b-c)H,(\frac{g-4}{2}c-\frac{g}{2}a-\frac{3g-6}{2}b)L)$. We have:

\begin{itemize}
   \item $0< \frac{(g-4)c-ga-(3g-6)b}{4g-4} -\beta_g(b-c) +\frac{a+2c}{2}(\beta^2_g-\alpha^2_g)\leq  \frac{\beta_g^2-\alpha^2_g}{2}$,
    
    \item $0> \frac{-ga-(3g-6)b}{4g-4} -\beta_gb +\frac{a}{2}(\beta^2_g-\alpha^2_g)$,
    
    \item $\frac{\beta_g a-b}{\frac{-ga-(3g-6)b}{4g-4} -\beta_gb +\frac{a}{2}(\beta^2_g-\alpha^2_g)}<\frac{1+2\beta_g}{\frac{g-4}{4g-4}+\beta_g+\beta^2_g-\alpha^2_g}$,

    \item $ \frac{2\beta_g}{\beta_g^2-\alpha_g^2}>\frac{\beta_g (a+2c)-(b-c)}{\frac{(g-4)c-ga-(3g-6)b}{4g-4} -\beta_g(b-c) +\frac{a+2c}{2}(\beta^2_g-\alpha^2_g)}>\frac{1+2\beta_g}{\frac{g-4}{4g-4}+\beta_g+\beta^2_g-\alpha^2_g}$,
    
    \item $a<0$,

    \item $-\frac{1}{2}a<b$.

\end{itemize}

When $g=8$, we have:

\begin{enumerate}
    \item $0< \frac{2332a+4081b+1458c}{12250}\leq \frac{2916}{6125}$,
    
    \item $b<-\frac{4}{7}a$,
    
    \item $-\frac{99632}{303389}a<b<-\frac{4}{7}a$,
    
    \item $-\frac{2989}{1458}> \frac{\frac{-122(a+2c)}{125}-(b-c)}{\frac{2332a+4081b+1458c}{12250}} >-\frac{5831}{729}$,
    
    \item $c>0, a<0$,
    
     \item $-\frac{1}{2}a<b$.

\end{enumerate}

Now $(b), (c),(e)$ and $(f)$ imply

\begin{equation} \label{g=8bcf}
    -\frac{1}{2}a<b<-\frac{4}{7}a.
\end{equation}

Also $(a)$ is equivalent to 

\begin{equation} \label{g=8a}
    -\frac{583(4a+7b)}{1458}<c\leq -\frac{583(4a+7b)-5832}{1458}.
\end{equation}

Thus \eqref{g=8bcf}, \eqref{g=8a}  and $(d)$ imply
\begin{align} \label{g=8fin1}
    &-\frac{110773008}{11281063}<a\leq -\frac{66479}{9558}, -\frac{a}{2}<b<-\frac{4(2654501a-13846626)}{32517071}, \\ &-\frac{213500a+115559b}{256419}<c\leq -\frac{583(4a+7b)-5832}{1458} \nonumber
\end{align}
or
\begin{equation} \label{g=8fin2}
    -\frac{66479}{9558}<a<0, -\frac{a}{2}<b<-\frac{4}{7}a, -\frac{213500a+115559b}{258066}<c\leq -\frac{583(4a+7b)-5832}{1458}.
\end{equation}
Thus we have $-9\leq a<0$, and it is not hard to see that there are no such integers $a,b,c$ satisfying either \eqref{g=8fin1} or \eqref{g=8fin2}.

When $g=10$ or $12$, the computation is similar to the $g=8$ case, so we omit the details.

Now we assume that $g=7$, thus $\ch_{\leq 2}(B)=(2a+5c,(b-2c)H,(-5a-6b)L)$. In this case we have:

\begin{enumerate}
    \item $0<  \frac{290a+348b+71c}{1008}  \leq \frac{355}{1008}$,
    
    \item $b<-\frac{5}{6}a$,
    
    \item $-\frac{3679}{4926}a<b<-\frac{5}{6}a$,
    
    \item $-\frac{2244}{71}<\frac{-\frac{71}{84}(2a+5c)-(b-2c)}{\frac{290a+348b+71c}{1008} }<-\frac{12}{5}$,
    
    \item  $c>0, a<0$,
    
    \item $-\frac{4}{5}a<b$.
    
\end{enumerate}

Now $(b),(c),(e)$ and $(f)$ imply

\begin{equation} \label{g=7bcef}
    -\frac{4}{5}a<b<-\frac{5}{6}a.
\end{equation}

Also $(a)$ is equivalent to

\begin{equation} \label{g=7a}
    -\frac{58(5a+6b)}{71}<c\leq -\frac{58(5a+6b)-335}{71}.
\end{equation}

Thus \eqref{g=7bcef}, \eqref{g=7a} and $(d)$ imply

\begin{equation} \label{g=7fin1}
    -\frac{804}{71}<a<0, -\frac{4}{5}a<b<-\frac{5}{6}a, \frac{-35a-6b}{72}<c\leq -\frac{58(5a+6b)-335}{71}
\end{equation}

or 

\begin{equation} \label{g=7fin2}
    -\frac{24120}{1309}<a\leq -\frac{804}{71}, -\frac{4}{5}a<b<\frac{4824-3679a}{4926}, \frac{-35a-6b}{72}<c\leq -\frac{58(5a+6b)-335}{71}.
\end{equation}

It is not hard to see that the only possible solution of \eqref{g=7fin1} and \eqref{g=7fin2} is $(a,b,c)=(-11,9,5)$, i.e. $\ch_{\leq 2}(B)=(3,-H,L)$. But since $B$ is $\sigma^0_{\alpha_g, \beta_g}$-semistable, this contradicts Lemma \ref{SBGg=7}.

Finally we assume that $g=9$. Then $\ch_{\leq 2}(B)=(a+3c,(b-c)H,(-3a-8b)L)$. In this case we have:

\begin{enumerate}

\item $0< \frac{33a+88b+31c}{320}\leq \frac{93}{320}$,

\item $b<-\frac{3}{8}a$,

\item $-\frac{509}{3703}a<b<-\frac{3}{8}a$,

\item $-\frac{8}{3}> \frac{-\frac{31}{40}(a+3c)-(b-c)}{\frac{33a+88b+31c}{320}} >-\frac{424}{31}$,

\item $c>0, a<0$,

\item $-\frac{1}{3}a<b$.

\end{enumerate}

Now $(b),(c),(e)$ and $(f)$ imply

\begin{equation} \label{g=9bcef}
    -\frac{1}{3}a<b<-\frac{3}{8}a.
\end{equation}

Also $(a)$ is equivalent to

\begin{equation} \label{g=9a}
    -\frac{11(3a+8b)}{31}<c\leq -\frac{11(3a+8b)-93}{31}.
\end{equation}

Thus \eqref{g=9bcef}, \eqref{g=9a} and $(d)$ imply

\begin{equation} \label{g=9fin1}
    -8<a<0, -\frac{1}{3}a<b<-\frac{3}{8}a, \frac{-15a-8b}{32}<c\leq \frac{-33a-88b+93}{31}
\end{equation}

or  

\begin{equation} \label{g=9fin2}
    -\frac{2976}{265}<a\leq -8, -\frac{1}{3}a<b<\frac{992-197a}{856}, \frac{-15a-8b}{32}<c\leq \frac{-33a-88b+93}{31}
\end{equation}

Then it is not hard to see that there are no integers $a,b,c$ satisfying either \eqref{g=9fin1} or \eqref{g=9fin2}.
\end{proof}

\begin{Lem} \label{g=6inequality}

Let the notation and assumptions be as in Case 2 of the proof of Proposition \ref{gneq6cone}. If we assume $[B]=av+bw+c[\cE]$, then the inequalities

\begin{itemize}
    \item $\Im(Z^0_{\alpha_g, \beta_g}(A))\geq 0, \Im(Z^0_{\alpha_g, \beta_g}(B))> 0$,
    
    \item $\Im(Z^0_{\alpha_g, \beta_g}(i^*(A)))\geq 0$,

    \item $\mu^0_{\alpha_g, \beta_g}(C)>\mu^0_{\alpha_g, \beta_g}(B)\geq  \mu^0_{\alpha_g, \beta_g}(\cE)$,
    
    \item $c>0, a<0$.

\end{itemize}

imply that $\frac{b}{a}\geq \mu_H(\cE)$, or $(a,b)=(-1,1),(-3,2)$.

\end{Lem}

\begin{proof}
We assume that $\frac{b}{a}<\mu_H(\cE)$, i.e. $-\frac{1}{2}a<b$. In this case we have $g=6$ and $\ch_{\leq 2}(B)=(a+2c, (b-c)H, (c-3a-6b)L)$. We have:

\begin{enumerate}
    \item $0< \frac{83a+240b+6c}{800} \leq \frac{323}{800}$,
    
    \item $b\leq \frac{305-83a}{240}$,
    
    \item $-\frac{720}{323}> \frac{-\frac{9}{10}(a+2c)-(b-c)}{\frac{83a+240b+6c}{800}} \geq -\frac{320}{3}$,
    
    \item $c>0, a<0$,
    
    \item $-\frac{1}{2}a<b$.
\end{enumerate}

 Now $(a)$, $(b)$, $(c)$, $(d)$ and $(e)$ imply
 
  \begin{equation}
     -5<a<0, -\frac{1}{2}a<b\leq \frac{305-83a}{240}, -\frac{216a+107b}{253}<c\leq -\frac{83a+240b-323}{6}
 \end{equation}
 
 or
 
 \begin{equation}
     -\frac{253}{32}<a\leq -5, -\frac{1}{2}a<b<\frac{253-61a}{186}, -\frac{216a+107b}{253}<c\leq -\frac{83a+240b-323}{6}.
 \end{equation}

It is not hard to see that the only possible values of $a,b \in \mathbb{Z}$ are $(a,b)=(-1,1)$ and $(a,b)=(-3,2)$.
\end{proof}

From a similar computation as in Lemma \ref{g=7wall1}, we have the following lemma:

\begin{Lem} \label{wall_3v-2w}
Let $X\coloneqq X_{10}$. Then there are no walls for the class $3v-2w$ on the line  $\beta=\beta_6=-\frac{9}{10}$ with respect to $\sigma_{\alpha,-\frac{9}{10}}$.
\end{Lem}

\end{appendix}

\bibliographystyle{alpha}
{\small{\bibliography{mybib}}}

@article {bayer2016space,
    AUTHOR = {Bayer, Arend and Macr\`i, Emanuele and Stellari, Paolo},
     TITLE = {The space of stability conditions on abelian threefolds, and
              on some {C}alabi--{Y}au threefolds},
   JOURNAL = {Invent. Math.},
  FJOURNAL = {Inventiones Mathematicae},
    VOLUME = {206},
      YEAR = {2016},
    NUMBER = {3},
     PAGES = {869--933},
      ISSN = {0020-9910,1432-1297},
   MRCLASS = {14F05 (14J32 14K05 18E30)},
  MRNUMBER = {3573975},
MRREVIEWER = {Colin\ Diemer},
       DOI = {10.1007/s00222-016-0665-5},
       URL = {https://doi.org/10.1007/s00222-016-0665-5},
}

@article {mukai95,
    AUTHOR = {Mukai, Shigeru},
     TITLE = {Curves and symmetric spaces. {I}},
   JOURNAL = {Amer. J. Math.},
  FJOURNAL = {American Journal of Mathematics},
    VOLUME = {117},
      YEAR = {1995},
    NUMBER = {6},
     PAGES = {1627--1644},
      ISSN = {0002-9327,1080-6377},
   MRCLASS = {14H45 (14J45 14M17)},
  MRNUMBER = {1363081},
MRREVIEWER = {Raquel\ Mallavibarrena},
       DOI = {10.2307/2375032},
       URL = {https://doi.org/10.2307/2375032},
}

@incollection {mukai:nonabelian,
    AUTHOR = {Mukai, Shigeru},
     TITLE = {Non-abelian {B}rill--{N}oether theory and {F}ano 3-folds
              [translation of {S}ugaku {\bf 49} (1997), no. 1, 1--24;
              {MR}1478148 (99b:14012)]},
      NOTE = {Sugaku Expositions},
   JOURNAL = {Sugaku Expositions},
  FJOURNAL = {Sugaku Expositions},
    VOLUME = {14},
      YEAR = {2001},
    NUMBER = {2},
     PAGES = {125--153},
      ISSN = {0898-9583,2473-585X},
   MRCLASS = {14D20 (14H60 14J45)},
  MRNUMBER = {1857462},
}

@article {bayer2017stability,
    AUTHOR = {Bayer, Arend and Lahoz, Mart\'i{} and Macr\`i, Emanuele and
              Stellari, Paolo},
     TITLE = {Stability conditions on {K}uznetsov components},
      NOTE = {With an appendix by Bayer, Lahoz, Macr\`i, Stellari and X.
              Zhao},
   JOURNAL = {Ann. Sci. \'Ec. Norm. Sup\'er. (4)},
  FJOURNAL = {Annales Scientifiques de l'\'Ecole Normale Sup\'erieure.
              Quatri\`eme S\'erie},
    VOLUME = {56},
      YEAR = {2023},
    NUMBER = {2},
     PAGES = {517--570},
      ISSN = {0012-9593,1873-2151},
   MRCLASS = {14F08 (14C34 18G80)},
  MRNUMBER = {4598728},
       DOI = {10.24033/asens.2539},
       URL = {https://doi.org/10.24033/asens.2539},
}

@article{bayer2024mukai,
  title={{Mukai bundles on Fano threefolds}},
  author={Bayer, Arend and Kuznetsov, Alexander and Macr{\`\i}, Emanuele},
  journal={arXiv preprint, arXiv:2402.07154},
  year={2024}
}

@article {feyzbakhsh2023new,
    AUTHOR = {Feyzbakhsh, Soheyla and Liu, Zhiyu and Zhang, Shizhuo},
     TITLE = {New perspectives on categorical {T}orelli theorems for
              del~{P}ezzo threefolds},
   JOURNAL = {J. Math. Pures Appl. (9)},
  FJOURNAL = {Journal de Math\'ematiques Pures et Appliqu\'ees. Neuvi\`eme
              S\'erie},
    VOLUME = {191},
      YEAR = {2024},
     PAGES = {Paper No. 103627, 39},
      ISSN = {0021-7824,1776-3371},
   MRCLASS = {14D20 (14D23 14J45)},
  MRNUMBER = {4821725},
       DOI = {10.1016/j.matpur.2024.103627},
       URL = {https://doi.org/10.1016/j.matpur.2024.103627},
}

@article {li2018stability,
    AUTHOR = {Li, Chunyi},
     TITLE = {Stability conditions on {F}ano threefolds of {P}icard number
              1},
   JOURNAL = {J. Eur. Math. Soc. (JEMS)},
  FJOURNAL = {Journal of the European Mathematical Society (JEMS)},
    VOLUME = {21},
      YEAR = {2019},
    NUMBER = {3},
     PAGES = {709--726},
      ISSN = {1435-9855,1435-9863},
   MRCLASS = {14F05 (14J45)},
  MRNUMBER = {3908763},
MRREVIEWER = {Izzet\ Coskun},
       DOI = {10.4171/JEMS/848},
       URL = {https://doi.org/10.4171/JEMS/848},
}

@article {jacovskis2024categorical,
    AUTHOR = {Jacovskis, Augustinas and Lin, Xun and Liu, Zhiyu and Zhang,
              Shizhuo},
     TITLE = {Categorical {T}orelli theorems for {G}ushel--{M}ukai
              threefolds},
   JOURNAL = {J. Lond. Math. Soc. (2)},
  FJOURNAL = {Journal of the London Mathematical Society. Second Series},
    VOLUME = {109},
      YEAR = {2024},
    NUMBER = {3},
     PAGES = {Paper No. e12878, 52},
      ISSN = {0024-6107,1469-7750},
   MRCLASS = {14F08 (14C34 14D20 14D22 14J45)},
  MRNUMBER = {4709831},
MRREVIEWER = {Joan\ Pons-Llopis},
       DOI = {10.1112/jlms.12878},
       URL = {https://doi.org/10.1112/jlms.12878},
}

@article {liu2023autoeq,
    AUTHOR = {Liu, Ziqi},
     TITLE = {On {F}ourier--{M}ukai type autoequivalences of {K}uznetsov
              components of cubic threefolds},
   JOURNAL = {Math. Nachr.},
  FJOURNAL = {Mathematische Nachrichten},
    VOLUME = {297},
      YEAR = {2024},
    NUMBER = {5},
     PAGES = {1866--1878},
      ISSN = {0025-584X,1522-2616},
   MRCLASS = {14F08 (14J30)},
  MRNUMBER = {4755740},
       DOI = {10.1002/mana.202300237},
       URL = {https://doi.org/10.1002/mana.202300237},
}

@article {kuznetsov:base-change,
    AUTHOR = {Kuznetsov, Alexander},
     TITLE = {Base change for semiorthogonal decompositions},
   JOURNAL = {Compos. Math.},
  FJOURNAL = {Compositio Mathematica},
    VOLUME = {147},
      YEAR = {2011},
    NUMBER = {3},
     PAGES = {852--876},
      ISSN = {0010-437X,1570-5846},
   MRCLASS = {14F05 (14A22 18E30)},
  MRNUMBER = {2801403},
MRREVIEWER = {Daniele\ Faenzi},
       DOI = {10.1112/S0010437X10005166},
       URL = {https://doi.org/10.1112/S0010437X10005166},
}

@book {huybrechts:geometry-of-moduli-space-of-sheaves,
	AUTHOR = {Huybrechts, D. and Lehn, M.},
	TITLE = {The geometry of moduli spaces of sheaves},
	SERIES = {Cambridge Mathematical Library},
	EDITION = {Second},
	PUBLISHER = {Cambridge University Press, Cambridge},
	YEAR = {2010},
	PAGES = {xviii+325},
	ISBN = {978-0-521-13420-0},
	MRCLASS = {14D20 (14F05)},
	MRNUMBER = {2665168},
	DOI = {10.1017/CBO9780511711985},
	URL = {http://dx.doi.org/10.1017/CBO9780511711985},
}

@article {ogrady:double-EPW-period,
    AUTHOR = {O'Grady, Kieran G.},
     TITLE = {Dual double {EPW}-sextics and their periods},
   JOURNAL = {Pure Appl. Math. Q.},
  FJOURNAL = {Pure and Applied Mathematics Quarterly},
    VOLUME = {4},
      YEAR = {2008},
    NUMBER = {2},
     PAGES = {427--468},
      ISSN = {1558-8599,1558-8602},
   MRCLASS = {14J35 (14C30)},
  MRNUMBER = {2400882},
MRREVIEWER = {Marco\ Andreatta},
       DOI = {10.4310/PAMQ.2008.v4.n2.a6},
       URL = {https://doi.org/10.4310/PAMQ.2008.v4.n2.a6},
}

@article {li2022derived,
    AUTHOR = {Li, Chunyi and Pertusi, Laura and Zhao, Xiaolei},
     TITLE = {Derived categories of hearts on {K}uznetsov components},
   JOURNAL = {J. Lond. Math. Soc. (2)},
  FJOURNAL = {Journal of the London Mathematical Society. Second Series},
    VOLUME = {108},
      YEAR = {2023},
    NUMBER = {6},
     PAGES = {2146--2174},
      ISSN = {0024-6107,1469-7750},
   MRCLASS = {14F08 (14J30 14J35 14J45 18G80)},
  MRNUMBER = {4673426},
MRREVIEWER = {Timothy\ De Deyn},
       DOI = {10.1112/jlms.12804},
       URL = {https://doi.org/10.1112/jlms.12804},
}

@article {kuznetsov:fano-threefolds,
    AUTHOR = {Kuznetsov, A. G.},
     TITLE = {Derived categories of {F}ano threefolds},
   JOURNAL = {Tr. Mat. Inst. Steklova},
  FJOURNAL = {Trudy Matematicheskogo Instituta Imeni V. A. Steklova},
    VOLUME = {264},
      YEAR = {2009},
    NUMBER = {Mnogomernaya Algebraicheskaya Geometriya},
     PAGES = {116--128},
      ISSN = {0371-9685},
   MRCLASS = {14F05 (14J10 14J45)},
  MRNUMBER = {2590842},
MRREVIEWER = {Alexandr V. Pukhlikov},
       DOI = {10.1134/S0081543809010143},
       URL = {https://doi.org/10.1134/S0081543809010143},
}

@book {huyb-book-FM,
    AUTHOR = {Huybrechts, D.},
     TITLE = {Fourier--{M}ukai transforms in algebraic geometry},
    SERIES = {Oxford Mathematical Monographs},
 PUBLISHER = {The Clarendon Press, Oxford University Press, Oxford},
      YEAR = {2006},
     PAGES = {viii+307},
      ISBN = {978-0-19-929686-6; 0-19-929686-3},
   MRCLASS = {14F05 (14-02 18E30)},
  MRNUMBER = {2244106},
MRREVIEWER = {Bal\'{a}zs Szendr\H{o}i},
       DOI = {10.1093/acprof:oso/9780199296866.001.0001},
       URL = {https://doi.org/10.1093/acprof:oso/9780199296866.001.0001},
}

@article {pirozhkov2020admissible,
    AUTHOR = {Pirozhkov, Dmitrii},
     TITLE = {Admissible subcategories of del {P}ezzo surfaces},
   JOURNAL = {Adv. Math.},
  FJOURNAL = {Advances in Mathematics},
    VOLUME = {424},
      YEAR = {2023},
     PAGES = {Paper No. 109046, 62},
      ISSN = {0001-8708,1090-2082},
   MRCLASS = {14F08 (14J26 18G80)},
  MRNUMBER = {4581971},
MRREVIEWER = {Joan\ Pons-Llopis},
       DOI = {10.1016/j.aim.2023.109046},
       URL = {https://doi.org/10.1016/j.aim.2023.109046},
}

@article {brambilla2013rank,
    AUTHOR = {Brambilla, Maria Chiara and Faenzi, Daniele},
     TITLE = {Rank-two stable sheaves with odd determinant on {F}ano
              threefolds of genus nine},
   JOURNAL = {Math. Z.},
  FJOURNAL = {Mathematische Zeitschrift},
    VOLUME = {275},
      YEAR = {2013},
    NUMBER = {1-2},
     PAGES = {185--210},
      ISSN = {0025-5874,1432-1823},
   MRCLASS = {14F05 (14J30 14J45)},
  MRNUMBER = {3101804},
MRREVIEWER = {Francesco\ Malaspina},
       DOI = {10.1007/s00209-012-1131-8},
       URL = {https://doi.org/10.1007/s00209-012-1131-8},
}

@article {Pertusi2021serreinv,
    AUTHOR = {Pertusi, Laura and Robinett, Ethan},
     TITLE = {Stability conditions on {K}uznetsov components of
              {G}ushel--{M}ukai threefolds and {S}erre functor},
   JOURNAL = {Math. Nachr.},
  FJOURNAL = {Mathematische Nachrichten},
    VOLUME = {296},
      YEAR = {2023},
    NUMBER = {7},
     PAGES = {2975--3002},
      ISSN = {0025-584X,1522-2616},
   MRCLASS = {14F08 (14J30 14J45 18G80)},
  MRNUMBER = {4626869},
MRREVIEWER = {Dylan\ Spence},
       DOI = {10.1002/mana.202200010},
       URL = {https://doi.org/10.1002/mana.202200010},
}

@article {bridgeland,
    AUTHOR = {Bridgeland, Tom},
     TITLE = {Stability conditions on triangulated categories},
   JOURNAL = {Ann. of Math. (2)},
  FJOURNAL = {Annals of Mathematics. Second Series},
    VOLUME = {166},
      YEAR = {2007},
    NUMBER = {2},
     PAGES = {317--345},
      ISSN = {0003-486X,1939-8980},
   MRCLASS = {14F05 (18E30)},
  MRNUMBER = {2373143},
MRREVIEWER = {Leovigildo\ M.\ Alonso Tarrio},
       DOI = {10.4007/annals.2007.166.317},
       URL = {https://doi.org/10.4007/annals.2007.166.317},
}

@article{zhang2020bridgeland,
  title={{Bridgeland Moduli spaces for Gushel--Mukai threefolds and Kuznetsov's Fano threefold conjecture}},
  author={Zhang, Shizhuo},
  journal={arXiv preprint arXiv:2012.12193},
  year={2020}
}

@inproceedings {mukai:icm,
    AUTHOR = {Mukai, Shigeru},
     TITLE = {Vector bundles on a {$K3$} surface},
 BOOKTITLE = {Proceedings of the {I}nternational {C}ongress of
              {M}athematicians, {V}ol. {II} ({B}eijing, 2002)},
     PAGES = {495--502},
 PUBLISHER = {Higher Ed. Press, Beijing},
      YEAR = {2002},
      ISBN = {7-04-008690-5},
   MRCLASS = {14J28 (14H51 14J60)},
  MRNUMBER = {1957059},
MRREVIEWER = {Mark\ Gross},
}

@article {mukai02,
    AUTHOR = {Mukai, Shigeru},
     TITLE = {New developments in {F}ano manifold theory related to the
              vector bundle method and moduli problems},
   JOURNAL = {Sūgaku},
  FJOURNAL = {Mathematical Society of Japan. Sūgaku (Mathematics)},
    VOLUME = {47},
      YEAR = {1995},
    NUMBER = {2},
     PAGES = {125--144},
      ISSN = {0039-470X,1883-6127},
   MRCLASS = {14J45 (14J10 14J60)},
  MRNUMBER = {1364825},
MRREVIEWER = {Takao\ Fujita},
}

@article {bayer2020desingularization,
    AUTHOR = {Bayer, Arend and Beentjes, Sjoerd Viktor and Feyzbakhsh,
              Soheyla and Hein, Georg and Martinelli, Diletta and Rezaee,
              Fatemeh and Schmidt, Benjamin},
     TITLE = {The desingularization of the theta divisor of a cubic
              threefold as a moduli space},
   JOURNAL = {Geom. Topol.},
  FJOURNAL = {Geometry \& Topology},
    VOLUME = {28},
      YEAR = {2024},
    NUMBER = {1},
     PAGES = {127--160},
      ISSN = {1465-3060,1364-0380},
   MRCLASS = {14D20 (14C34 14F08 14J30 14J45)},
  MRNUMBER = {4711835},
MRREVIEWER = {Francesco\ Bottacin},
       DOI = {10.2140/gt.2024.28.127},
       URL = {https://doi.org/10.2140/gt.2024.28.127},
}

@article {macri2007stability,
    AUTHOR = {Macr\`i, Emanuele},
     TITLE = {Stability conditions on curves},
   JOURNAL = {Math. Res. Lett.},
  FJOURNAL = {Mathematical Research Letters},
    VOLUME = {14},
      YEAR = {2007},
    NUMBER = {4},
     PAGES = {657--672},
      ISSN = {1073-2780},
   MRCLASS = {18E30 (14F05 16G20)},
  MRNUMBER = {2335991},
       DOI = {10.4310/MRL.2007.v14.n4.a10},
       URL = {https://doi.org/10.4310/MRL.2007.v14.n4.a10},
}

@article {debarre2012period,
    AUTHOR = {Debarre, Olivier and Iliev, Atanas and Manivel, Laurent},
     TITLE = {On the period map for prime {F}ano threefolds of degree 10},
   JOURNAL = {J. Algebraic Geom.},
  FJOURNAL = {Journal of Algebraic Geometry},
    VOLUME = {21},
      YEAR = {2012},
    NUMBER = {1},
     PAGES = {21--59},
      ISSN = {1056-3911,1534-7486},
   MRCLASS = {14J45 (14J30 14K30 14M22)},
  MRNUMBER = {2846678},
MRREVIEWER = {Marco\ Andreatta},
       DOI = {10.1090/S1056-3911-2011-00594-8},
       URL = {https://doi.org/10.1090/S1056-3911-2011-00594-8},
}

@article {kuznetsov2018derived,
    AUTHOR = {Kuznetsov, Alexander and Perry, Alexander},
     TITLE = {Derived categories of {G}ushel--{M}ukai varieties},
   JOURNAL = {Compos. Math.},
  FJOURNAL = {Compositio Mathematica},
    VOLUME = {154},
      YEAR = {2018},
    NUMBER = {7},
     PAGES = {1362--1406},
      ISSN = {0010-437X,1570-5846},
   MRCLASS = {14F05 (14E08 14J45)},
  MRNUMBER = {3826460},
MRREVIEWER = {P.\ E.\ Newstead},
       DOI = {10.1112/s0010437x18007091},
       URL = {https://doi.org/10.1112/s0010437x18007091},
}

@book {shafa,
     TITLE = {Algebraic geometry. {V}},
    SERIES = {Encyclopaedia of Mathematical Sciences},
    VOLUME = {47},
    EDITOR = {Parshin, A. N. and Shafarevich, I. R.},
      NOTE = {Fano varieties,
              A translation of {\it Algebraic geometry. 5} (Russian), Ross.
              Akad. Nauk, Vseross. Inst. Nauchn. i Tekhn. Inform., Moscow},
 PUBLISHER = {Springer-Verlag, Berlin},
      YEAR = {1999},
     PAGES = {iv+247},
      ISBN = {3-540-61468-0},
   MRCLASS = {14J45},
  MRNUMBER = {1668575},
MRREVIEWER = {Takao\ Fujita},
}

@article {KPS,
    AUTHOR = {Kuznetsov, Alexander G. and Prokhorov, Yuri G. and Shramov,
              Constantin A.},
     TITLE = {Hilbert schemes of lines and conics and automorphism groups of
              {F}ano threefolds},
   JOURNAL = {Jpn. J. Math.},
  FJOURNAL = {Japanese Journal of Mathematics},
    VOLUME = {13},
      YEAR = {2018},
    NUMBER = {1},
     PAGES = {109--185},
      ISSN = {0289-2316,1861-3624},
   MRCLASS = {14J45 (14C05 14J30 14J50)},
  MRNUMBER = {3776469},
MRREVIEWER = {Alexandr\ V.\ Pukhlikov},
       DOI = {10.1007/s11537-017-1714-6},
       URL = {https://doi.org/10.1007/s11537-017-1714-6},
}

@article {PY20,
    AUTHOR = {Pertusi, Laura and Yang, Song},
     TITLE = {Some remarks on {F}ano threefolds of index two and stability
              conditions},
   JOURNAL = {Int. Math. Res. Not. IMRN},
  FJOURNAL = {International Mathematics Research Notices. IMRN},
      YEAR = {2022},
    NUMBER = {17},
     PAGES = {12940--12983},
      ISSN = {1073-7928,1687-0247},
   MRCLASS = {14J45},
  MRNUMBER = {4475280},
MRREVIEWER = {Sho\ Tanimoto},
       DOI = {10.1093/imrn/rnaa387},
       URL = {https://doi.org/10.1093/imrn/rnaa387},
}

@article {dimitrov2019bridgeland,
    AUTHOR = {Dimitrov, George and Katzarkov, Ludmil},
     TITLE = {Bridgeland stability conditions on wild {K}ronecker quivers},
   JOURNAL = {Adv. Math.},
  FJOURNAL = {Advances in Mathematics},
    VOLUME = {352},
      YEAR = {2019},
     PAGES = {27--55},
      ISSN = {0001-8708,1090-2082},
   MRCLASS = {14F05 (16G20)},
  MRNUMBER = {3959651},
MRREVIEWER = {Chrysostomos\ Psaroudakis},
       DOI = {10.1016/j.aim.2019.05.032},
       URL = {https://doi.org/10.1016/j.aim.2019.05.032},
}

@article {FeyzbakhshPertusi2021stab,
    AUTHOR = {Feyzbakhsh, Soheyla and Pertusi, Laura},
     TITLE = {Serre-invariant stability conditions and {U}lrich bundles on
              cubic threefolds},
   JOURNAL = {\'Epijournal G\'eom. Alg\'ebrique},
  FJOURNAL = {\'Epijournal de G\'eom\'etrie Alg\'ebrique. EPIGA},
    VOLUME = {7},
      YEAR = {2023},
     PAGES = {Art. 1, 32},
      ISSN = {2491-6765},
   MRCLASS = {14F08 (14J30 14J45 14J60 18G80)},
  MRNUMBER = {4545360},
MRREVIEWER = {Pieter\ Belmans},
       DOI = {10.46298/epiga.2022.9611},
       URL = {https://doi.org/10.46298/epiga.2022.9611},
}

@article {brambilla2014vector,
    AUTHOR = {Brambilla, Maria Chiara and Faenzi, Daniele},
     TITLE = {Vector bundles on {F}ano threefolds of genus 7 and
              {B}rill--{N}oether loci},
   JOURNAL = {Internat. J. Math.},
  FJOURNAL = {International Journal of Mathematics},
    VOLUME = {25},
      YEAR = {2014},
    NUMBER = {3},
     PAGES = {1450023, 59},
      ISSN = {0129-167X,1793-6519},
   MRCLASS = {14J60 (14D20 14F05 14H30 14J45)},
  MRNUMBER = {3189780},
MRREVIEWER = {Carlo\ Giovanni\ Madonna},
       DOI = {10.1142/S0129167X14500232},
       URL = {https://doi.org/10.1142/S0129167X14500232},
}

@article {brambilla2011moduli,
    AUTHOR = {Brambilla, Maria Chiara and Faenzi, Daniele},
     TITLE = {Moduli spaces of rank-2 {ACM} bundles on prime {F}ano
              threefolds},
   JOURNAL = {Michigan Math. J.},
  FJOURNAL = {Michigan Mathematical Journal},
    VOLUME = {60},
      YEAR = {2011},
    NUMBER = {1},
     PAGES = {113--148},
      ISSN = {0026-2285,1945-2365},
   MRCLASS = {14J60 (14C17 14J30 14J45)},
  MRNUMBER = {2785867},
MRREVIEWER = {Ernesto\ C.\ Mistretta},
       DOI = {10.1307/mmj/1301586307},
       URL = {https://doi.org/10.1307/mmj/1301586307},
}

@article {kuznetsov2003derived,
    AUTHOR = {Kuznetsov, A. G.},
     TITLE = {Derived category of a cubic threefold and the variety
              {$V_{14}$}},
   JOURNAL = {Tr. Mat. Inst. Steklova},
  FJOURNAL = {Trudy Matematicheskogo Instituta Imeni V. A. Steklova},
    VOLUME = {246},
      YEAR = {2004},
     PAGES = {183--207},
      ISSN = {0371-9685,3034-1809},
   MRCLASS = {14J30 (14J45 14N05)},
  MRNUMBER = {2101293},
MRREVIEWER = {Ivan\ Chel\cprime tsov},
}

@article {ku:fractional-cy,
    AUTHOR = {Kuznetsov, Alexander},
     TITLE = {Calabi--{Y}au and fractional {C}alabi--{Y}au categories},
   JOURNAL = {J. Reine Angew. Math.},
  FJOURNAL = {Journal f\"ur die Reine und Angewandte Mathematik. [Crelle's
              Journal]},
    VOLUME = {753},
      YEAR = {2019},
     PAGES = {239--267},
      ISSN = {0075-4102,1435-5345},
   MRCLASS = {14F08 (14A30 14J32 18G80)},
  MRNUMBER = {3987870},
MRREVIEWER = {Clemens\ Koppensteiner},
       DOI = {10.1515/crelle-2017-0004},
       URL = {https://doi.org/10.1515/crelle-2017-0004},
}

@incollection {kuznetsov2010derivedcubicfourfolds,
    AUTHOR = {Kuznetsov, Alexander},
     TITLE = {Derived categories of cubic fourfolds},
 BOOKTITLE = {Cohomological and geometric approaches to rationality
              problems},
    SERIES = {Progr. Math.},
    VOLUME = {282},
     PAGES = {219--243},
 PUBLISHER = {Birkh\"auser Boston, Boston, MA},
      YEAR = {2010},
      ISBN = {978-0-8176-4933-3},
   MRCLASS = {14F05 (14E05)},
  MRNUMBER = {2605171},
MRREVIEWER = {Paolo\ Stellari},
       DOI = {10.1007/978-0-8176-4934-0\_9},
       URL = {https://doi.org/10.1007/978-0-8176-4934-0_9},
}

@article {bernardara2012categorical,
    AUTHOR = {Bernardara, Marcello and Macr\`i, Emanuele and Mehrotra,
              Sukhendu and Stellari, Paolo},
     TITLE = {A categorical invariant for cubic threefolds},
   JOURNAL = {Adv. Math.},
  FJOURNAL = {Advances in Mathematics},
    VOLUME = {229},
      YEAR = {2012},
    NUMBER = {2},
     PAGES = {770--803},
      ISSN = {0001-8708,1090-2082},
   MRCLASS = {14F05 (14C34 14J30)},
  MRNUMBER = {2855078},
       DOI = {10.1016/j.aim.2011.10.007},
       URL = {https://doi.org/10.1016/j.aim.2011.10.007},
}

@article {bondal2001reconstruction,
    AUTHOR = {Bondal, Alexei and Orlov, Dmitri},
     TITLE = {Reconstruction of a variety from the derived category and
              groups of autoequivalences},
   JOURNAL = {Compositio Math.},
  FJOURNAL = {Compositio Mathematica},
    VOLUME = {125},
      YEAR = {2001},
    NUMBER = {3},
     PAGES = {327--344},
      ISSN = {0010-437X,1570-5846},
   MRCLASS = {18E30 (14F05)},
  MRNUMBER = {1818984},
MRREVIEWER = {Richard\ P.\ Thomas},
       DOI = {10.1023/A:1002470302976},
       URL = {https://doi.org/10.1023/A:1002470302976},
}

@article {pertusi2022categorical,
    AUTHOR = {Pertusi, Laura and Stellari, Paolo},
     TITLE = {Categorical {T}orelli theorems: results and open problems},
   JOURNAL = {Rend. Circ. Mat. Palermo (2)},
  FJOURNAL = {Rendiconti del Circolo Matematico di Palermo. Second Series},
    VOLUME = {72},
      YEAR = {2023},
    NUMBER = {5},
     PAGES = {2949--3011},
      ISSN = {0009-725X,1973-4409},
   MRCLASS = {14F08 (14C34 14J28 14J45 18G80)},
  MRNUMBER = {4609668},
MRREVIEWER = {Shintarou\ Yanagida},
       DOI = {10.1007/s12215-022-00796-x},
       URL = {https://doi.org/10.1007/s12215-022-00796-x},
}

@article {happel1996tilting,
    AUTHOR = {Happel, Dieter and Reiten, Idun and Smal\o, Sverre O.},
     TITLE = {Tilting in abelian categories and quasitilted algebras},
   JOURNAL = {Mem. Amer. Math. Soc.},
  FJOURNAL = {Memoirs of the American Mathematical Society},
    VOLUME = {120},
      YEAR = {1996},
    NUMBER = {575},
     PAGES = {viii+ 88},
      ISSN = {0065-9266,1947-6221},
   MRCLASS = {16D90 (16G10 18E10 18E40)},
  MRNUMBER = {1327209},
MRREVIEWER = {Jeremy\ Rickard},
       DOI = {10.1090/memo/0575},
       URL = {https://doi.org/10.1090/memo/0575},
}

@article {bayer2011bridgeland,
    AUTHOR = {Bayer, Arend and Macr\`i, Emanuele and Toda, Yukinobu},
     TITLE = {Bridgeland stability conditions on threefolds {I}:
              {B}ogomolov--{G}ieseker type inequalities},
   JOURNAL = {J. Algebraic Geom.},
  FJOURNAL = {Journal of Algebraic Geometry},
    VOLUME = {23},
      YEAR = {2014},
    NUMBER = {1},
     PAGES = {117--163},
      ISSN = {1056-3911,1534-7486},
   MRCLASS = {14F05 (14J30)},
  MRNUMBER = {3121850},
MRREVIEWER = {Adrian\ Langer},
       DOI = {10.1090/S1056-3911-2013-00617-7},
       URL = {https://doi.org/10.1090/S1056-3911-2013-00617-7},
}

@article {bm99,
    AUTHOR = {Bridgeland, Tom and Maciocia, Antony},
     TITLE = {Fourier--{M}ukai transforms for {$K3$} and elliptic fibrations},
   JOURNAL = {J. Algebraic Geom.},
  FJOURNAL = {Journal of Algebraic Geometry},
    VOLUME = {11},
      YEAR = {2002},
    NUMBER = {4},
     PAGES = {629--657},
      ISSN = {1056-3911,1534-7486},
   MRCLASS = {14D06 (14D20 14J28 14J32)},
  MRNUMBER = {1910263},
       DOI = {10.1090/S1056-3911-02-00317-X},
       URL = {https://doi.org/10.1090/S1056-3911-02-00317-X},
}

@article {kuznetsov2005derived,
    AUTHOR = {Kuznetsov, A. G.},
     TITLE = {Derived categories of the {F}ano threefolds {$V_{12}$}},
   JOURNAL = {Mat. Zametki},
  FJOURNAL = {Matematicheskie Zametki},
    VOLUME = {78},
      YEAR = {2005},
    NUMBER = {4},
     PAGES = {579--594},
      ISSN = {0025-567X,2305-2880},
   MRCLASS = {14F05 (14J45 18E30)},
  MRNUMBER = {2226730},
MRREVIEWER = {Adrian\ Langer},
       DOI = {10.1007/s11006-005-0152-6},
       URL = {https://doi.org/10.1007/s11006-005-0152-6},
}

@incollection {Kuznetsov_2016,
    AUTHOR = {Kuznetsov, Alexander},
     TITLE = {Derived categories view on rationality problems},
 BOOKTITLE = {Rationality problems in algebraic geometry},
    SERIES = {Lecture Notes in Math.},
    VOLUME = {2172},
     PAGES = {67--104},
 PUBLISHER = {Springer, Cham},
      YEAR = {2016},
      ISBN = {978-3-319-46208-0; 978-3-319-46209-7},
   MRCLASS = {14F05 (14E07 14E08 14J45)},
  MRNUMBER = {3618666},
MRREVIEWER = {Sofia\ Tirabassi},
}

@article {kuznetsov2015categorical,
    AUTHOR = {Kuznetsov, Alexander and Lunts, Valery A.},
     TITLE = {Categorical resolutions of irrational singularities},
   JOURNAL = {Int. Math. Res. Not. IMRN},
  FJOURNAL = {International Mathematics Research Notices. IMRN},
      YEAR = {2015},
    NUMBER = {13},
     PAGES = {4536--4625},
      ISSN = {1073-7928,1687-0247},
   MRCLASS = {14F05 (14B05 18E30)},
  MRNUMBER = {3439086},
MRREVIEWER = {Andreas\ Krug},
       DOI = {10.1093/imrn/rnu072},
       URL = {https://doi.org/10.1093/imrn/rnu072},
}

@article {bayer2022kuznetsov,
    AUTHOR = {Bayer, Arend and Perry, Alexander},
     TITLE = {Kuznetsov's {F}ano threefold conjecture via {K}3 categories
              and enhanced group actions},
   JOURNAL = {J. Reine Angew. Math.},
  FJOURNAL = {Journal f\"ur die Reine und Angewandte Mathematik. [Crelle's
              Journal]},
    VOLUME = {800},
      YEAR = {2023},
     PAGES = {107--153},
      ISSN = {0075-4102,1435-5345},
   MRCLASS = {14F08 (14J30 14J45)},
  MRNUMBER = {4609825},
MRREVIEWER = {Dylan\ Spence},
       DOI = {10.1515/crelle-2023-0021},
       URL = {https://doi.org/10.1515/crelle-2023-0021},
}

@article {bayer2021stability,
    AUTHOR = {Bayer, Arend and Lahoz, Mart\'i{} and Macr\`i, Emanuele and
              Nuer, Howard and Perry, Alexander and Stellari, Paolo},
     TITLE = {Stability conditions in families},
   JOURNAL = {Publ. Math. Inst. Hautes \'Etudes Sci.},
  FJOURNAL = {Publications Math\'ematiques. Institut de Hautes \'Etudes
              Scientifiques},
    VOLUME = {133},
      YEAR = {2021},
     PAGES = {157--325},
      ISSN = {0073-8301,1618-1913},
   MRCLASS = {14F08 (14J42)},
  MRNUMBER = {4292740},
MRREVIEWER = {Hao\ Max\ Sun},
       DOI = {10.1007/s10240-021-00124-6},
       URL = {https://doi.org/10.1007/s10240-021-00124-6},
}

@article {liu2021note,
    AUTHOR = {Liu, Zhiyu and Zhang, Shizhuo},
     TITLE = {A note on {B}ridgeland moduli spaces and moduli spaces of
              sheaves on {$X_{14}$} and {$Y_3$}},
   JOURNAL = {Math. Z.},
  FJOURNAL = {Mathematische Zeitschrift},
    VOLUME = {302},
      YEAR = {2022},
    NUMBER = {2},
     PAGES = {803--837},
      ISSN = {0025-5874,1432-1823},
   MRCLASS = {14F08 (14D20 14D23 14J45)},
  MRNUMBER = {4480212},
MRREVIEWER = {Emma\ Lepri},
       DOI = {10.1007/s00209-022-03074-9},
       URL = {https://doi.org/10.1007/s00209-022-03074-9},
}

\end{document}